\documentclass[11pt,a4paper]{article}
\usepackage{packages}

\begin{document}
\title{A Riemannian approach for PDE--constrained shape optimization over the diffeomorphism group using outer metrics}

\author[1,2]{Estefania Loayza-Romero}
\author[3]{Lidiya Pryymak}
\author[3]{Kathrin Welker}

\affil[1]{University of Strathclyde, Department of Mathematics and Statistics, 26 Richmond St. Glasgow, United Kingdom}
\affil[2]{MODEMAT Research Center in Mathematical Modeling and Optimization, Quito, Ecuador; \textsc{estefania.loayza-romero@strath.ac.uk}}
  
\affil[3]{Helmut-Schmidt-University / University of the Federal Armed Forces Hamburg, Faculty of Mechanics and Civil Engineering, Holstenhofweg~85, 22043 Hamburg, Germany; \textsc{pryymak@hsu-hh.de}, \textsc{welker@hsu-hh.de}}

\maketitle             
\begin{abstract}
In this paper, we study the use of outer metrics, in particular Sobolev-type metrics on the diffeomorphism group in the context of PDE-constrained shape optimization. 
Leveraging the structure of the diffeomorphism group, we analyze the connection between the push-forward of a function on the diffeomorphism group and the classical shape derivative. 
This connection and the representation of a shape through diffeomorphisms allow us to formulate PDE-constrained shape optimization problems as optimization problems over a Riemannian manifold. 
The different Sobolev metrics considered in this work are then used to formulate different variants of the Riemannian steepest descent method.
We consider in particular, two representative examples of PDE-constrained shape optimization. The first is inspired by electric impedance tomography, and the second is a two-dimensional bridge optimization problem.
These problems are discretized using finite elements and solved with the different variants of the method. 
For comparison reasons, we also solve the problem using the previously proposed Steklov-Poincaré metric. 
\end{abstract}
\\[0.5cm]

\maketitle

\section{Introduction}
\label{sec:Introduction}

In essence, shape optimization involves looking for geometries that enhance performance. The focus in this paper lies in applications tied to partial differential equations (PDE) due to their broad scope. For instance, one might seek the most rigid structure within a given volume~\cite{Allaire2014shape}, optimize the distribution between void and material to enhance an electric motor's performance~\cite{Gangl2015shape}, or identify inclusions in certain material using the principles of impedance tomography~\cite{Ito-Kunisch-Peichl}, among other possibilities.
A key challenge in shape optimization is the suitable choice of the underlying shape space. Typically, the shape space is selected to align well with the optimization problem to be solved and the solving algorithms. A prototypical algorithm to solve PDE-constrained shape optimization problems is the steepest descent method, for which a gradient needs to be computed. 
There are several ways to approach this, including the level-set method~\cite{Allaire2014shape} and the method of mappings~\cite{Haubner2021continous}. However, in this paper, we focus on the 
Riemannian version of the gradient descent method which can be found in ~\cite[Alg.~4.1]{Boumal2023intromanifolds}. Various variants of the method can be devised depending on the choice of metric and retraction.

The choice of inner product used to identify a descent direction from the shape derivative has received considerable attention. In the Hilbertian extension-regularization procedure surveyed in~\cite{Allaire2021shape}, descent directions are Riesz representatives with respect to an inner product on a linear space of vector fields over a fixed hold-all domain, and the influence of that choice on the behaviour of the resulting method is well documented, going back to the work of Schulz \cite{Schulz2014riemannian} and coauthors, e.g., \cite{Schulz2016computational,Welker2016efficient}. 
The same viewpoint underlies the line of work in~\cite{Hiptmair2015shape,Paganini2018higher} but two features distinguish that setting from ours: first, the metrics used there are of first order, i.e., the $H^1_0$ product on the hold-all domain in~\cite{Hiptmair2015shape} and in~\cite[Sec.~7.1]{Paganini2018higher}, and the Cauchy--Riemann-modified elasticity product in~\cite[Sec.~7.2]{Paganini2018higher}. 
It is also worth mentioning that the ``higher-order'' from~\cite{Paganini2018higher} concerns the polynomial degree of the discretization rather than the order of the operator defining the gradient. 
Second, and more substantially, the set of admissible transformations carries no evident manifold structure, right translation on the group is used in~\cite{Hiptmair2015shape}, explicitly through the identity $T_{V+W}=T_{W\circ T_V^{-1}}\circ T_V$ and implicitly through the composition update and the transported basis~\cite{Paganini2018higher}, but the inner product is a regularizing device rather than a Riemannian metric. 

In order to make a connection between classical shape calculus and the differential geometrical structure of the considered shape space, we analyze the relation between the push-forward of a smooth shape functional in the Riemannian framework and the actual shape derivative, as an Eulerian semi-derivative.  

In the last decade, one common shape space used in shape optimization is the space of simple closed one-dimensional infinitely smooth curves (cf.~\cite{Michor2004riemannian}) defined by
\begin{equation}
\label{eq:Be}
B_e:=B_e(S^1, \mathbb{R}^2):= \textup{Emb}(S^1, \mathbb{R}^2)/\textup{Diff}(S^1),
\end{equation}
where $S^1$ denotes the one-dimensional unit sphere, $\textup{Emb}(S^1, \mathbb{R}^2)$ denotes the set of all embeddings of $S^1$ into $\R^2$ and $\textup{Diff}(S^1)$ denotes the set of all diffeomorphisms of $S^1$. 
The use of this shape space in the context of PDE-constrained shape optimization was first investigated in~\cite{Schulz2014riemannian} and in a number of texts (cf., e.g., \cite{Geiersbach2021stochastic,Welker2016efficient,Schulz2016computational,Schulz2016efficient}), with particular emphasis on an inner metric, i.e., a metric induced by right-invariant metrics from $\textup{Emb}(S^1, \mathbb{R}^2)$. In \cite{Schulz2016efficient}, the so-called Steklov--Poincar\'e (SP) metric is established as an inner metric on the space $B_e$. 
\begin{align}
\label{eq:SPmetric}
     G^S \colon H^{\nicefrac12 }(u) \times H^{\nicefrac12 }(u) \to \R ,\, (v,w) \mapsto \int_u  v(S^{\textup{pr}})^{-1}w \, \textup{d} s.
\end{align}

Here $u\in B_e$ denotes a shape, that is, a simple closed curve in $\R^2$, and $H^{\nicefrac12}(u)$ is the trace space of scalar functions on $u$ with arc-length element $\dd s$. 
The arguments $v$ and $w$ are elements of that space, to be understood as normal deformations of the curve: a scalar function on $u$ prescribing, at each point, the displacement in the direction of the unit normal $\bm n$. 
Finally, $S^{\textup{pr}}$ denotes the projected Poincar\'e--Steklov operator, which maps a scalar force on $u$ to the normal trace on $u$ of the vector field it generates through an auxiliary elliptic problem on the hold-all domain, usually chosen as the Lam\'e system. As its inverse appears in~\eqref{eq:SPmetric}, a normal deformation of the curve is measured by the energy of the ambient field producing it. 
We refer to~\cite[Sec.~3]{Schulz2016efficient} for the precise definition, which we do not need in what follows.

Unfortunately, using inner metrics, we obtain shape gradients which only deform the shape itself and not the ambient space, for instance $\R^2$. 

\begin{figure}[h]
    \centering
    \includegraphics[width= 0.6\textwidth]{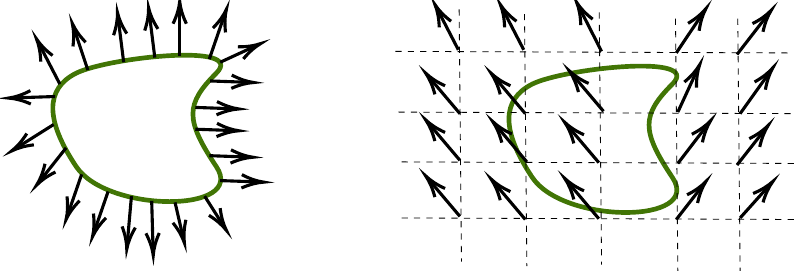}
    \caption{Tangent vectors on $B_e$, with inner metrics (left) and outer metrics where the vector field deforms the ambient space and induces a deformation of the shape (right).} 
    \label{fig:illustration}
\end{figure}

Inspired by the analysis proposed in~\cite{Bauer2011new}, where outer and inner metrics were compared in surface registration problems, 
we now propose a novel approach by considering outer metrics in the context of PDE-constrained shape optimization instead. An outer metric on $B_e$ is given by a right-invariant metric on the space 
\begin{equation}
    \label{Diffcn}
\Diffcn \coloneqq \{\varphi\in \smooth(\R^n, \R^n)|\, \varphi^{-1} \in \smooth(\R^n, \R^n), \,\operatorname{supp}(\varphi - \id) \text{ is compact}\}
\end{equation}
for $n=2$, inducing a metric on  $\text{Emb}(S^1, \mathbb{R}^2)$ by left action (cf. \cite{Bauer2013overview}). Hereby, $\id$ denotes the identity map and $\text{supp}$ the support of a function. For clarity, we note that, with a slight abuse of notation, the term `outer metric' will be used throughout this paper to refer to a metric on $\Diffc$. Using such outer metrics one expects to immediately obtain shape gradients which deform the ambient space $\mathbb{R}^2$, and thus inducing a deformation on the shape, as depicted in \cref{fig:illustration}.  
Representative metrics on $\Diffc$ are given by Sobolev-type metrics of order $s \geq 0$. 
As described in~\cite{Michor2007overview}, these are defined by
\begin{equation}
\label{eq:outerMetric}
H^s(\bm{V},\bm{W}):= \int_{\mathbb{R}^2} \langle L\bm{V},\bm{W} \rangle \ \textup{d}z, \qquad\text{with } L=(\textup{id}-A\Delta)^s \ \text{and } A>0
\end{equation}
for any two vector fields $\bm{V}= \bm{V}(x,y),\bm{W} = \bm{W}(x,y)$ in $\mathbb{R}^2$ with compact support, and $z = (x,y)^\top$. 
Hereby, $\Delta$ denotes the Laplacian operator. 
Additionally, the inverse of the operator $L$ is an integral operator with a kernel given by a classical Bessel function applied to the distance in $\R^2$ between two points on the curve, see \cite{Michor2007overview} for more details.
There are various important properties of this metric which we aim to exploit in this and later publications. 
For diffeomorphisms with Sobolev regularity $s\in \N$, the Riemannian manifold $\left(\textup{Diff}(\mathbb{R}^2),H^s\right)$ is geodesically and metrically complete for integer $s>2$; see, e.g.,~\cite{Bruveris2017completeness,Younes2010shapes}. Moreover, there exists an explicit expression for the geodesic equation associated to this metric (cf.~\cite{Michor2007overview}).

The contributions of this paper are the following.
First, we formulate PDE-constrained shape optimization on $(\Diffc,H^s)$ as optimization on an infinite-dimensional Lie group with a right-invariant Riemannian metric, and establish the connection between the classical shape derivative and the push-forward in this setting.
Second, we show how the $H^s$ gradient can actually be computed. Its defining problem is of order $2s$, and we reformulate it, inspired by the stable splitting of~\cite{Gallistl2017stable}, as a chain of $s$ coupled second-order problems. 
\Cref{prop:split-equivalence} characterises the field this chain produces, identifying its regularity, the boundary conditions it inherits, and its relation to the $H^s$ inner product; the reformulation makes $H^s$ gradients accessible to piecewise linear finite elements, the discretization for which discrete and continuous shape derivatives coincide~\cite{Berggren2010unified}. 
Third, we compare these gradients systematically on two PDE-constrained problems undergoing large deformations. 
In the interface identification problem the $H^s$ variants with $s\ge2$ keep the mesh-quality measure essentially constant along the optimization, where the Steklov--Poincar\'e and $H^1$ variants degrade it by several orders of magnitude. 
In the compliance problem, where remeshing is admitted and the attained objective is the basis of comparison, $H^3$ and $H^4$ variant reach lower objective values in fewer iterations and with fewer PDE solves, while $H^2$ and the SP variants reach the maximum number of iterations.

The paper is organized as follows.
In \cref{sec:RiemanStrucDiff}, we introduce our notation and the basics  on differential geometry regarding the two above-mentioned shape spaces necessary to follow the subsequent sections. \Cref{sec:ShapeOpt} focuses on shape optimization, summarizing the classical and Riemannian approach and presenting a general gradient descent algorithm. Beyond that, we introduce our new approach of shape optimization on the shape space $\Diffc$ with an appropriate Riemannian metric. 
This section is also devoted to present one of our main results, the connection of the classical shape derivative with the push-forward and
the examination of the existence of shape gradients for the considered metrics. 
Building upon the obtained results, in \cref{sec:modelProb}, we consider two examples of PDE-constrained optimization problems as optimization problems on the diffeomorphism group, namely,
a tracking-type shape optimization approach for electric impedance tomography and a compliance minimization problem for a two-dimensional bridge. 
Finally, in \cref{sec:NumericalEx}, we numerically solve both problems with the Riemannian steepest descent method while varying the Riemannian metric with which we compute the descent directions. 
\Cref{sec:conclusion} concludes the paper with a summary of our main results and discussion of future research lines.

\section{Riemannian structure  on $\Diffc$}
\label{sec:RiemanStrucDiff}

We will introduce basic notions on smooth manifolds in \cref{subsec:PrelDiffLie}. This will lead to an investigation of existing results on the diffeomorphism group---in particular, its tangent structure (\cref{subsec:PropDiffGroup}) and Riemannian structure when equipped with Sobolev-type metrics as in~\eqref{eq:outerMetric} (\cref{subsec:SobtypeMetric}). For further details on infinite-dimensional smooth manifolds, we refer the reader to the literature, e.g., \cite{Kriegl1997convenient,Kriegl1991aspects,Michor1980manifolds,Michor2020manifolds,Michor2013zoo}.

\subsection{Preliminaries on differential geometry}
\label{subsec:PrelDiffLie}

Spaces of interest in this paper are infinite-dimensional manifolds modeled on $C^{\infty}$-open subsets of convenient vector spaces $E$. 
As usual, for any infinite-dimensional smooth manifold $M$ of this kind, one can distinguish between the operational and the kinematic tangent space. 
We define the operational tangent bundle $DM$ as the set of bounded derivations of smooth functions.
One may, however, consider the kinematic tangent bundle  as the space
$T_pM=\{ c\colon \R \rightarrow M\mid c\text{ differentiable, } c(0) =p\}/\sim$
with the following equivalence relation: for two differentiable curves  $c,\tilde{c}\colon\R \rightarrow M$  with $c(0) = \tilde{c}(0) =p$, we consider $
c \sim \tilde{c} \Leftrightarrow \dfrac{d}{d t} \big|_{t=0}\phi_{\alpha}(c(t))  =\dfrac{d}{d t} \big|_{t=0} \phi_{\alpha}(\tilde{c}(t))
$ 
for all $\alpha$ with $p \in U_\alpha$, where $\{(U_\alpha, \phi_\alpha)\}_\alpha$ is an  atlas of $M$. 
Note that while $TM \subset DM$ as a subbundle (cf. \cite[Sec.~28.12]{Kriegl1997convenient}), in general, there exist operational tangent vectors which are not kinematic ones. 
However, given certain conditions on the underlying vector space (cf.~\cite[Sec.~28.7]{Kriegl1997convenient}), the spaces are identical. In particular, $TM=DM$, if $E$ is a reflexive, nuclear Fr\'echet or LF  space \cite{Kriegl1991aspects}.
The space of smooth kinematic vector fields is defined as
$\X(M):=\{\bm{V}:M \ra TM \ \vert \ \bm{V}(f) \in \smooth(M) \ \fa f \in \smooth(M)\}.$

A Riemannian metric on a smooth manifold $M$ is a positive definite family  $g=\{g_p\}_{p \in M}$ of inner products on $T_pM$, where $p \mapsto g_p$ is smooth for all $p\in M$. We call $(M,g)$ a Riemannian manifold. For finite-dimensional Riemannian manifolds, the map $TM \rightarrow TM^{*}, \ v \mapsto g_p(v,.)$ where $v \in T_pM$ is an isomorphism and therefore $g$ meets the requirements for the definition of a so called strong Riemannian metric. In infinite dimensions, however, this is not necessarily the case. Indeed, a Riemannian metric on an infinite-dimensional manifold can only be strong, if $M$ is a Hilbert manifold (cf. \cite[Proposition 4.5]{Schmeding2022introduction}).

For a $C^1$-curve $c\colon [a,b]\to M$, the length of $c$ is defined as $\textup{L}(c) = \int_a^b \left\|c'(t) \right\|_{g} \textup{d}t$, where $\left\| v \right\|_{g} \coloneqq \sqrt{g_p(v,v)} \ \forall v \in T_pM , p \in M$. 
The geodesic distance between two points $p, \ q \in M$ is defined as 
\[  
		\mathrm{d}(p,q) 
	:= \inf \{\textup{L}(c) |\,   c\colon [0,1]\to M  \text{ piecewise smooth curve} \text{ with } c(0)=p, \ c(1)=q\}.
\] 
On arbitrary infinite-dimensional Riemannian manifolds, there may exist points $p,q \in M$, such that $\mathrm{d}(p,q)=0$ and $p \neq q$. 
In particular, $\mathrm{d}$ may not define a distance function on $M$. 
If, however, we assume that $g$ is a strong Riemannian metric, then $\mathrm{d}$ defines an actual metric on $M$ (cf. \cite[Theorem 1.9.5]{Klingenberg1995riemannian}).

Curves that are critical points of the energy functional $\textup{E}(c)= \frac{1}{2} \int_a^b g_{c(t)}(c'(t), c'(t)) \ \textup{d}t$ are called geodesics. 
Equivalently, geodesics can be defined as solutions to $\nabla_{\gamma^{'}}\gamma^{'}=0$, where $\nabla$ denotes the Levi-Civita connection associated with $(M,g)$.
Note that under certain conditions there exists a neighborhood of the origin $V \subset TM$ such that the map $\ex\colon V \ra M, v \mapsto c_v(1)$ exists and is smooth, where $c_v\colon [0,1] \ra M$ is the unique geodesic for the initial condition $c'(0)=v$. 
The exponential map in $p \in M$ is then defined as the restriction of $\ex$ to the set $V \cap T_pM$ and denoted by $\ex_p$. However, this definition is not always possible, as the existence of geodesics is not guaranteed. Again, we obtain existence of the exponential map for strong Riemannian manifolds. Moreover, in that case there exists a $V \subset TM$, such that the exponential map is a local diffeomorphism (cf. \cite[Theorem 1.8.15]{Klingenberg1995riemannian}). Often, geodesics are described as locally shortest paths between points on a manifold. This description is not valid for arbitrary infinite-dimensional Riemannian manifolds.  If $(M,g)$ is a strong Riemannian manifold, geodesics are locally distance-minimizing paths, i.e., for a geodesic $\tilde{c}\colon [0,1] \ra \R$, where $\tilde{c}(0)=p$, there exists a $\tilde{t}>0$, such that it holds $\mathrm{d}(p,c(t))= \textup{L}(\tilde{c})$ for each $t \in [0, \tilde{t})$. 

For a smooth function $F\colon M \ra N$ between smooth manifolds $M$ and $N$, the differential or the push-forward of $F$ is defined as a function $dF_p\colon T_pM \ra T_{F(p)}N$ for each $p \in M$ such that 
\begin{equation}
\label{eq:pushForward}
dF_p(v)(f):=v(f \circ F) \qquad \forall v \in T_pM,
f \in C^{\infty}(N).
\end{equation}
Equivalently, $dF_p(c):=\frac{\mathrm{d}}{\mathrm{d}t}\vert_{t=0} F(c(t)) = (F \circ c)'(0),$ where $c \in T_pM$. 
Using the definition of the push-forward, the gradient $\nabla F \in \X(M)$ of a smooth function $F\colon M \ra \R$ on a Riemannian manifold $(M,g)$ is defined as the solution of 
\begin{equation}
\label{eq:gradient}
g_p(\nabla F(p), v) = dF_p(v) \ \ \forall p \in M \text{ and } \forall v \in T_pM.
\end{equation}

\subsection{Properties of the diffeomorphism group}
\label{subsec:PropDiffGroup}

We consider the diffeomorphism group of functions equal to the identity outside a compact set, i.e., the space $\Diffcn$ defined in \eqref{Diffcn}. 
This space is a smooth submanifold of $\smooth(\R^n, \R^n)$ as has been proved in \cite[Theorem 43.1]{Kriegl1997convenient}. 
Moreover, it holds that  $T \Diffcn \cong \Diffcn \times T_{\id} \Diffcn \cong \Diffcn \times \Xc(\R^n)$. At this point, we want to look further into the identification between $T_{\varphi} \Diffcn$ and $\Xc(\R^n)$ for some $\varphi \in \Diffc$. 
Hereby, $\Xc(\R^n) \owns \bm{V} \cong \gamma \in T_{\varphi}\Diffcn$ if and only if
\begin{equation}
\label{eq:tangentDiff}
F_t(x)  \coloneq\gamma(t,x)  \quad  \text{satisfies}\quad \dfrac{d}{dt}\Big|_{t=0} F_t(x)=\bm{V}(F_{0}(x)), \ F_0(x)=\varphi(x).
\end{equation}
In this paper, we focus on the case $n=2$. Consider now the space 
$\text{Diff}_c(\R^2, \innerCirc):= \{ \varphi \in \Diffc \vert \varphi(\innerCirc)=\innerCirc\} ,$
where $\innerCirc\coloneq \{(x,y) \in \R^2 \vert \ x^2+y^2=1\} \subset \R^2$ denotes the unit circle.  There is a well-known (cf. \cite{Michor2007overview}) connection to the shape space $B_e$ introduced in \eqref{eq:Be}:
\begin{equation}
\label{eq:diffDelta}
\Diffc/\text{Diff}_c(\R^2, \innerCirc) \cong B_e ,    
\end{equation}
where the isomorphism is given by $ \Diffc\ni \varphi\mapsto \varphi \circ i$ with $i\colon S^1 \ra \R^2,\, t \mapsto (\sin t, \cos t)$. Hereby, $S^1= \R / (2 \pi \Z) $.

\subsection{Sobolev-type metrics}
\label{subsec:SobtypeMetric}

We consider the manifold $\Diffc$ equipped with a Sobolev-type metric of order $s \in \N$ defined for two vector fields $\bm{V},\bm{W}\colon \R^2\to\R^2 $, $\bm{V}=\bm{V}(x,y), \bm{W}=\bm{W}(x,y)$ with $z=(x,y)^{\top}$  by~\eqref{eq:outerMetric}
(cf. \cite{Bauer2013overview,Bauer2020sobolev,Bruveris2017completeness}).
Equivalently, the metric can be expressed by 
\begin{equation}
    \label{eq:metric2}
    H^s(\bm{U},\bm{V})=\sum\limits_{\substack{i,j \geq 0\\i+j\leq s}}\frac{A^{i+j}s!}{i!j!(s-i-j)!}\int_{\R^2}\langle \partial_x^i\partial_y^j\bm{U},\partial_{x}^i\partial_{y}^j\bm{V} \rangle \textup{d}z 
   \end{equation} 
   (cf., e.g., \cite{Michor2007overview}). Note that for $s \in \R$, this metric can be defined via Fourier transforms (cf. e.g. \cite{evans2010}), where the definition coincides with~\eqref{eq:outerMetric} for $s \in \N$.
   While for $s= 0.5$, the defined metric is a weak one, as geodesic distance vanishes for distinct elements in $\Diffc$, the geodesic distance is proved to be positive for $s \geq 1$  in \cite[Theorem 4.1]{Bauer2012geodesic}.

Finally, we want to address the quotient metric on the quotient manifold defined in \eqref{eq:diffDelta} given by $H^s$. For a Riemannian manifold $(M,g)$ and a Lie group $G$ acting isometrically, freely and properly on $M$, a Riemannian metric $\bar{g}$ on $M/G$ can be defined via
\begin{align}
\label{eq:quomet}
\bar{g}_{[p]}(u,v)=g_p(\bar{u}, \bar{v}).
\end{align}
Note that this is the unique metric, such that the quotient map $\pi: M \mapsto M/G$ is a Riemannian submersion (cf. \cite[Proposition 2.28]{Gallot2004riemannian}). Hereby, $T(M/G) \rightarrow TM, v \mapsto \bar{v}$ is given by the horizontal lift. More precisely, let $\pi \colon M \rightarrow M/G$ be the canonical projection and hence $d \pi_p \colon T_pM \rightarrow T_{[p]}(M/G)$ for each $p \in M$, where $[p] \coloneqq \pi(p)$. Since $\pi$ is a submersion, $d \pi_p:(\textup{ker} \ d \pi_p)^{\perp} \rightarrow T_{[p]}(M/G)$ is an isomorphism. Roughly speaking, the horizonzal lift removes the components corresponding to motions within an equivalence class. 

\section{Gradient-based shape optimization on $(\Diffc, H^s)$}
\label{sec:ShapeOpt}

This section discusses the general concepts of shape optimization and our choice of shape space. 
We start by setting up the notation and terminology of the basic concepts on classical shape optimization and on shape optimization on Riemannian manifolds, and presenting a generic Riemannian gradient descent method (\cref{subsec:BasicsShapeCalc}). 
In \cref{subsec:DiffGroupShapeSpace}, results regarding the diffeomorphism group as shape space in optimization are given. The main contributions are the connection between the push-forward and the classical shape derivative as well as the well-posedness of $H^s$ gradients, where the $H^s$ metric was introduced in \eqref{eq:outerMetric}. 
\subsection{Basics on shape calculus}
\label{subsec:BasicsShapeCalc}
In the following, we give notation and terminology of basic shape optimization concepts. 
For a detailed introduction into shape calculus, the reader is referred to the monographs \cite{Delfour2001shapes,Sokolowski1992introduction}. 

One of the key concepts in shape optimization are shape functionals. Formally speaking, they are functions of the type $J\colon \mathcal{D} \rightarrow \R$ that assign to each shape $u \in \mathcal{D}$ a real number. 
The choice of the set $\mathcal{D}$---commonly known as the shape space---plays a crucial role in the design of this kind of problems.

An unconstrained shape optimization problem is given by
\begin{equation}
\label{minproblem}
\min_{u \in \mathcal{D}} J(u).
\end{equation}
Often, shape optimization problems are constrained by equations describing the behavior of the materials and structures or their surroundings, e.g., equations involving an unknown function of two or more variables and at least one partial derivative of this function. 
In this case, the objective functional depends not only on the shapes $u$ but also on a state variable $y$, where the state variable is the solution of the underlying constraint. In other words, one has a shape functional of the form $\hat{J}\colon \mathcal{D} \times \mathcal{Y} \rightarrow \R$ and an operator $e\colon \mathcal{D} \times \mathcal{Y} \rightarrow \mathcal{W}$, where $\mathcal{Y}$ and $\mathcal{W}$ are Banach spaces. 
A constrained shape optimization problem reads as follows:
\begin{align}
&\min_{(u,y)\in \mathcal{D}\times \mathcal{Y}} \hat{J}(u,y), \text{ s.t. }  e(u,y)=0 \label{minproblem2}
\end{align}
The optimization problem is called PDE-constrained, if $e$ in \eqref{minproblem2} represents a PDE.
\begin{remark}
    \label{remarkMinProblems}
    Formally, if the PDE has a (unique) solution given any choice of $u$, then the so-called control-to-state operator $S\colon \mathcal{D} \rightarrow \mathcal{Y}$, $u \mapsto y$ is well-defined. With $J(u) := \hat{J}(u, S u)$ one obtains an unconstrained optimization problem of the form \eqref{minproblem}. This observation justifies working with \eqref{minproblem} in what remains of the section. Although later in \cref{sec:modelProb}, problems of the form \eqref{minproblem2} are presented.
\end{remark}
\paragraph{Shape calculus} 
In classical shape optimization, one usually sets 
$\mathcal{D}$ to be the subset of measurable sets of the power set $P(\holdall)$ of some domain $\holdall$ in $\R^2$. For optimizing with respect to a shape, the concept of shape derivatives is needed.  
In order to define these derivatives, one considers a
family $\{F_t\}_{t\in[0,T]}$ of mappings $F_t\colon \overline{\holdall}\to \mathbb{R}^2$ 
 such that $F_0=\id$, where $\overline{\holdall}$ denotes the closure of $\holdall$ and $T>0$.
This family transforms shapes $u$ into new perturbed shapes of the form
$F_t(u) = \{  F_t(x)\colon x\in u\} . $
Such a transformation can be described by the velocity method or by the perturbation of identity (cf.~\cite[pages 45 and 49]{Sokolowski1992introduction}). The perturbation of identity method is defined by considering
$F_t(x)=F_t^{\bm{W}}(x):= \id (x)+t\bm{W}(x) $ for $x\in u$, where $\bm{W}:\Omega \to \R^2$ is a vector field which guarantees $F_t(u)$ obtained with the perturbation of identity are elements in $\mathcal{D}$. Often, this means that $F_t(u)$ have the same topology as $u$.

Alternatively, one can consider, for $x\in u$, the transformation defined by the velocity method given by $F_t(x)=F_t^{\bm{W}}(x):=\gamma(t,x)$  with $\gamma$ solving
\begin{equation}
\label{eq:velocity}
\partial_t \gamma\left( t,x\right) =\bm{W}(t,\gamma(t,x)) \text{  and } \gamma\left( 0,x\right) =x, 
\end{equation}
where $\bm{W}\colon [0,\tau]\times \overline{\holdall}\to\mathbb{R}^2$  denotes a non-autonomous vector field, which is at least continuous in the first argument for some $\tau>0$ and Lipschitz continuous in the second argument.

For computational reasons, one often considers $\bm{W}$ to be an element in either $W^{1,\infty}(\holdall,\R^2)$ or $C^k_0(\Omega, \R^2)$ for some $k \in \N \cup \infty$ (details in Section \ref{subsubsec:shape-grad}). Setting $F_t=F_t^{\bm{W}}$ to be given either by the perturbation of identity or the velocity method and letting $\mathcal{H} \coloneqq  W^{1,\infty}(\holdall,\R^2)$ or $\mathcal{H} \coloneqq 
 C^k_0(\Omega, \R^2)$, we define, in what follows, the shape derivative. 
\begin{definition}[Shape derivative]\label{def:shapeDer}
	Let $\holdall\subset \mathbb{R}^{n}$ be open and $u\subset \holdall$ measurable.
The Eulerian derivative of a shape functional $J$ at $u$ in direction $\bm{W} \in \mathcal{H}$ is defined by 
\begin{equation}
\label{eulerian}
DJ(u)[\bm{W}]:= \lim\limits_{t \fabove 0}\frac{J(F_t^{\bm{W}}(u))-J(u)}{t}. 
\end{equation}
If for all directions $\bm{W}$ the Eulerian derivative \eqref{eulerian} exists and the mapping 
 $DJ(u) \colon T \mathcal{D} \to \mathbb{R}, \ \bm{W}\mapsto DJ(u)[\bm{W}]$
is linear and continuous, the expression $DJ(u)[\bm{W}]$ is called the shape derivative of $J$ at $u$ in direction $\bm{W}$. In this case, $J$ is called shape differentiable at $u$. 
\end{definition}

\begin{remark}\label{timeInd}
 Note that in the definition of the shape derivative we only consider time-independent vector fields $\bm{W}$ which is not in line with the definition of $F_t$ using the velocity method. 
 As discussed in~\cite[Section 2.2]{Delfour1991velocity}, applications typically focus on cases where the shape derivative depends only on the initial velocity field $\bm{W}(0,x)$ and hence on $\gamma$ solving 
 $$\dfrac{d}{dt}\Big|_{t=0} \gamma(t,x)=\bm{W}(\gamma(0,x))=\bm{W}(x). $$
\end{remark}

\begin{remark}\label{remark:perDomains}
As was discussed in \cite[Section 2.4]{Delfour1990shape}, first-order shape derivatives computed via the perturbation of identity method and the velocity method are identical to each other. 
\end{remark}

The proof of existence of shape derivatives can be done via different approaches like the Lagrangian~\cite{Sturm2013lagrange}, min-max~\cite{Delfour2001shapes}, chain rule~\cite{Sokolowski1992introduction}, rearrangement~\cite{Ito-Kunisch-Peichl} methods, among others.
In this work, we will consider only problems for which the existence of shape derivatives is well-defined and can be achieved by any of the previously mentioned methods, and thus its explicit proof will be omitted.

\paragraph{Optimization on Riemannian shape manifolds} In contrast to classical shape optimization, where $\mathcal{D}$ is usually chosen as a subset of the power set, one can also consider $\mathcal{D}$ as a Riemannian manifold. 
Thanks to the Riemannian metric, the gradient needed for gradient-based optimization techniques is specified. 

In this paper, we consider $\mathcal{D}$ as an infinite-dimensional Riemannian  (shape) manifold $(M,g)$. 
Thus, we will describe how to optimize on a general Riemannian shape manifold $(M,g)$ in the following. An unconstrained and a constrained shape optimization problem on $(M,g)$ is then given by inserting $\mathcal{D}=M$ in \eqref{minproblem} and \eqref{minproblem2}, respectively.

Now, we are ready to formulate and describe the steepest descent algorithm on $(M,g)$ to solve \eqref{minproblem} with $\mathcal{D}=M$ and we present it in \cref{alg:gradAlg} (cf. \cite[Algorithm 1]{Geiersbach2022pde}). 
Hereby, the shape gradient is computed according to \eqref{eq:gradient}. To update the shape iterates, the exponential map is usually used. 
However, since the exponential map on arbitrary Riemannian manifolds $(M,g)$ only exists locally, additional properties of the underlying shape space might be required. An alternative approach to use retractions on Riemannian manifolds can be considered. We will follow this approach in \cref{subsec:retraction}. 
	\begin{algorithm2e}
       \caption{Riemannian steepest descent method on $(M,g)$}
        \label{alg:gradAlg}
	    \DontPrintSemicolon

		\textbf{Require:} Objective function $J$ on $(M,g)$ \;
		
             \textbf{Input:} Initial shape $u^0\in M$, maximum number of iterations $N_{\textup{max}}$ \;
		        
        \While{a suitable termination criterion is not satisfied and $k\leq N_{\textup{max}}$}{ %
        Compute the Riemannian gradient $\bm{V}^k\in T_{u^k}M$ as in~\eqref{eq:gradient}. \nllabel{step:gradient}  \label{lst:algGrad} \; 

    Choose a suitable stepsize $t^k$.\;
    
     Set $u^{k+1}\coloneqq \operatorname{exp}_{u^k}(-t^k \bm{V}^k)$.\nllabel{step:update} \; 

     Set $ k \gets k + 1 $. \;
        }
	\end{algorithm2e}
We now want to connect the definitions from classical shape calculus, in particular the definition of the shape derivative in \cref{def:shapeDer}, to the Riemannian shape optimization approach depicted in \cref{alg:gradAlg}. In this algorithm, the gradient in \cref{lst:algGrad} is computed using the push forward. In particular, the gradient is an element in the tangent space of our considered Riemannian shape manifold $(M,g)$.
On the contrary, in the shape optimization community, it is customary to extract good descent directions from the shape derivative, i.e. we seek a vector field $\bm{W}$ such that $DJ(u)[-\bm{W}]< 0$.
In \cite{Allaire2021shape}, this process is referred to as ``Hilbertian extension-regularization procedure'', and yields more regular deformation fields.  
The key idea is to leverage the structure of an inner product space by computing a deformation field as the Riesz representative of the shape derivative. This deformation field is often referred as \emph{shape gradient}. 

In brief, one can understand the notion of shape gradient as follows:
Let $H$ be a Hilbert space, with inner product $\langle\cdot,\cdot\rangle_H$. The shape gradient $\bm{V}$ is computed by solving the following identification problem
\begin{equation}\label{eq:shapeGrad_def}
\text{Find } \bm{V} \in H \ \text{ such that } \ \langle \bm{V},\bm{W}\rangle_H =DJ(u)[\bm{W}]\quad \forall \bm{W}\in H.
\end{equation}  
The reader should also notice the relation between the identification problem used to compute the shape gradient and the definition of the Riemannian gradient given in equation~\eqref{eq:gradient}, where the Riemannian metric plays the role of the inner product and the Hilbert space $H$ is given by the tangent space of the shape manifold at $u$. However, there is still a discrepancy between the shape derivative and the push-forward, which we will evaluate in the next section for the diffeomorphism group.

\subsection{Diffeomorphism group as shape space in optimization}
\label{subsec:DiffGroupShapeSpace}
\begin{figure}
	\centering
	\includesvg[width=10cm]{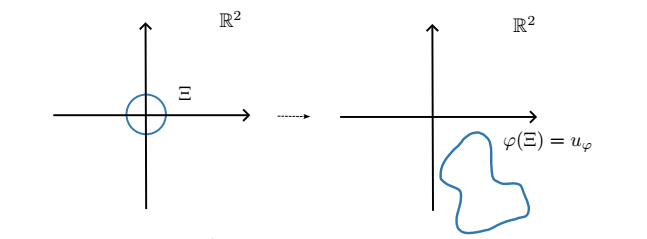}
		\caption{A shape $\varphi \in \Diffc$ is identified with its image of the unit circle $u_{\varphi}=\varphi(\innerCirc)$}
		\label{fig:formdef}%
	\end{figure}

Section~\ref{subsec:BasicsShapeCalc} introduces an optimization approach (cf.~\cref{alg:gradAlg}) applicable to arbitrary Riemannian manifolds. 
In this section, we set $(M,g)=(\Diffc, H^s)$, where $H^s$ is the metric described in~\eqref{eq:outerMetric} or~\eqref{eq:metric2},  extended from the identity to all of $\Diffc$ by right translation, so that $g$ is a right-invariant Riemannian metric on $M$. 
In classical shape optimization, a shape is given by a subset $u \subset \holdall \subset \R^2$, where $\holdall$ is some domain in $\R^2$ while in the optimization on a Riemannian shape manifold $(M,g)$, a shape is given by an element $u \in M$ (cf. \cref{subsec:BasicsShapeCalc}). 
As we consider the shape space $(\Diffc, H^s)$, the link between both definitions becomes apparent through identifying $\varphi$ with its image of the unit circle $\varphi(\innerCirc)$ (cf. \cref{fig:formdef}). 
Note that this identification is only well-defined if we consider $\varphi \in  \Diffc/\text{Diff}(\R^2, \innerCirc)$, where the quotient space was introduced in \eqref{eq:diffDelta}. 
As usual, we define a shape functional on $(\Diffc, H^s)$ as a function $j\colon \Diffc \ra \R$ and denote by $\mathcal{A}(\Diffc)$ the space of all shape functionals on $\Diffc$. 

Consider $G\coloneqq\text{Diff}_c(\R^2,\innerCirc)$, acting on $M$ on the right by composition, $\varphi\cdot\psi\coloneqq\varphi\circ\psi$, and let $\pi\colon M\ra M/G$ be the canonical projection.
Additionally, consider
\begin{align*}
\mathcal{J} \coloneqq\bigl\{\,j \in \mathcal{A}(\Diffc)\ \bigl|\ \exists\,
J \colon \mathcal{D}\subset P(\R^2) \ra \R \text{ s.t. }
j(\varphi)=J(\varphi(\innerCirc))\ \ \forall \varphi \in \Diffc \,\bigr\}.
\end{align*}
Every $j\in\mathcal J$ is thus represented by a shape functional $J$ from classical shape optimization theory, and is invariant under the action of $G$. 
Indeed, for $\psi\in G$ one has $(\varphi\circ\psi)(\innerCirc)=\varphi(\psi(\innerCirc)) =\varphi(\innerCirc)$, so that $j(\varphi\circ\psi)=j(\varphi)$.
Consequently $j$ factors as $j=\tilde\jmath\circ\pi$ for a unique functional $\tilde\jmath$ on $M/G$.

The following lemma shows that for a shape functional $j \in \mathcal{J}$ and the Riemannian steepest descent method described in~\cref{alg:gradAlg} with $(M,g)=(\Diffc, H^s)$, the computed gradients have no vertical components and their projection coincides with the gradient on the quotient space $\Diffc/\text{Diff}(\R^2, \innerCirc)$.

\begin{lemma}
\label{lem:submersion-gradient}
Let $(M,g)$ be a Riemannian manifold and let the Lie group $G$ act on $M$ isometrically, freely and properly, so that $M/G$ carries the
quotient metric $\tilde g$ of~\eqref{eq:quomet} and $\pi\colon M\ra M/G$ is a Riemannian submersion; write $v\mapsto\bar v$ for the horizontal lift, that is, for the inverse of the isomorphism $\dd\pi_\varphi\colon(\ker \dd\pi_\varphi)^{\perp}\ra T_{[\varphi]}(M/G)$.
Let $j\in\mathcal A(M)$ be $G$-invariant, with $j=\tilde\jmath\circ\pi$, and let $\varphi\in M$ be such that the gradient $\nabla j(\varphi)$ exists in $T_\varphi M$. 
Then
\begin{enumerate}[label=(\roman*)]
\item
$\nabla j(\varphi)\in(\ker \dd\pi_\varphi)^{\perp}$, that is,
$\nabla j(\varphi)$ has no vertical component;
\item
the gradient of $\tilde\jmath$ at $[\varphi]\coloneqq\pi(\varphi)$
exists and satisfies
$\nabla\tilde\jmath([\varphi])=\dd\pi_\varphi\bigl(\nabla j(\varphi)\bigr)$.
\end{enumerate}
\end{lemma}

\begin{proof}
For the proof of statement~(i), let $v\in\ker \dd\pi_\varphi$.
Then
\[
g_{\varphi}\bigl(\nabla j(\varphi),v\bigr)
=\dd j_{\varphi}(v)
=\dd(\tilde\jmath\circ\pi)_{\varphi}(v)
=\dd\tilde\jmath_{\pi(\varphi)}\bigl(\dd\pi_{\varphi}(v)\bigr)
=0 ,
\]
the last equality because $\dd\pi_\varphi(v)=0$.
As $v$ was arbitrary in $\ker \dd\pi_\varphi$, the gradient is orthogonal to that subspace.

For~(ii), let $v\in T_{[\varphi]}(M/G)$.
By~(i) the field $\nabla j(\varphi)$ is horizontal and is carried by $\dd\pi_\varphi$ to $\dd\pi_\varphi(\nabla j(\varphi))$, so it is the horizontal lift of the latter.
Using the definition~\eqref{eq:quomet} of the quotient metric,
\[
\tilde g_{[\varphi]}\bigl(\dd\pi_{\varphi}(\nabla j(\varphi)),v\bigr)
=g_{\varphi}\bigl(\nabla j(\varphi),\bar v\bigr)
=\dd j_{\varphi}(\bar v)
=\dd\tilde\jmath_{[\varphi]}\bigl(\dd\pi_{\varphi}(\bar v)\bigr)
=\dd\tilde\jmath_{[\varphi]}(v) ,
\]
the last equality because $\dd\pi_\varphi(\bar v)=v$ by definition of
the horizontal lift.
Since $v$ was arbitrary, $\dd\pi_\varphi(\nabla j(\varphi))$ represents
$\dd\tilde\jmath_{[\varphi]}$ with respect to $\tilde g$, which is the
assertion.
\end{proof}

Notice that \cref{lem:submersion-gradient} assumes the existence of $\nabla j(\varphi)$.
This is a hypothesis rather than a consequence of the Riemannian structure. 
The metric~\eqref{eq:outerMetric} on $\Diffc$ is indeed weak, since the tangent spaces are modelled on $\Xc(\R^2)$, which is not complete in the $H^s$ norm, so the induced map from $T_\varphi M$ to its dual is injective but need not be surjective, and a differential need not be represented by a tangent vector.
The gradients used below are therefore representatives in a completion of $\Xc(\R^2)$, and the identification of shapes with elements of $\Diffc$ has to be understood accordingly.

\begin{remark}\label{rem:conventions}
From this point on and throughout the paper, we only consider shape functionals $j \in \mathcal{J}$ and will identify them with a classical shape functional $J \colon \mathcal{D}  \subset P(\R^2) \ra \R$. Therefore, the shape functional $J$ is invariant under reparametrizations by elements in $\text{Diff}_c(\R^2, \innerCirc)$. Moreover, we will identify a shape $\varphi \in \Diffc$ with its image of the unit circle $u_{\varphi}=\varphi({\innerCirc})$, as depicted in \cref{fig:formdef} and write $u=u_{\varphi} $ when the context is clear. In that setting, by slight abuse of notation, it holds that $u \in \Diffc$ and $u \subset \R^2$. 
By \cref{lem:submersion-gradient} the descent method may then be run directly on $(\Diffc,H^s)$, the iterates being independent of the representative chosen within each equivalence class.
In this setting, it should be emphasized that that due to the previous discussions, we are justified in perfoming the Riemannian descent method directly on $(\Diffc, H^s)$.
\end{remark} 

\subsubsection{Equivalence of shape derivative and push-forward} 
\label{subsubsec:ShapeDerPushForw}
In the context of the previous section, it should be noted that one can consider using the shape derivative introduced in~\eqref{eulerian} to compute the shape gradient introduced in \eqref{eq:shapeGrad_def} instead of using the push-forward as done in \cref{alg:gradAlg}. We show the equivalence of both approaches for our considered shape space.

\begin{theorem} \label{thm:derivate}
Let $j \in \mathcal{J}$ be a shape functional such that
$j(\varphi)=J(u_{\varphi})$ for all $ \varphi \in \Diffc$
for some $J\colon \mathcal{D}  \subset P(\R^2) \ra \R$. Then, for any tangent curve $\gamma\colon \R \ra \Diffc$ in $\gamma(0)=\varphi$, it holds that
$
dj_{\varphi}(\gamma)=DJ(u_{\varphi})[\bm{W}],
$
where $\gamma \cong \bm{W} \in \Xc(\R^2) \subset C^{\infty}(\R^2, \R^2)$ as in \eqref{eq:tangentDiff}.
\end{theorem}

\begin{proof}
Initially, we observe that for the identification of $\bm{W}$ with $\gamma$ described in \eqref{eq:tangentDiff} and by setting $\tilde{\gamma}(t,x) \coloneqq (\gamma \circ \varphi^{-1})(t,x)$ we obtain
\begin{align*}
     \dfrac{d}{dt}\Big|_{t=0} \tilde{\gamma} \left( t,x\right) =\dfrac{d}{dt}\Big|_{t=0} \gamma(t, \varphi^{-1} (x)) = \bm{W}(\gamma(t,\varphi^{-1}(0))) =\bm{W}(\tilde{\gamma}(0,x))\text{  and } \tilde{\gamma}\left( 0,x\right) =x.
	\end{align*}
Considering \cref{timeInd}, we obtain perturbation functions $F_t$ for the vector field $\bm{W} $ as in \eqref{eq:velocity} by setting $F_t(x):=\tilde{\gamma}(t,x)$. 
 From the definition of the push-forward, we follow directly
\begin{align*}
dj_{\varphi}(\gamma)&=  \dfrac{d}{dt} \bigg|_{t=0} j(\gamma(t)) =\liml_{t \fabove 0} \dfrac{j(\gamma(t))-j(\gamma(0))}{t}=\liml_{t \fabove 0} \dfrac{J(\gamma(t, \innerCirc))-J(u_{\varphi})}{t}\\
&=\liml_{t \fabove 0} \dfrac{J(\tilde{\gamma}(t,u_{\varphi}))-J(u_{\varphi})}{t}=\liml_{t \fabove 0} \dfrac{J(F_t(u_{\varphi}))-J(u_{\varphi})}{t} =DJ(u_{\varphi})[\bm{W} ].
 \end{align*}
\end{proof}

\begin{remark}
  From \cref{remark:perDomains}, we know that the computation of the shape derivative does not depend on the choice between the velocity method and the perturbation of identity. Inferentially, the result of \cref{thm:derivate} can be transferred to the computation of the shape derivative using the perturbation of identity method.
\end{remark}

For clarity, we recall that by \Cref{thm:derivate}, together with the use of diffeomorphisms as our shape space, the shape gradient is computed by the following identification problem 
\[
\text{Find }\bm{V}\in \Xc(\R^2) \ \text{such that} \  H^s(\bm{V},\bm{W}) = DJ(u_\varphi)[\bm{W}]\qquad \forall \ \bm{W}\in \Xc(\R^2),
\]
where $H^s(\cdot,\cdot)$ is the Riemannian metric given in \eqref{eq:metric2}, and the right-hand side is given by the shape derivative (cf. \cref{def:shapeDer}). 
As previously mentioned, in classical shape optimization one aims to compute the shape gradient as the Riesz representative of the shape derivative. However, $\Xc(\R^2)$ equipped with the $H^s$ metric is not a Hilbert space. This issue will be discussed in the next section, where we will also provide an explanation of the computation of the shape gradient.

So far the metric and the push-forward have been formulated on all of $\R^2$, whereas the applications of \cref{sec:modelProb} are posed on a bounded domain. We therefore fix the following conventions, which are used for the remainder of the paper and which determine the boundary behaviour of the vector fields appearing throughout.

\begin{remark}
\label{rem:compSup}
For the remainder of the paper we consider shapes as elements of $\Diffc$ equal to the identity outside the bounded domain $\holdall\subset\R^2$, a set we denote by $\DiffcOm$; their tangent vectors then belong to $C_0^\infty(\holdall,\R^2)$.
In virtue of \cref{rem:conventions} we identify a shape with the image of $\innerCirc$ under $\varphi\in\Diffc$, which requires $\innerCirc\subset\holdall$. 
This assumption imposes no restriction as the circle may be scaled and translated to suit the hold-all domain fixed by the application.
The framework extends to several shapes at once by identifying a multi-shape with the image of disjoint circles, see~\cite[Sec.~2]{Geiersbach2022pde}.
\end{remark}

\subsubsection{Shape gradients with respect to the $H^s$ metric}
\label{subsubsec:shape-grad}

Having identified the push-forward with the classical shape derivative, we now turn to the computation of the shape gradient, that is, to representing $DJ(u)$ with respect to the metric~\eqref{eq:outerMetric}. 
Two difficulties arise. 
The metric is defined on $\Xc(\R^2)$, whereas the applications of \cref{sec:modelProb} are posed on a bounded hold-all domain; and the operator it induces is of order $2s$, so that its inversion is not directly accessible to standard finite element software. 
This subsection addresses the first difficulty and prepares the second.

We begin with the functional setting. 
The completion of $\Xc(\R^2)$ is the Sobolev space $H^s(\R^2,\R^2)$ (cf.\ e.g.~\cite{adams2003sobolev}), so that $\Xc(\R^2)$ is continuously embedded in $H^s(\R^2,\R^2)$ for all $s\in\N$. 
Using the $L^2$ inner product one identifies the dual of $H^s(\R^2,\R^2)$ with $H^{-s}(\R^2,\R^2)$, and the operator $(\textup{id}-A\Delta)^s\colon H^s(\R^2,\R^2)\to H^{-s}(\R^2,\R^2)$ is then a bounded linear operator. 
Under periodicity assumptions this operator is invertible and can be inverted efficiently in the Fourier domain (see e.g.~\cite{Beg2005computing}); we do not pursue that route here, since the hold-all domains of \cref{sec:modelProb} are bounded and carry boundary conditions.

Following \cite[Sec.~5.1]{Michor2007overview}, the $H^s$ inner product on $\mathfrak{X}_c(\R^2)$ admits the three equivalent expressions
\begin{align}
\label{eq:metric_whole_space}
H^s(\bm V,\bm W)
&=\sum_{\substack{i,j\ge0\\ i+j\le s}}
\frac{A^{i+j}s!}{i!\,j!\,(s-i-j)!}
\int_{\R^2}\big\langle\partial_x^i\partial_y^j\bm V,\
\partial_x^i\partial_y^j\bm W\big\rangle\dd x \nonumber\\
&=\sum_{\substack{i,j\ge0\\ i+j\le s}}
\frac{(-A)^{i+j}s!}{i!\,j!\,(s-i-j)!}
\int_{\R^2}\big\langle\partial_x^{2i}\partial_y^{2j}\bm V,\bm W\big\rangle\dd x \nonumber\\
&=\int_{\R^2}\big\langle(\textup{id}-A\Delta)^s\bm V,\bm W\big\rangle\dd x .
\end{align}
The equalities are obtained by using integration by parts, and it is exact only because the fields have compact support, so that every boundary term vanishes. 
This is consistent with \cref{rem:compSup}, where tangent vectors to $\Diffc$ are taken in $C_0^\infty(\R^2,\R^2)$.

On a bounded hold-all the three expressions cease to agree for general fields, since integration by parts now produces boundary terms whose vanishing is not guaranteed only by belonging to $H^s(\holdall,\R^2)$. 
Rather than carry~\eqref{eq:metric_whole_space} over to the bounded domain, we take the split system introduced below as the definition of the shape gradient on $\holdall$. 
\cref{prop:split-equivalence} identifies the problem it solves and shows that it reproduces~\eqref{eq:metric_whole_space} on fields supported away from $\partial\holdall$.

Throughout, $\holdall\subset\R^2$ is a bounded Lipschitz domain, so that traces on $\partial\holdall$ are defined.
Its boundary splits as $\partial\holdall=\Gamma\sqcup\Sigma$, with $\Gamma$ held fixed during the optimisation and $\Sigma$, possibly empty, free to move.
We take both to be unions of connected components of $\partial\holdall$, so that $\ol\Gamma\cap\ol\Sigma=\emptyset$; this separation is all that is used below.
In \cref{subsec:interface} the moving boundary is an interior interface, so $\Gamma=\partial\holdall$ and $\Sigma=\emptyset$, while in \cref{subsec:compliance} the shape consists of four holes in the material domain, $\Gamma$ being the outer boundary and $\Sigma$ the union of the four moving components.
Deformation fields are accordingly sought in
\[
H^1_{\Gamma,0}(\holdall,\R^2)
\coloneqq\big\{\bm W\in H^1(\holdall,\R^2)\colon \gamma_0\bm W=0 \text{ on }
\Gamma\big\},
\]
where $\gamma_0\colon H^1(\holdall,\R^2)\to H^{1/2}(\partial\holdall,\R^2)$ denotes the trace operator, acting componentwise on vector fields, which is bounded and surjective since $\holdall$ is a bounded Lipschitz domain.
Because $\Gamma$ and $\Sigma$ are unions of connected components of $\partial\holdall$, each is relatively open and closed there and restriction to it is well defined; we write
\[
\gamma_\Gamma\bm W\coloneqq\gamma_0\bm W|_{\Gamma},
\qquad
\gamma_\Sigma\bm W\coloneqq\gamma_0\bm W|_{\Sigma},
\]
so that the condition defining $H^1_{\Gamma,0}(\holdall,\R^2)$ reads $\gamma_\Gamma\bm W=0$ in $H^{1/2}(\Gamma,\R^2)$, no condition being imposed on the remaining portion.
Notice that $H^1_{\Gamma,0}(\holdall,\R^2)$ reduces to $H^1_0(\holdall,\R^2)$ when $\Sigma=\emptyset$.
The same notation is used for other exponents. In other words, for $1<q<\infty$ we write $W^{1,q}_{\Gamma,0}(\holdall,\R^2)\coloneqq\{\bm W\in W^{1,q}(\holdall,\R^2): \gamma_\Gamma\bm W=0\}$, the trace being understood in $W^{1-1/q,q}(\Gamma,\R^2)$.

Each stage of the split system introduced below is a second-order problem, tested against a single field $\bm W\in H^1_{\Gamma,0}(\holdall,\R^2)$, and the stage carrying the shape derivative has right-hand side $DJ(u)[\bm W]$.
For that stage to be well posed, $DJ(u)$ must act boundedly on $H^1_{\Gamma,0}(\holdall,\R^2)$, which is more than the shape derivative provides in general. 
We therefore assume the following throughout.

\begin{assumption}
\label{assmp:dj-bounded}
The shape functional $J$ is shape differentiable at $u$, and $DJ(u)$ extends
to a bounded linear functional on $H^1_{\Gamma,0}(\holdall,\R^2)$.
\end{assumption}

Drawing inspiration from Gallistl's work~\cite{Gallistl2017stable}, where a stable splitting of $(2m)$-th order elliptic partial differential equations into $2(m-1)$ Poisson-type problems and one generalized Stokes problem was proposed, we compute the deformation field on $\holdall$ by splitting the operator $(\textup{id}-A\Delta)^s$ into $s$ second-order problems. Introducing the auxiliary fields $\bm X_1,\dots,\bm X_{s-1}$ and writing $\bm X_0\coloneqq\bm V$, we consider the following system:
\bigskip

For $j=0,1,\dots,s-1$, set $\bm X_0:=\bm V$ and seek $\bm X_0,\bm X_1,\dots,\bm X_{s-1}\in H_{\Gamma,0}^1(\holdall,\R^2)$ such that, for every $\bm W\in H_{\Gamma,0}^1(\holdall,\R^2)$,
\begin{equation}
\label{eq:grad_system_reindexed}
\int_{\holdall}\bm X_j\cdot\bm W\,\dd x + A\int_{\holdall}\nabla\bm X_j:\nabla\bm W\,\dd x
=
\begin{cases}
\displaystyle\int_{\holdall}\bm X_{j+1}\cdot\bm W\,\dd x, & j=0,\dots,s-2,\\[2mm]
DJ(u)[\bm W], & j=s-1.
\end{cases}
\end{equation}
Only the condition on $\Gamma$ is imposed, through the choice of test and trial space; the conditions on $\Sigma$ are natural.

Under \cref{assmp:dj-bounded}, the right-hand side of every equation in~\eqref{eq:grad_system_reindexed} is bounded on $H^1_{\Gamma,0}(\holdall,\R^2)$.
For $j=s-1$ this is precisely what the assumption provides, while for $j\le s-2$ the right-hand side is the field $\bm X_{j+1}$ obtained at the previous stage, which lies in $L^2(\holdall,\R^2)$.
The right-hand side there is tested against every $\bm W\in H^1_{\Gamma,0}(\holdall,\R^2)$, including fields whose trace on $\Sigma$ does not vanish, and \cref{assmp:dj-bounded} supplies $DJ(u)$ on exactly that space.
Knowing $DJ(u)$ only on $C_0^\infty(\holdall,\R^2)$ would not be enough: a distribution on the open set $\holdall$ is blind to anything supported on $\partial\holdall$, so two data differing by a term of the form $\bm W\mapsto\int_\Sigma\bm g\cdot\gamma_0\bm W\,\ds$ induce the same distribution while producing different solutions of~\eqref{eq:grad_system_reindexed}.
The assumption is what gives the conormal condition of \cref{prop:split-equivalence}\,(i) its meaning.
A distribution on $\holdall$ does not determine a conormal derivative on $\Sigma$, and the two objects differ by exactly the freedom in extending it to the test space. 
Thus, we use \cref{assmp:dj-bounded} to remove this degree of freedom by making the extension part of the datum.
Solving the stages in decreasing order of $j$ then determines the whole tuple, and each stage improves on the regularity delivered by the preceding one, as quantified in \cref{prop:split-equivalence}\,(i).

The following result guarantees that the fields produced by~\eqref{eq:grad_system_reindexed} are more regular than the formulation suggests.
It also shows that $\bm V$ inherits $s$ boundary conditions of Navier type on $\Gamma$, from the vanishing traces of the $\bm X_j$, together with $s$ further conditions on $\Sigma$ carried by the natural conditions of the stages.
It characterises the system through an equivalent strong formulation, consisting of the equation $(\id-A\Delta)^s\bm V=DJ(u)$ together with both families of boundary conditions.
Finally, it shows that when tested against fields supported in $\holdall$, the field $\bm V$ satisfies the multi-index identity of~\eqref{eq:metric_whole_space} exactly, so that the system does compute the $H^s$ gradient of $J$ away from the boundary of the hold-all domain.

The conditions on $\Sigma$ absent from~\eqref{eq:grad_system_reindexed} are recovered in \cref{prop:split-equivalence}\,(i), and all but one of them take the classical form $\partial_n\bm X_j=0$, where $\partial_{n}$ denotes the normal derivative.  
The exception is the one belonging to the stage $j=s-1$.

\cref{assmp:dj-bounded} asks of the shape derivative exactly what the stage $j=s-1$ needs, namely boundedness on $H^1_{\Gamma,0}(\holdall,\R^2)$, and the regularity it delivers matches; in other words,  $\bm X_{s-1}\in H^1(\holdall,\R^2)$, with $\nabla\bm X_{s-1}\in L^2(\holdall,\R^{2\times2})$.
However, given the regularity of $\nabla \bm X_{s-1}$, the conormal derivative must be defined as a variational object rather than in terms of traces, and \cref{def:gen-conormal} is the form it takes.

\begin{definition}
\label{def:gen-conormal}
Let $F\in\big(H^1_{\Gamma,0}(\holdall,\R^2)\big)^{*}$ and let $\bm Z\in H^1(\holdall,\R^2)$ satisfy $(\id-A\Delta)\bm Z=F$ in $\left(C^\infty_0(\holdall, \R^2)\right)^\ast$ with $A>0$.
The \emph{generalised conormal derivative of $\bm Z$ relative to $F$} is the unique element $\partial_n^{F}\bm Z\in H^{-1/2}(\Sigma,\R^2)$ such that
\begin{equation}
\label{eq:gen-conormal}
A\,\big\langle \partial_n^{F}\bm Z,\ \gamma_\Sigma\bm W\big\rangle_{H^{-1/2}(\Sigma),H^{1/2}(\Sigma)}
=\int_\holdall \bm Z\cdot\bm W\,\dd x
+A\int_\holdall\nabla\bm Z:\nabla\bm W\,\dd x
-F[\bm W]
\end{equation}
for every $\bm W\in H^1_{\Gamma,0}(\holdall,\R^2)$.
\end{definition}

Notice that when $\Sigma=\emptyset$ this notion is not required anymore. 
Additionally, for the operator $\id-A\Delta$, the conormal derivative coincides with $A\,\partial_n$, the coefficient $A$ being scalar.  
We keep the term ``conormal'' because it is standard in the literature cited below.
Existence and uniqueness, and the fact that $\partial_n^{F}\bm Z=\partial_n\bm Z$ whenever $\bm Z\in H^2(\holdall,\R^2)$ and $F\in L^2(\holdall,\R^2)$ is shown in \cref{lem:gen-conormal}.

Defining a conormal derivative through the first Green identity, for fields whose Laplacian is only a functional, goes back to Costabel~\cite{costabel1988boundary} and is treated for systems on Lipschitz domains in~\cite[Ch.~4]{mclean2000strongly}.
Our $\partial_n^{F}\bm Z$ coincides with the generalised conormal derivative relative to the datum, at $s=1$. 
The object depends on the extension of the right-hand side from $C_0^\infty(\holdall,\R^2)$ to the test space, and so is not determined by the equation alone; \cref{assmp:dj-bounded} fixes that choice here by supplying $DJ(u)$ already as a functional on $H^1_{\Gamma,0}(\holdall,\R^2)$.

\begin{theorem}
\label{prop:split-equivalence}
Let $s\ge1$ and $A>0$, and let $\holdall\subset\R^2$ be a bounded Lipschitz domain whose boundary decomposes as $\partial\holdall=\Gamma\sqcup\Sigma$, where $\Gamma$ and $\Sigma$ are (possibly empty) unions of connected components of $\partial\holdall$; in particular $\ol\Gamma\cap\ol\Sigma=\emptyset$.
Assume that every connected component of $\partial\holdall$ is of class $C^{2s-2,1}$ when $s\ge2$, and let \cref{assmp:dj-bounded} hold.
Then~\eqref{eq:grad_system_reindexed} has a unique solution $(\bm X_0,\dots,\bm X_{s-1})\in\big(H^1_{\Gamma,0}(\holdall,\R^2)\big)^s$, and, writing $\bm V\coloneqq\bm X_0$:
\begin{enumerate}[label=(\roman*)]
\item
$\bm X_j\in H^{2(s-1-j)+1}(\holdall,\R^2)$ for $j=0,\dots,s-1$; in particular
$\bm V\in H^{2s-1}(\holdall,\R^2)\hookrightarrow H^{s}(\holdall,\R^2)$, with
\begin{equation}
\label{eq:prop-apriori}
\|\bm V\|_{H^{2s-1}(\holdall,\R^2)}
\le C\,\big\|DJ(u)\big\|_{(H^1_{\Gamma,0}(\holdall,\R^2))^{*}},
\end{equation}
where $C$ is a constant depending on $s,A,\holdall,\Gamma,\Sigma$.
Moreover $\bm X_j=(\id-A\Delta)^j\bm V$ as an identity in $H^{2(s-1-j)+1}(\holdall,\R^2)$, and $\bm V$ satisfies the $2s$ boundary conditions
\begin{align}
\label{eq:bc-dirichlet}
\gamma_\Gamma\big((\id-A\Delta)^j\bm V\big)&=0
&&\text{on }H^{1/2}(\Gamma,\R^2), && j=0,\dots,s-1,\\
\label{eq:bc-neumann}
\partial_n\big((\id-A\Delta)^j\bm V\big)&=0
&&\text{on }H^{1/2}(\Sigma,\R^2), && j=0,\dots,s-2,\\
\label{eq:bc-neumann-top}
\partial_n^{DJ(u)}\big((\id-A\Delta)^{s-1}\bm V\big)&=0
&&\text{on }H^{-1/2}(\Sigma,\R^2), &&
\end{align}
the last in the sense of \cref{def:gen-conormal}.
\item
A field $\bm V\in H^{2s-1}(\holdall,\R^2)$ coincides with the component $\bm X_0$ of
the solution of~\eqref{eq:grad_system_reindexed} if and only if
\[
(\id-A\Delta)^s\bm V=DJ(u)\ \text{ in }\big(C_0^\infty(\holdall,\R^2)\big)^{*}
\qquad\text{and}\qquad
\eqref{eq:bc-dirichlet}\text{--}\eqref{eq:bc-neumann-top}\ \text{hold},
\]
and such a $\bm V$ is unique.
\item
For every $\bm W\in C_0^\infty(\holdall,\R^2)$,
\[
\sum_{\substack{i,j\ge0\\ i+j\le s}}
\frac{A^{i+j}s!}{i!\,j!\,(s-i-j)!}
\int_\holdall\big\langle\partial_x^i\partial_y^j\bm V,\
\partial_x^i\partial_y^j\bm W\big\rangle\dd x
= DJ(u)[\bm W].
\]
That is, $\bm V$ is the $H^s$ gradient of $J$ with respect to~\eqref{eq:metric_whole_space} when tested against fields supported in $\holdall$. 
\end{enumerate}
\end{theorem}
\noindent
The proof is given in \cref{sec:proof_prop_2}.

When $\Sigma=\emptyset$ the conditions~\eqref{eq:bc-neumann}--\eqref{eq:bc-neumann-top} impose no conditions, $H^1_{\Gamma,0}(\holdall,\R^2)=H^1_0(\holdall,\R^2)$, $C_0^\infty(\holdall,\R^2)$ is dense in that space, and~(ii) in \cref{prop:split-equivalence} reduces to the equation together with the $s$ Navier-type conditions~\eqref{eq:bc-dirichlet} on the whole of $\partial\holdall$.
For $\Sigma\ne\emptyset$ that density fails, so the equation alone does not see the action of $DJ(u)$ through traces on $\Sigma$, and~\eqref{eq:bc-neumann-top} is what restores it.
The conditions~\eqref{eq:bc-dirichlet} on $\Gamma$ are of Navier type: only the member $j=0$ constrains $\bm V$ itself, the remaining $s-1$ constraining successive powers of $(\id-A\Delta)$ applied to it.
We remark that they are not the clamped conditions $\bm V=\partial_n\bm V=0$ that a splitting into a different family of second-order problems would encode.

The two hypotheses of \cref{prop:split-equivalence} carry different weight. 
Unique solvability of~\eqref{eq:grad_system_reindexed} follows from coercivity of the stage form, which holds for any decomposition $\partial\holdall=\Gamma\sqcup\Sigma$ and needs no regularity of $\partial\holdall$ beyond that required for traces; the deformation field is therefore well defined on any bounded Lipschitz hold-all domain, and \cref{alg:gradAlg} is unaffected as it operates on~\eqref{eq:grad_system_reindexed} alone.
The higher regularity of $\partial\holdall$ buys the regularity of the auxiliary fields, and with it the identity $\bm X_j=(\id-A\Delta)^j\bm V$ on which the boundary conditions and the strong formulation rest.
It enters only through the iterated elliptic regularity argument, at whose last step
$C^{2s-2,1}$ is required. 
Despite the regularity hypothesis not being satisfied by any of the applications of \cref{sec:modelProb}, for $s\ge2$,  part~(iii) in \cref{prop:split-equivalence} does hold in any case. 
Its proof reduces the left-hand side to the duality pairing $\big\langle(\id-A\Delta)^s\bm V,\bm W\big\rangle$, which requires only $\bm V\in H^{s}_{\mathrm{loc}}(\holdall,\R^2)$ and is supplied by interior elliptic regularity; the equation itself follows from testing the stage $j=s-1$ against $C_0^\infty(\holdall,\R^2)$.

Part~(iii) of \cref{prop:split-equivalence} holds only for test fields supported in $\holdall$, and that restriction cannot be removed.
Already for $s=2$ with $\Sigma=\emptyset$, eliminating the auxiliary field from~\eqref{eq:grad_system_reindexed} shows that $\bm V$ represents $DJ(u)$ on $H^2(\holdall,\R^2)\cap H^1_0(\holdall,\R^2)$ for an inner product whose second-order term involves the Laplacian, whereas~\eqref{eq:metric_whole_space} restricted to $\holdall$ involves the full Hessian.
On that space the two forms differ by a boundary integral weighted by the curvature of $\partial\holdall$ \cite[Lem.~3.1.1.1]{grisvard2011elliptic}, so they are equivalent but not equal, and the deformation field is the gradient for the first.
The boundary term is governed by the normal derivative of the field on $\partial\holdall$, which is also why the extension of $\bm V$ by zero does not lie in $\Xc(\R^2)$, as noted in \cref{rem:compSup}.

Two features of the exponent in~(i) are worth noting.
It is sharp under \cref{assmp:dj-bounded}: the stage $j=s-1$ absorbs one derivative, its datum being a functional rather than a density, and the $s-1$ stages below it gain two apiece. It is also more than the metric requires: had each stage gained only $\delta$ derivatives, the outcome would be $\bm V\in H^{1+(s-1)\delta}(\holdall,\R^2)$, and $1+(s-1)\delta\ge s$ already for $\delta\ge1$, so a gain of one derivative per stage suffices for $\bm V\in H^{s}(\holdall,\R^2)$.

\section{Applications}
\label{sec:modelProb}

In this section, we formulate two shape optimization models, which we use to demonstrate the performance of \cref{alg:gradAlg} using different Riemannian metrics (inner and outer metrics) on the corresponding shape space.
For the first experiment, we consider an interface identification problem studied in a number of texts such as \cite{Buttazzo2003optimal,Ito-Kunisch-Peichl,Sokolowski1992introduction}. Electrical impedance tomography, where the material distribution of electrical properties such as electric conductivity and permittivity inside the body is to be determined \cite{Cheney1999electrical,Kwon2002magnetic}, is a motivation example for this model. 
The second example is motivated by structural mechanics, and aims to solve the compliance minimization of a two-dimensional bridge. This problem has been solved for example in \cite{Allaire2014shape,Herzog2024discretize}. 

Throughout this section, we consider the Riemannian manifold $(\DiffcOm, H^s)$ as the underlying shape space, where the space $\DiffcOm$ was introduced in Remark \ref{rem:compSup}. 
For ease of notation, from now on, we will denote as $(M,g)$ the previously mentioned shape space.
We remain in the setting of \cref{subsec:DiffGroupShapeSpace}, in particular considering Remark \ref{rem:conventions}, Remark \ref{rem:compSup} and the findings in \cref{subsubsec:shape-grad}. Note that due to the assumptions in Remark \ref{rem:compSup}, it is justified to consider a shape $u \in \DiffcOm$ to be a subset of our chosen hold-all domain $\holdall$. Construction of the hold-all domain $\holdall$ will be done for each application individually. 
\subsection{Interface identification problem}
\label{subsec:interface}
This section describes an interface identification problem inspired by electric impedance tomography. 
Briefly speaking, the main goal is to determine the unknown interface between different regions based on measurements such that it minimizes a certain objective functional. 
We assume we have access to measurements distributed on the domain we are interested in. 
Notice that we intentionally use the notation $\holdall$ despite this being used to describe the hold-all domain, since without loss of generality, we assume they coincide. 

The hold-all domain $\holdall$ is a bounded domain representing the body under investigation, and the shape $u\in M$  is a simple closed curve compactly contained in $\holdall$, which separate it into an interior region $\Omega_{\textup{in}}$, enclosed by the curve, and the remaining region $\Omega_{\textup{out}}$, such that $\holdall = \Omega_{\textup{in}} \sqcup u \sqcup \Omega_{\textup{out}}$. 
We recall that, in virtue of \cite[Sec.~2]{Geiersbach2022pde}, this framework can be extended to a multi-shape problem where $u=\bigcup_{i=1}^{N} u_i$ consists of finitely many pairwise disjoint simple closed curves compactly contained in $\holdall$; however, we do not pursue this approach here. 
The conductivity $\kappa$ is piecewise constant with respect to this partition, taking one value on each region, and $u$ is the interface across which it jumps. 
A current $g$ is applied on $\partial\holdall$, the resulting potential is measured there, and the shape is sought so that the measured and predicted potentials agree.

In the decomposition $\partial\holdall=\Gamma\sqcup\Sigma$ of \cref{subsubsec:shape-grad}, this is the case $\Gamma=\partial\holdall$ and $\Sigma=\emptyset$. 
The shape is an interior interface rather than a portion of the boundary of the hold-all, which is held fixed in its entirety. 
The particular configuration used in our computations, i.e., the hold-all domain, the number and initial position of the interfaces, the two conductivity values and the applied current, is specified in \cref{subsec:Expe1}.

\begin{figure}[h]
\tikzset{every picture/.style={line width=0.55pt}} 
\centering
    \resizebox{!}{0.3\textheight}{%
\begin{tikzpicture}[x=0.75pt,y=0.75pt,yscale=-1,xscale=1]

\draw  [fill={rgb, 255:red, 177; green, 105; blue, 50 }  ,fill opacity=1 ] (188,21) -- (350,21) -- (350,171) -- (188,171) -- cycle ;
\draw  [color={rgb, 255:red, 255; green, 255; blue, 255 }  ,draw opacity=1 ][fill={rgb, 255:red, 63; green, 126; blue, 171 }  ,fill opacity=1 ] (242,62) .. controls (262,52) and (327,45) .. (307,65) .. controls (287,85) and (287,95) .. (307,125) .. controls (327,155) and (252,147) .. (232,117) .. controls (212,87) and (222,72) .. (242,62) -- cycle ;
\node[above] at (270,20)  {$\Gamma = \partial \Omega$};
\node[below] at (220,100)  {\textcolor{white}{$u$}};
\node[above] at (250,165)  {$\Omega_{\textup{out}}$};
\node[above] at (250,100)  {$\Omega_{\textup{in}}$};
\end{tikzpicture}}
\caption{Geometry of the interface identification problem of \cref{subsec:interface}. The fixed boundary $\Gamma$ is the outer boundary of the body, shown in black, and carries the applied current. The shape $u$, shown in white, is an interior interface separating the inner region $\Omega_{\textup{in}}$, in blue, from the outer region $\Omega_{\textup{out}}$, in orange, on which the conductivity $\kappa$ takes its two values. Here the moving portion of $\partial\holdall$ is empty, so that $\Gamma=\partial\holdall$.}
\label{fig:geometry_impedance}
\end{figure}

Let $\bar{y}\colon \holdall\to\R$ denote the above-mentioned measurements and $\nu>0$ be a given constant.  The outward normal vector to $\holdall$ and the outward normal vector to $\Omega_{\text{in}}$ are both denoted by $\bm{n}$ in the following.
Then, the mathematical model of the interface identification problem is given by 
\begin{align}
\min_{u \in M} &\quad  \frac{1}{2} \int_{\holdall} (y(x)  - \bar{y}(x))^2 \, \textup{d}x+\nu \int_{u}  \ds \label{eq:problem} \\
\text{s.t.} & \quad - \nabla \cdot (\kappa \nabla y) =  0,\quad \text{in } \holdall,  \label{eq:PDE1}\\
& \quad \qquad\quad  \kappa \frac{\partial y}{\partial \bm{n}} = g, \quad \text{on } \Gamma. \label{eq:PDE2}
\end{align}

Note that~\eqref{eq:PDE1}--\eqref{eq:PDE2} is a pure Neumann problem, so we seek $y$ in $H^1_{\textup{av}}(\holdall) = \{  v \in H^1(\holdall) | \int_\holdall v \, \textup{d}x = 0\}$, which fixes the additive constant through the zero-average condition. 
In the implementation this constraint is enforced weakly by a scalar Lagrange multiplier, which also supplies the constant volume correction that makes the problem solvable for any $g\in H^{-1/2}(\Gamma)$, without a compatibility condition; the adjoint equation~\eqref{eq:adjoint-weak} is treated in the same way.

Thus, the boundary input function is $g\colon\Gamma  \rightarrow \R$ and the material coefficient $\kappa$ is a function $ \holdall  \rightarrow \R$, such that $ x \mapsto \kin\mathbbm{1}_{\Omega_{\textup{in}}}(x)+\kout\mathbbm{1}_{\Omega_{\textup{out}}}(x)$, where $\kappa_i\colon \Omega_i \rightarrow \R$ are constants and $\mathbbm{1}_{\Omega_i}$ denotes the indicator function of the set $\Omega_i$, for $i\in\{\text{in}, \text{out}\}$.

The objective function contains two terms. The first one is a tracking-type  functional, where the model is fitted to data measurements $\bar{y}$; we will denote this term as $J^{\textup{trk}}(u)$. The second term is a perimeter regularization (denoted by $J^{\textup{reg}}(u)$) often required for well-posedness; see for instance \cite[Section 1.1]{Sokolowski1992introduction}.
Moreover, we consider the continuity conditions for the state and flux at the interface $\left\llbracket \kappa \frac{\partial y}{\partial \bm{n}} \right\rrbracket = 0$, and $\llbracket y \rrbracket  = 0$, on $u$. The jump symbol $\left\llbracket\cdot\right\rrbracket$ is defined on the interface $u$ by $\llbracket y \rrbracket := y_{\text{in}} - y_{\text{out}}$, where  $y_{ \text{in}} := \operatorname{tr}_{\text{in}}(y|_{\Omega_{\textup{in}}})$ and $y_{ \text{out}} := \tr_{\text{out}}(y |_{\Omega_{\textup{out}}})$, and $\tr_{\text{in}}\colon \Omega_{\textup{in}} \rightarrow u$, $\tr_{\text{out}}\colon \Omega_{\textup{out}} \rightarrow u$ are trace operators.

Following the same arguments as in \cite{Geiersbach2021stochastic}, the shape derivative and shape differentiability of \eqref{eq:problem}--\eqref{eq:PDE2} can be obtained.
The weak formulation of the shape derivative of the tracking-term as
\begin{equation}
\label{shape_derivative_J_obj}
\begin{split}
D J^{\textup{trk}}(u)[\bm{W}] &= \dfrac12\int_{\holdall}  \divv(\bm{W}) (y-\bar{y})^2 \, \textup{d}x - \int_{\holdall}(y-\bar{y})\nabla\bar{y}\cdot \bm{W} \, \textup{d}x\\ 
 &\quad\quad + \int_{\holdall} \kappa \nabla y ^\top (\divv(\bm{W})\id - \nabla \bm{W} - \nabla \bm{W}^{\top})\nabla p \, \textup{d}x, \\
\end{split}
\end{equation}  
where $y\in \hH := \{  v \in H^1(\holdall) | \int_\holdall v \, \textup{d}x = 0\}$ is the weak solution of~\eqref{eq:PDE1}--\eqref{eq:PDE2} and $p\in \hH$ solves the adjoint equation 
\begin{equation}\label{eq:adjoint-weak}
\int_{\holdall} \kappa \nabla \varphi \cdot \nabla p\, \textup{d}x = - \int_{\holdall} (y-\bar{y})\varphi \, \textup{d}x \quad  \forall \varphi \in \hH.
\end{equation} 

On the other hand, the shape derivative of the regularization term $J^{\text{reg}}$ can be found for example in \cite{Novruzi2002structure} and is given by $\textup{d}J^\text{reg}(u)[\bm{W}] =\nu\int_{u}\zeta\left<\bm{W},n\right> ds$ with $\zeta$ denoting the mean curvature of $u$.
Finally, the shape derivative of $J$ is given by 
\begin{equation}
   \label{eq:shape_derivative}
    D J(u)[\bm{W}] = D J^{\text{trk}}(u)[\bm{W}] 
    + 
   D J^\text{reg}(u)[\bm{W}].
\end{equation}

Finally, it is worth mentioning that the volume term in~\eqref{eq:shape_derivative} is bounded on $H_{\Gamma,0}^1(\holdall,\R^2)$: since $\kappa$ is piecewise constant with respect to the partition $\Omega_{\textup{in}}\sqcup u\sqcup\Omega_{\textup{out}}$ and the interface $u$ is smooth, the state and adjoint equations have piecewise $H^2$ solutions, so that $\nabla y,\nabla p\in L^4(\holdall)$ and the tensor field multiplying $\nabla\bm W$ lies in $L^2(\holdall,\R^{2\times2})$. 
The perimeter term $DJ^{\textup{reg}}(u)[\bm W] =\nu\int_u\zeta\langle\bm W,\bm n\rangle\,\ds$ is bounded on the same space for a different reason, namely the boundedness of the trace operator $H^1(\holdall)\to L^2(u)$ together with $\zeta\in L^2(u)$. 
Hence \cref{assmp:dj-bounded} is satisfied.

\subsection{Two-dimensional optimal bridge}
\label{subsec:2d-bridge}
In this section, we consider a compliance minimization problem. As usual, we are interested in minimizing the compliance of the system while keeping the area of the structure as low as possible, in other words, we will solve the same experiment as in \cite{Allaire2014shape, Herzog2024discretize}.

Compliance minimization is classically approached by topology optimization, and the optimal designs are known to involve microstructure. 
That is not the problem we solve here. 
The method developed here deforms a shape by diffeomorphisms, so the topology of the initial configuration is preserved throughout, and changes of topology are outside its scope. 
We therefore assume that a topology optimization stage has been carried out beforehand and has fixed the number and approximate position of the holes, and we use the present method to refine their geometry. 
This sequential practice is often used in structural optimization~\cite{Dev2022sequential,houta2026sequential}, where topology optimization supplies the configuration and shape optimization refines the geometry.
Additionally, it makes the example a natural test case for our purposes: the deformations required are large, the mesh is placed under real strain, and the metric is what distinguishes the variants.

We describe the geometry before stating the problem. 
The hold-all domain $\holdall$ is the material occupied by the structure; that is, a bounded domain whose boundary decomposes as $\partial\holdall=\Gamma\sqcup\Sigma$ in the sense of~\cref{subsubsec:shape-grad}. 
The fixed portion $\Gamma$ is the outer contour of the structure, itself partitioned into a portion $\Gamma_D$ where the structure is clamped, a portion $\Gamma_N$ where a traction $\bm g$ is applied, and a remaining traction-free portion. 
The moving portion $\Sigma$ consists of the boundaries of finitely many disjoint holes in the material, and it is these that constitute the shape,
\[
u=\bigcup_{i=1}^{N} u_i ,
\]
each $u_i$ being a simple closed curve, with the $u_i$ pairwise disjoint and enclosed by $\Gamma$. In the terms of~\cref{rem:compSup}, $u$ is a multi-shape, identified with the image of $N$ non-intersecting circles under an element of $\DiffcOm$. 
The multi-shape optimization framework stablished in~\cite[Sec.~2]{Geiersbach2022pde} allow us to consider this generalization. 
The optimization therefore moves the holes within the structure, while the outer contour, the clamped portion and the applied load remain fixed. 
We refer the reader to \cref{fig:geometry_compliance} for an example of such a configuration. 

The particular configuration used in our computations, the shape of the outer contour, the number and initial position of the holes, and the location of $\Gamma_D$ and $\Gamma_N$, is specified in~\cref{subsec:compliance}.

\begin{figure}[ht]
    \centering
    \resizebox{!}{0.3\textheight}{%
    \begin{tikzpicture}
    \draw[thick] (0,0) -- (0,1) -- (2.5,4) -- (5,5) -- (7.5,4) -- (10,1) -- (10,0) -- (9,0) --  (5.5,0) -- (4.5,0) -- (1,0) --cycle;
    \draw[line width = 4pt, color = NavyBlue] (0,0)--(1,0);
    \draw[line width = 4pt, color = NavyBlue] (9,0)--(10,0);
    \draw[->, thick] (4.5,0) -- (4.5,-0.5);
    \draw[->, thick] (5,0) -- (5,-0.5);
    \draw[->, thick] (5.5,0) -- (5.5,-0.5);
    \draw[line width = 4pt, color = OliveGreen] (4.5,0)--(5.5,0);
    \node at (5,-0.7)  {$\bm{g}$};
    \draw[thick, orange] (2.5,1) circle (0.5);
    \draw[thick , orange] (3.5,3) circle (0.5);
    \draw[thick , orange] (6.5,3) circle (0.5);
    \draw[thick , orange] (7.5,1) circle (0.5);
   \node[right] at (9,2.5) {$\Gamma$};	
   \node[below,orange] at (5,2) {$\Sigma = u_1 \sqcup u_2 \sqcup u_3 \sqcup u_4$};
   \node[below, orange] at (7.5,0.5) {$u_2$};
   \node[above, orange] at (6.5,3.5) {$u_3$};
   \node[below, orange] at (2.5,0.5) {$u_1$};
   \node[above, orange] at (3.5,3.5) {$u_4$};
    \node[below] at (0.5,0)  {$\Gamma^D$};
    \node[below] at (9.5,0)  {$\Gamma^D$};
    \node[above] at (5,0)  {$\Gamma^N$};
    \end{tikzpicture}}
\caption{Geometry of the compliance problem of \cref{subsec:2d-bridge}. The fixed boundary $\Gamma$ is the outer contour of the structure, shown in blue on the clamped portion $\Gamma^D$, in green on the loaded portion $\Gamma^N$, and in black where it is traction-free. The moving boundary $\Sigma$, shown in orange, consists of the hole boundaries and constitutes the shape $u$.}
    \label{fig:geometry_compliance}
\end{figure}

Mathematically this problem can be formulated as 
\begin{align}
    \label{prob:compliance}
    &\min_{u\in M} \quad  J(u) = \hat{J}(u,\bm{y}) \coloneqq \int_{\holdall} \bm{f}\cdot \bm{y} \, \textup{d} x + \int_{\Gamma_N} \bm{g} \cdot \bm{y} \, \textup{d} s +
    \ell \int_{\holdall} \textup{d}x \\
    &\text{s.t. } \text{for all }  \bm{V}\in H^1_{\Gamma_D}(\Omega) \times H^1_{\Gamma_D}(\Omega) \text{ the following equation holds:} 
    \\ & \qquad  2\mu \int_{\holdall} \varepsilon(\bm{y})\colon \varepsilon(\bm{V}) \, \textup{d} x 
    + \lambda \int_{\holdall} \operatorname{tr}(\varepsilon(\bm{y})) \operatorname{tr}(\varepsilon(\bm{V})) \, \textup{d} x 
    = 
    \int_{\holdall} \bm{f}\cdot \bm{V} \, \textup{d} x + \int_{\Gamma_N} \bm{g} \cdot \bm{V} \, \textup{d} s
\end{align}
where $\bm{f}\in L^2(\R^2)\times L^2(\R^2)$ and $\bm{g}\in H^1(\R^2)\times H^1(\R^2)$ are given volume and boundary loads, respectively. The subspace $H^1_{\Gamma_D}(\Omega)$ contains the functions which are zero on $\Gamma_D$. 

The computation of the shape derivative has been carried out following the variational approach suggested in \cite{Ito-Kunisch-Peichl} and its volume formulation is given by the following expression
\begin{align}
\nonumber  DJ(u)[\bm{W}]  & = 
  2 \, \int_{\holdall} (\nabla\bm{f} \bm{y}) \cdot \bm{W} + (\bm{f}\cdot \bm{y})\operatorname{div}(\bm{W}) \, \textup{d}x 
  \\
  \label{eq:shape_der_compliance}
  &  - \int_{\holdall} \mathcal{H}\varepsilon(\bm{y})\colon \varepsilon(\bm{y})\operatorname{div}(\bm{W}) \, \textup{d} x + 2 \, \int_{\holdall} \mathcal{H}\operatorname{sym}(\nabla\bm{y}\nabla \bm{W})\colon \varepsilon(\bm{y}) \, \textup{d} x \\
\nonumber  & + 2 \, \int_{\Gamma_N}(\bm{g}\cdot \bm{y}) (\operatorname{div}(\bm{W}) - \bm{n}^\top \operatorname{sym}(\nabla \bm{W})\bm{n}) \, \textup{d} s
  + \ell \int_{\holdall} \operatorname{div}(\bm{W}) \, \textup{d} x,
\end{align}
where the linearized strain tensor is given by $\varepsilon(\bm{W}) = (0.5)(\nabla \bm{W} + \nabla\bm{W}^\top)$, the matrix operator $\mathcal{H}$ maps $\R^{2\times 2}$ into $\R^{2\times 2}$ matrices such that $\mathcal{H} G = 2\mu G + \lambda \operatorname{tr}(G) I$, with $I$ the identity matrix. The symmetrization of a matrix will be denoted by $\operatorname{sym}(G) = (0.5)( G + G^\top)$, and $\mu,\lambda$ are the so-called Lam\'e parameters. The outward normal vector to $\holdall$ is denoted by $\bm{n}$, and $A \colon B$ denotes the Frobenius inner product between the matrices $A$, $B$. 
The shape differentiability of this problem is proven using~\cite[Thm.~1]{Allaire2014shape}, as this case is a specific instance of their more general result.

\Cref{assmp:dj-bounded} requires more care here than in \cref{subsec:interface}.
The boundary term over $\Gamma_N$ in~\eqref{eq:shape_der_compliance} causes no difficulty, since for deformation fields with zero trace on the outer boundary one has $\nabla\bm W=\partial_n\bm W\otimes\bm n$ there, so that $\operatorname{div}(\bm W)=\bm n^\top\operatorname{sym}(\nabla\bm W)\bm n$ and the integrand vanishes.
The remaining terms are quadratic in $\varepsilon(\bm y)$, and pairing them with $\nabla\bm W\in L^2(\holdall)$ requires $\nabla\bm y\in L^4(\holdall)$, which is a question of corner regularity.
The transitions between $\Gamma_D$ and the adjacent traction-free portions occur in the interior of a straight edge; for such a flat angle the singular exponents of the Lam\'e system have real part $\operatorname{Re}\alpha=k+\nicefrac12$, $k\in\Z$, and the solution lies in $W^{2,p}$ near the transition point precisely for $p<\nicefrac43$~\cite[Rem.~6.5]{grisvard1989singularites}.
By $W^{2,p}(\holdall)\hookrightarrow W^{1,q}(\holdall)$ with $\nicefrac1q=\nicefrac1p-\nicefrac12$ this gives $\nabla\bm y\in L^q(\holdall)$ for every $q<4$, but not for $q=4$.
At the two endpoints of $\Gamma_N$ the condition remains of traction type on both sides and only the datum jumps, which produces a weaker singularity and leaves $\nabla\bm y\in L^q(\holdall)$ locally for every $q<\infty$.

\Cref{assmp:dj-bounded} therefore fails, though only just and only at the two Dirichlet--traction transitions: $DJ(u)$ is bounded on $W^{1,q}_{\Gamma,0}(\holdall,\R^2)$ for every $q<4$, and not on $H^1_{\Gamma,0}(\holdall,\R^2)$.
Fortunately, the splitting system does hold in this weaker setting.
It suffices to pose the stage $j=s-1$ in $W^{1,q}_{\Gamma,0}(\holdall,\R^2)$; the Lax--Milgram theorem is not available there, and what replaces it is the isomorphism property of $\id-A\Delta$ between $W^{1,q'}_{\Gamma,0}(\holdall,\R^2)$ and the dual of $W^{1,q}_{\Gamma,0}(\holdall,\R^2)$, $q'=q/(q-1)$, which for mixed conditions on a Lipschitz domain holds for $q$ in a neighbourhood of $2$ whose size depends on the geometry, see e.g.,~\cite{meyers1963p}; whether $q<4$ lies in that range depends on the interior angles of $\Gamma$ and we do not verify it here.
Since $q'<2$, the resulting field lies in $W^{1,q'}_{\Gamma,0}(\holdall,\R^2)\hookrightarrow L^{2q'/(2-q')}(\holdall,\R^2)\hookrightarrow L^2(\holdall,\R^2)$ on the bounded domain $\holdall$, and the remaining stages are unchanged.

\section{Numerical investigations with inner and outer metrics}
\label{sec:NumericalEx}
This section is focused on describing the performance of the Riemannian steepest descent method introduced in~\cref{subsec:BasicsShapeCalc} while using different metrics to compute the shape gradients. 
In particular, we will focus on the $H^s$ metric on $\DiffcOm$ introduced and discussed in \cref{sec:ShapeOpt} and \cref{sec:modelProb}.  We will consider different values of the parameters $s\in \mathbb{N}$. Note that for the use of the Sobolev metrics $H^s$ in our numerical applications, the constant $A>0$ must be chosen. In this work, $A$ is chosen empirically in order to ensure good algorithmic performance. We will compare the results of this outer metric approach to results using inner metrics, in particular the SP metric introduced in \eqref{eq:SPmetric}.  In this paper, and as suggested in~\cite{Schulz2016efficient}, we set the bilinear form defining the SP metric to be the one associated to the variational formulation of the linear elasticity problem. Regarding the computational aspects of this metric, we will again follow the recommendations in~\cite{Schulz2016computational} and set $\lambda = 0$, and the parameter $\mu$ will be the solution of a Poisson problem decreasing smoothly from the  values $\mu_{\textup{min}}$ to $\mu_{\textup{max}}$.
These parameters are chosen ad hoc for each experiment. 

We numerically solve the two problems described in \cref{sec:modelProb}: the problem inspired in the electrical impedance tomography, and the two-dimensional optimal bridge. 
In both experiments, very large deformations of the meshes  are expected. This will allow us to compare how the different metrics handle such deformations in terms of the mesh quality at the final meshes. 

As previously discussed, both problems will be solved using the Riemannian steepest descent method. However, we will employ different metrics to generate different variants of the Riemannian gradient. While the method remains the same in the broader sense, the use of different metrics creates distinct approaches. For clarity, we will use the term ``variant'' to differentiate among the resulting Riemannian steepest descent methods, each one characterized by its unique Riemannian metric.  
\subsection{Implementation details}
\label{subsec:implementation}			

In this section, we briefly describe various aspects used in our numerical investigations and the discretization of both, the shapes and state variables used in this paper. 

\subsubsection{Riemannian retraction}
\label{subsec:retraction}
One decisive aspect of~\cref{alg:gradAlg} is how one updates the iterates. Commonly, this is done by using the exponential map on a Riemannian manifold. However, this implies solving a highly non-linear second order systems of ODEs, the geodesic equation. For the inner SP metric, the associated geodesic equation has not been derived yet, and this is one advantage of using outer metrics instead. However, we have delayed the numerical approximation of the geodesic equation on this manifold for a later publication. In this paper, we will use instead the retraction $\mathcal{R}_u$ defined from $T_u M \to M$,  such that $\bm{W} \mapsto (\id + \bm{W})(u)$.

We call $\mathcal R_u$ a retraction because it satisfies the two defining properties, $\mathcal R_u(0)=u$ and $\dd\mathcal R_u(0)=\id_{T_uM}$, both of which hold by construction for any $s$; what depends on $s$ is whether $\mathcal R_u(-t\bm V)$ is an admissible shape.
For $s\ge2$ it is, by the following argument. 
\cref{prop:split-equivalence} gives $\bm V\in H^{2s-1}(\holdall,\R^2)$, and $2s-1>1+\nicefrac{d}{2}$ with $d=2$ precisely when $s\ge2$, so the Sobolev embedding $H^{2s-1}(\holdall,\R^2)\hookrightarrow C^1(\holdall,\R^2)$ gives $\bm V$ a bounded continuous gradient.
For $t\|\nabla\bm V\|_{L^\infty}<1$ the map $\id-t\bm V$ is a small Lipschitz perturbation of the identity, hence a bijection of $\R^2$ by the Banach fixed-point theorem and a $C^1$-diffeomorphism by the inverse function theorem, its Jacobian $I-t\,D\bm V$ being invertible under the same condition; restricted to $u$, this leaves $\mathcal R_u(-t\bm V)$ a smoothly embedded, non-self-intersecting curve for $t$ small enough.
For $s=1$ the gradient is only $H^1$ and no such conclusion follows.
The same applies to the SP variant, whose gradient solves a linear elasticity problem tested against $H^1_0(\holdall,\R^2)$ and therefore lies in $H^1(\holdall,\R^2)$.
Among the variants compared in \cref{sec:NumericalEx}, the $H^s$ metrics with $s\ge2$ are thus the only ones for which admissibility follows from the regularity of the gradient rather than being enforced at the mesh level.

We do not turn this into a theorem, since two points would need care. The regularity $\bm V\in H^{2s-1}(\holdall,\R^2)$ holds globally only under the boundary hypothesis of \cref{prop:split-equivalence}, which neither hold-all of \cref{sec:modelProb} satisfies, the fixed portion $\Gamma$ being polygonal in both cases; the conclusion survives up to the moving portion $\Sigma$, where the deformation is applied, but the fixed-point argument needs a bound on $\|\nabla\bm V\|_{L^\infty(\holdall)}$, and near the corners of $\Gamma$ that would have to be argued separately.
The threshold on $t$ would also have to be made explicit.
In practice admissibility is enforced numerically, by halving a trial step whenever it would drop the mesh-quality measure~\eqref{eq:mesh_quality} below a fixed floor; among the variants reported here only SP and $H^1$ required this intervention.

Therefore, the iterate update described in line~\ref{step:update} of the algorithm can be replaced by $\holdall(u^{k+1}) \coloneqq \{ x\in \holdall \colon x = x^k + t^k \bm{W}^{k} \}$,  where $x^k \in \holdall(u^k)$ and $\bm{W}^k$ solves~\eqref{eq:grad_system_reindexed} and $t^k$ is an appropriate step length. As already noted in~\cite{Geiersbach2021stochastic,Herzog2024discretize} and many other works, while using the perturbation of identity in a context of triangular meshes (meaning elements from the finite element discretization) the mesh destruction is almost unavoidable when the shape gradient is ``too large''. However, there are many ways to avoid this, one of them is by adapting the stepsize $t^k$, which is the approach we follow. 

Within the framework established in the previous sections,~\cref{alg:gradAlg} establishes different variants of the Riemannian steepest descent method, one per each value of $s\in \mathbb{N}$.
The shape derivative appearing in \cref{step:gradient} is given by either \eqref{eq:shape_derivative} or \eqref{eq:shape_der_compliance} depending on the problem under consideration.

\begin{algorithm2e}
       \caption{Riemannian steepest descent method on $(\DiffcOm, H^s)$}
        \label{alg:gradAlg_Diff}
	
    \DontPrintSemicolon

		\textbf{Require:} Objective function $J$, given by \eqref{eq:problem} and \eqref{prob:compliance} respectively, on $(\DiffcOm, H^s)$ \;
		
             \textbf{Input:} Initial shape $u^0\in \DiffcOm$, maximum number of iterations $N_{\textup{max}}$ \;
		        
        \While{a suitable termination criterion is not satisfied and $k\leq N_{\textup{max}}$}{

    Compute the deformation field $\bm{V}^k$ as the solution of~\eqref{eq:grad_system_reindexed}.\; \label{step:gradient}

    Choose a suitable stepsize $t^k$.\;
    
     Set $u^{k+1}\coloneqq \mathcal{R}_{u^k}(-t^k \bm{V}^k)$.\nllabel{step:update} \; 

     Set $ k \gets k + 1 $. \;
        }
	\end{algorithm2e}

\subsubsection{Discretization of the problem}
\label{subsubsec:FEM}
To numerically solve the problems, we use the finite element method (FEM) because of various reasons. 
First, triangulations allow explicit representation of the shapes, and therefore of the corresponding shape deformations. Moreover, when using FEM only one mesh is needed to solve the state and adjoint equations and the shape gradient's computation.

As noted in \cite{Berggren2010unified} using this approach one is restricted to use piecewise linear elements to obtain the equivalence between the continuous and discrete shape derivative. 
This is an additional reason to consider the reformulation of the shape gradient problem \eqref{eq:metric_whole_space} given in \eqref{eq:grad_system_reindexed} as opposed to using higher-order finite elements.
As usual in FEM, the domain is covered by a triangular mesh which will be denoted by $\holdall_h$. 
Notice that the creation of the triangular mesh on $\holdall$ should provide a discretization for the unknown shape $u$.  
We denote by $\mathcal{S}^1(\holdall_h)$ the finite element space of piecewise linear globally continuous functions defined over $\holdall_h$, and we will also consider particular subspaces of functions with zero Dirichlet boundary conditions $\mathcal{S}_{\Gamma_D}^1(\holdall_h)$, which we will specify in each application. 

The moving boundary is resolved as a piecewise linear facet chain of a conforming simplicial mesh; no isoparametric geometry representation is used, so the geometric error in representing the interface is $\mathcal{O}(h^2)$ and does not depend on $s$. 
Mesh sizes and element counts for each experiment are reported in \cref{subsec:Expe1,subsec:compliance}. No quadrature degree is prescribed: all forms are assembled with the automatic degree estimation of the FEniCS form compiler, which is exact here, since with piecewise linear fields the state, adjoint and shape-derivative forms have piecewise constant integrands and the mass term of each stage of~\eqref{eq:grad_system_reindexed} is quadratic. 
Note finally that \cref{assmp:dj-bounded} is a hypothesis on the continuous problem, in which the interface is smooth; the discrete interface is polygonal and converges to it as $h\to0$.
We note finally that \cref{assmp:dj-bounded} and \cref{prop:split-equivalence} are statements about the continuous problem, in which the moving boundary is smooth; the discrete moving boundary is a polygon, and converges to it as $h\to0$. 
The hypothesis $\partial\holdall\in C^{2s-1,1}$ is likewise a continuous one, and is in any case not satisfied by the hold-alls used here, as discussed after \cref{prop:split-equivalence}.

All of our experiments were implemented in the finite element software FEniCS with a CPU Apple M2 (8 cores) and 8.0 GB RAM.
The code and data supporting the findings of this study are available from the corresponding author upon request.

\subsection{Experiment 1: interface identification problem}
\label{subsec:Expe1}

\begin{figure}[h]
    \begin{center}
    \begin{subfigure}{0.31\textwidth}
    \centering
    \raisebox{5pt}{%
    \includegraphics[width=0.73\textwidth
    ]{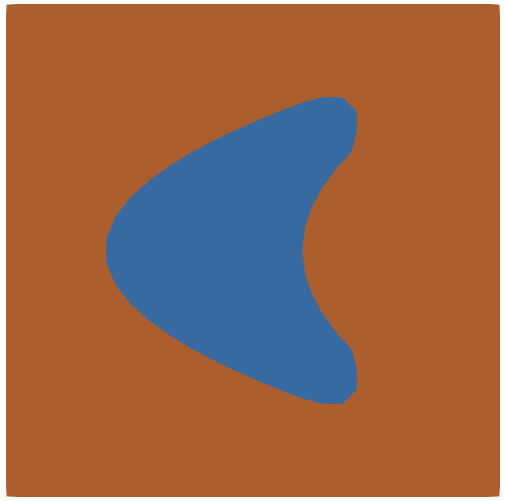}}
    \caption{Optimal solution}
    \label{fig:optimal_shape}
    \end{subfigure}
    ~
    \begin{subfigure}{0.315\textwidth}
    \centering
    \includegraphics[width=\textwidth]{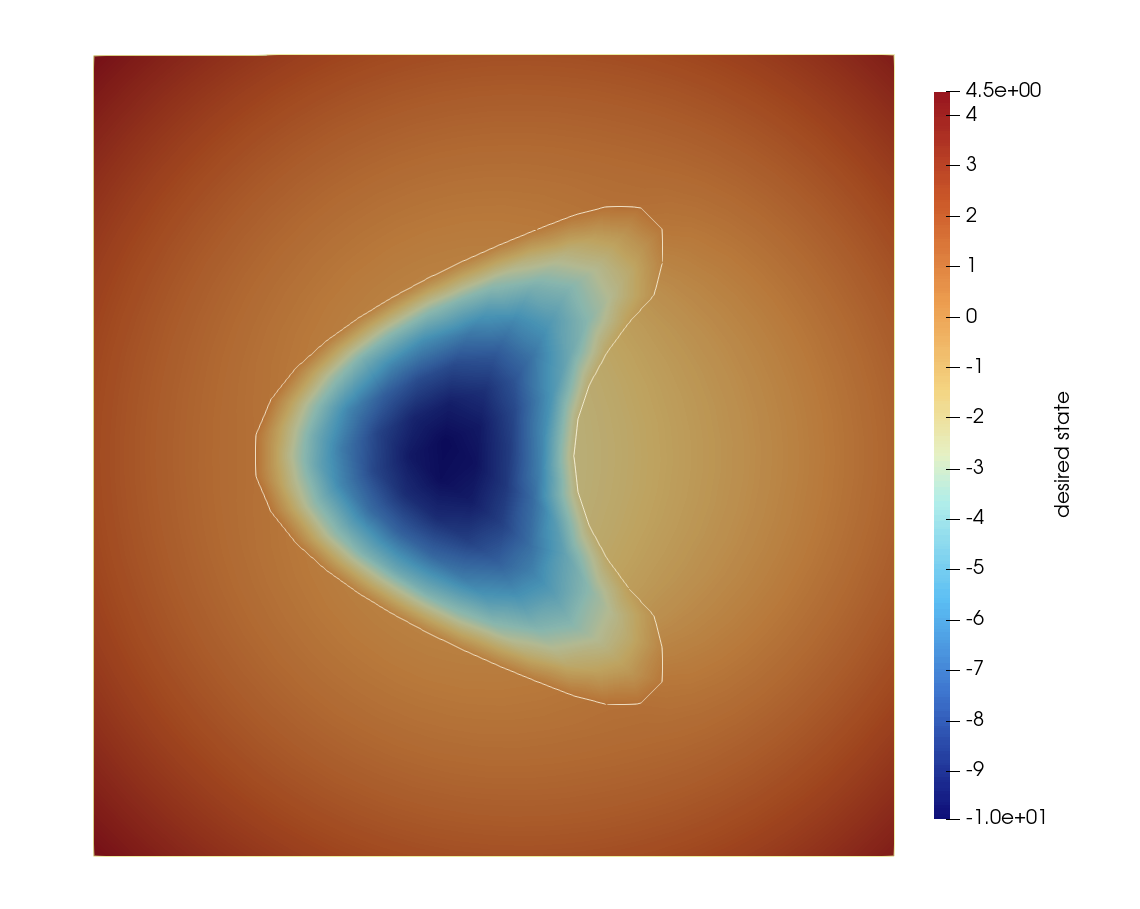}
    \caption{Target state $\bar{y}$}
    \label{fig:desired_state}
    \end{subfigure}
    ~
    \begin{subfigure}{0.31\textwidth}
    \centering
    \raisebox{5pt}{%
    \includegraphics[width=0.73\textwidth
    ]{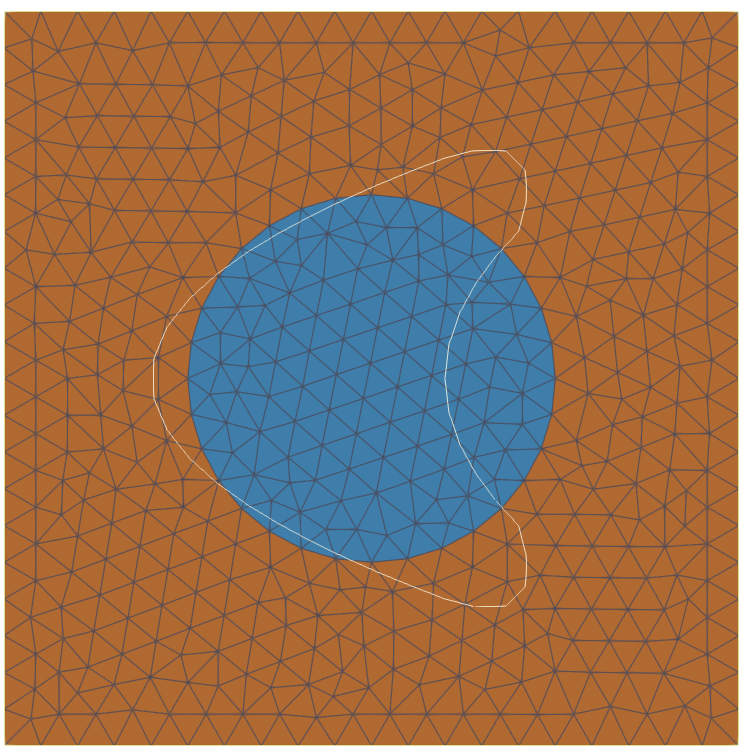}}
    \caption{Initial mesh $u^0$}
    \label{fig:initial_mesh}
    \end{subfigure}
    \end{center}
    \caption{Data used in the experiment described in \cref{subsec:Expe1}.}
\end{figure}

This experiment is inspired by~\cite{Geiersbach2021stochastic}, and aims to represent the identification of a human lung, where the target $\bar{y}$ has been computed using electrical impedance tomography. The hold-all domain is set to $\holdall = [-1, 0]\times [-0.5,0.5]$ with $\Gamma = \partial\holdall $ and $\Sigma = \emptyset$, the material coefficients are set to $\kappa_{\textup{out}} = 1$ and $\kappa_{\textup{in}} = 0.05$, the boundary input is a constant function equal to $10$ and the perimeter regularization parameter $\nu=0$. 
The target state $\bar{y}$ is generated using the optimal shape depicted in~\cref{fig:optimal_shape}. 
A key point is that the algorithm does not know the actual shape, the only information provided is the value of $\bar{y}$ shown in~\cref{fig:desired_state}.

Setting $\nu=0$ in~\eqref{eq:problem} switches off the perimeter regularisation, and is deliberate for this experiment.
The target state $\bar y$ is generated from an exactly recoverable interface using the same finite element model as the inversion, with no noise added, so that differences between the reconstructed interfaces are attributable to the choice of metric alone rather than to regularisation or measurement error.
It also keeps the shape derivative in volume form.
This should not be read as a claim that $\nu=0$ is generally admissible: the regularisation is in general essential for well-posedness of the reconstruction problem, see e.g.~\cite[Th.~2.9]{Sokolowski1992introduction}.

\Cref{fig:initial_mesh} shows the initial mesh, which consists of 544 nodes and 1118 elements.

As a first step, we verify numerically whether the shape gradients obtained with different metrics differ from one another. To this end, we compute the shape derivative associated with the initial shape depicted in \cref{fig:initial_mesh}. Note that this is the only iteration of the algorithm for which all shape derivatives coincide.

\begin{figure}[h]
    \noindent
\begin{tikzpicture}[trim left=(SP.west)]
\node[inner sep=0pt] (SP) at (0,0)
    {\includegraphics[width=.25\textwidth]{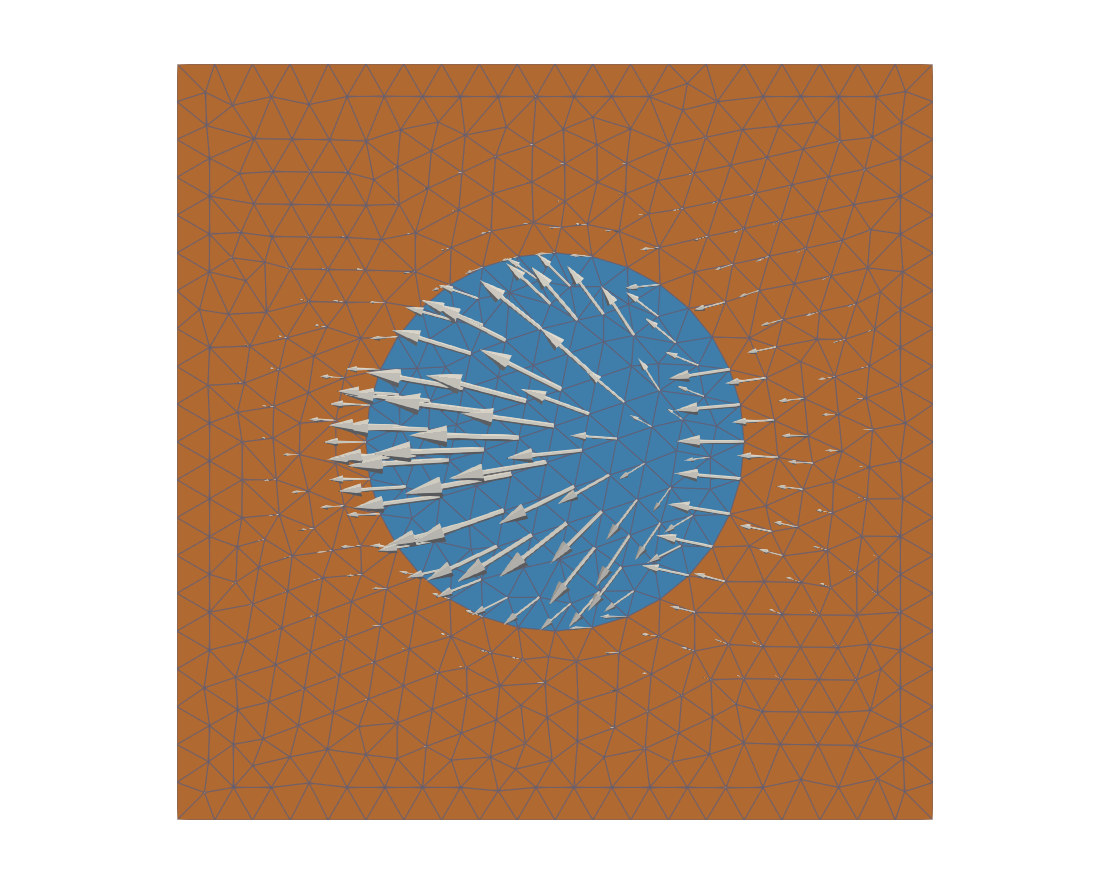}};
\node at (0,-2.1) {SP metric};
\node[inner sep=0pt] (H1) at (3.2,0)
    {\includegraphics[width=.25\textwidth]{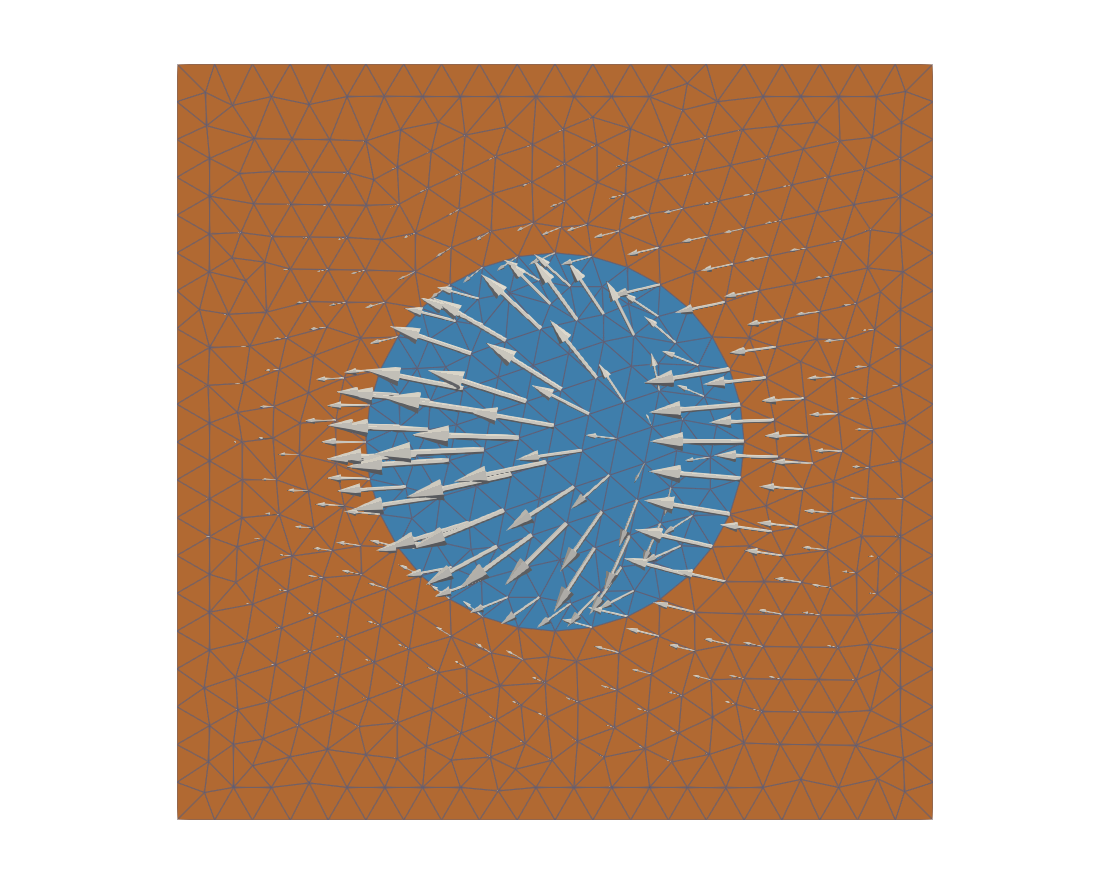}};
\node at (3.2,-2.1) {$H^1$ metric};
\node[inner sep=0pt] (H2) at (6.4,0)
    {\includegraphics[width=.25\textwidth]{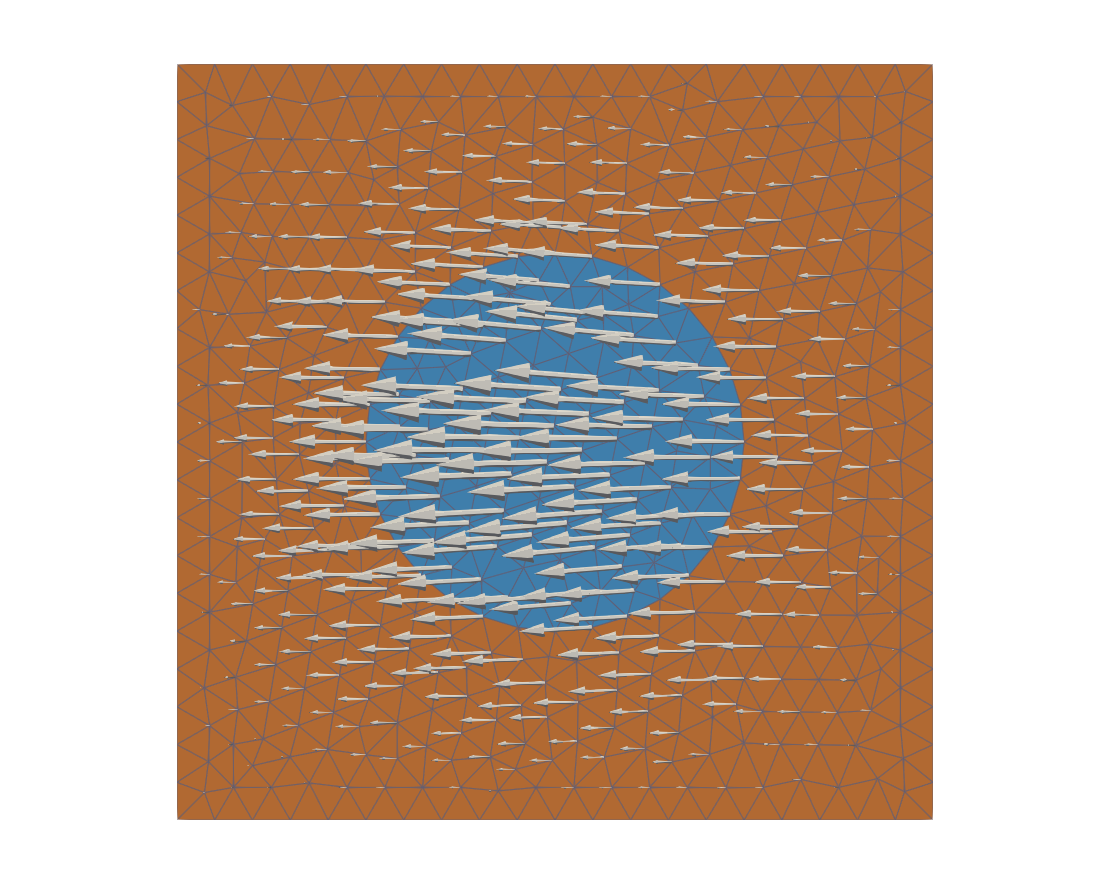}};
\node at (6.4,-2.1) {$H^2$ metric};
\node[inner sep=0pt] (H3) at (9.6,0)
    {\includegraphics[width=.25\textwidth]{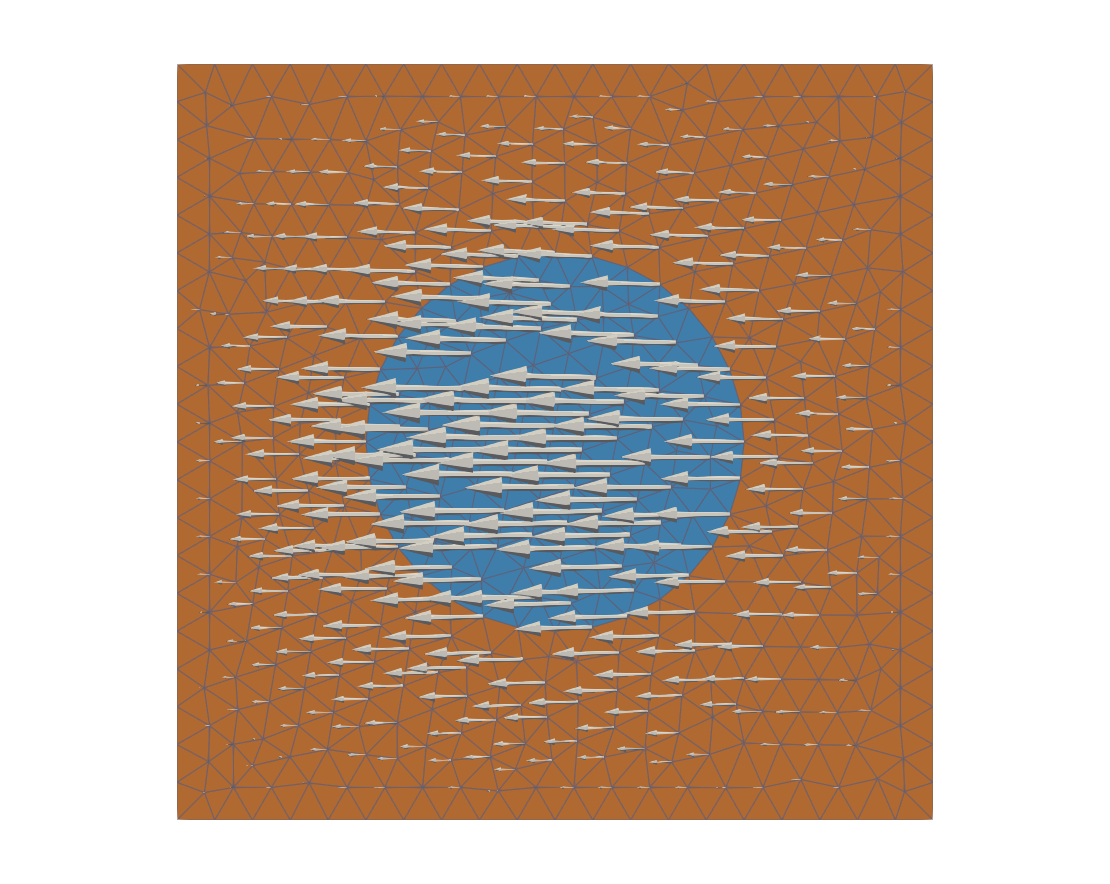}};
\node at (9.6,-2.1) {$H^3$ metric};
\node[inner sep=0pt] (H4) at (12.8,0)
    {\includegraphics[width=.25\textwidth]{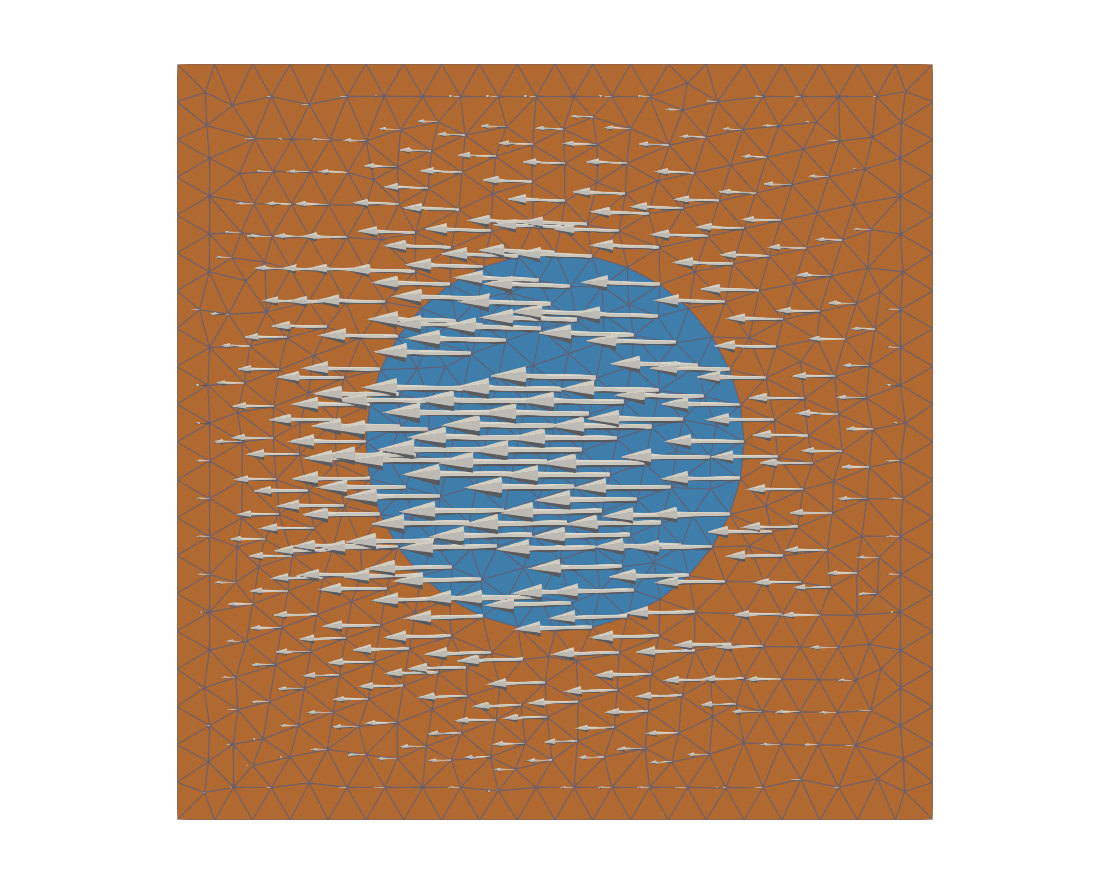}};
\node at (12.8,-2.1) {$H^4$ metric};
\end{tikzpicture}
    \caption{Different shape gradients from the experiment described in \cref{subsec:Expe1}.}
    \label{fig:exp_1_gradients}
\end{figure}

\Cref{fig:exp_1_gradients} displays the shape gradients, showing that they differ in both magnitude and direction. Consequently, each variant of the method produces different subsequent iterates. This can be understood as each variant exploring the energy landscape differently, which leads to different performance of the algorithm and even to find different local minima as we show in what follows.

Moreover, \Cref{tab:diff_gradient} presents the computation time required to obtain a single shape gradient for each metric.
Although computing gradients for the $H^s$ metrics requires solving a system of $s$ coupled elliptic PDEs, our results show that this does not significantly increase the computational burden. 
Interestingly, the time required to compute the SP gradient is comparable to that of the $H^3$ gradient. We attribute this partly to the SP metric's requirement of solving a Poisson equation to obtain more accurate values of the Lam\'e parameter $\mu$.
For each $H^s$ metric, we selected parameter $A$ values to optimize algorithm performance. For the computation of the shape gradient w.r.t. the SP metric we use the values $\mu_{\textup{min}} = 5$ and $\mu_{\textup{max}} = 20$ where $\mu$ decreases from $\mu_{\textup{min}}$ on $u$ to $\mu_{\textup{max}}$ on the outer boundary.

The stepsize $t^k$ used in the retraction $\mathcal{R}_{u^k}(-t^k\bm{V}^k)$ is fixed per metric variant at the nominal value reported in \cref{tab:experiment1}. At every iteration this nominal value is attempted first; if the resulting mesh would drop the quality measure~\eqref{eq:mesh_quality} below a fixed floor ($10^{-2}$), the step is halved and the trial repeated, up to a bounded number of times (eight) before the best available step is accepted regardless. 
Each iteration re-attempts the nominal stepsize from scratch rather than continuing from whatever reduced value the previous iteration ended on. 
We report that among all the variants reported in this experiment only the $SP$ and $H^1$ required this intervention.
Unlike the compliance problem of \cref{subsec:compliance}, 
Additionally, no remeshing was allowed in this experiment. 
The values of $A$ in \cref{tab:diff_gradient,tab:experiment1} were obtained by a coarse manual search over a small set of candidates for each $s$, retaining the value giving the fastest stable decrease of both the objective and the gradient norm without degrading the mesh-quality measure over a short trial run; $H^2$ and $H^3$ required minimal such tuning, whereas $H^4$ needed the most careful joint adjustment of $A$ and $t^k$.

\cref{fig:exp_1_gradients} shows the different gradients obtained for the same shape derivative. The vector fields have been scaled to improve the visibility, since we are interested in the direction of each gradient rather than their magnitude. To give the reader an idea of the actual magnitude of each gradient we present the $L^2$ norm in \cref{tab:diff_gradient}, and we specify each used scale factor.

We can observe in \cref{fig:exp_1_gradients} that the SP- and $H^1$ gradients look very similar; however, their magnitudes are completely different. The direct consequence of this behavior is that when using a retraction (perturbation of identity) smaller stepsizes will be required to avoid mesh destruction. 
On the other hand, the $H^s$ gradients with $s\geq 2$ look very similar among them but very different from the SP-gradient, showing that is definitely worth studying the performance of each version of the Riemannian steepest descent method with these metrics.

\pgfplotstableset{
  empty cells with={--},
  columns/metric/.style={column type={c},
    column name={metric}, string type,
    string replace*={H1}{$H^1$},
    string replace*={H2}{$H^2$},
    string replace*={H3}{$H^3$},
    string replace*={H4}{$H^4$},
  },
  columns/A/.style={column type={c},
    column name={$A$}, fixed, precision=4, dec sep align,
  },
  columns/scale_factor/.style={column type={c},
    column name={scale factor}, fixed, precision=4, dec sep align,
  },
  columns/norm_V_L2/.style={
    column name={$\|\bm{V}\|_{L^2}$}, fixed, precision=3, dec sep align,
  },
  columns/time_s/.style={column type={r},
    column name={time (s)}, fixed, precision=4, dec sep align,
    assign cell content/.code={%
				\pgfkeyssetvalue{/pgfplots/table/@cell content}%
				{\SI{##1}{\second}}%
		}},
}
\begin{table}[ht]
\centering
\makebox[\linewidth][c]{%
\begin{filecontents*}{Data/it_table_1.txt}
metric,A,scale_factor,norm_V_L2,time_s
SP,,1.0,5.1196e-02,0.0346
H1,0.0625,0.0085,6.7344e+00,0.0117
H2,0.5,0.65,1.0369e-01,0.0329
H3,0.2,0.65,1.0077e-01,0.0557
H4,0.05,0.1,7.8409e-01,0.0843
\end{filecontents*}
\pgfplotstabletypeset[
  col sep=comma,
  every head row/.style={before row=\toprule, after row=\midrule},
    every row no 0/.style={after row=\midrule},
  every last row/.style={after row=\bottomrule}]{Data/it_table_1.txt}}
\caption{Parameters used to illustrate the shape gradients from \cref{subsec:Expe1}.}
\label{tab:diff_gradient}
\end{table}

Now, we report the performance of the Riemannian steepest descent method when we use the different metrics. We measure the value of the objective function at each iteration, the norm of the gradients, and the mesh quality measured as
\begin{equation}
    \label{eq:mesh_quality}
\varphi(\holdall_h^k) = \min_{T \in \Sigma_h^k} 2\frac{\textup{inrad}(T)}{\textup{circumrad}(T)} 
\end{equation}
where $\Sigma_h^k$ denotes the triangulation associated to the mesh describing $\holdall_h^k$, $T$ denotes a triangle from the triangulation, and $\textup{inrad}$ and $\textup{circumrad}$ are the radius of the incircle and circumcircle of triangle $T$, respectively. This quality measure has range zero to one where close-to zero values indicate the presence of at least one degenerate element.
 
As stopping criterion, we use the $L^2$-norm of the gradient with a tolerance of $2\cdot 10^{-4}$. It is worth mentioning that in the Riemannian context the $L^2$ may not be the best way to measure the optimality of the solution. In fact, this should be measured using the appropriate Riemannian metric. 
However, since we are comparing the different variants of the method associated with different metrics, we have found this is the fairest possible way. 
It is worth mentioning that we set the maximum number of iterations to $500$. 
Failure to converge within 500 iterations does not imply that the method
would not converge later.
In fact, all the variants of the method seem to recover the solutions relatively well. The main difference relies on the mesh quality of the obtained solutions.

\Cref{tab:cost_experiment1} complements \cref{tab:diff_gradient} (which reports only the cost of a single gradient evaluation) with the total cost of each complete optimization run: the iteration count already given in \cref{tab:experiment1}, the total wall-clock time, and the total number of linear systems solved. Each iteration requires solving the state equation, the adjoint equation, and the shape-gradient system~\eqref{eq:grad_system_reindexed}; for the $H^s$ metrics the latter is assembled as a single, monolithically coupled mixed-space linear system regardless of $s$ (three solves per iteration in total), while the SP metric additionally solves a Poisson problem for the Lam\'e parameter $\mu$ (four solves per iteration in total). Finally, we note for transparency that the synthetic measurement $\bar{y}$ used in this experiment is generated by solving the same finite element model (the same conductivities $\kappa_{\textup{lung}}, \kappa_{\textup{ext}}$ and the same choice of finite elements) as is used throughout the inversion, albeit on an independently generated mesh with a different triangulation and a different number of vertices than the one the optimization is initialized on. No noise was added to $\bar{y}$. Using the same underlying model for data generation and inversion is a mild form of what is sometimes called an inverse crime; a more stringent test, generating $\bar{y}$ on a substantially finer and fully independent mesh and/or with additive noise, is left for future work.

\pgfplotstableset{
  columns/iterations/.style={column type={c},
    column name={iterations $\bar k$},
  },
  columns/total_wall_clock_s/.style={column type={r},
    column name={total time}, fixed, precision=3,
    assign cell content/.code={%
      \pgfkeyssetvalue{/pgfplots/table/@cell content}%
      {\SI{##1}{\second}}%
    }},
  columns/s_per_iter/.style={column type={r},
    column name={time / iteration}, fixed, precision=4,
    assign cell content/.code={%
      \pgfkeyssetvalue{/pgfplots/table/@cell content}%
      {\SI{##1}{\second}}%
    }},
  columns/total_pde_solves/.style={column type={c},
    column name={total PDE solves}, /pgf/number format/1000 sep={},
  },
}
\begin{table}[ht]
\centering
\makebox[\linewidth][c]{%
\begin{filecontents*}{Data/it_total_cost_table.txt}
metric,iterations,total_wall_clock_s,s_per_iter,total_pde_solves
SP,500,36.389,0.0728,2000
H1,500,39.514,0.0790,1500
H2,184,11.007,0.0598,552
H3,279,21.975,0.0788,837
H4,500,53.221,0.1064,1500
\end{filecontents*}
\pgfplotstabletypeset[
  col sep=comma,
  every head row/.style={before row=\toprule, after row=\midrule},
    every row no 0/.style={after row=\midrule},
  every last row/.style={after row=\bottomrule}]{Data/it_total_cost_table.txt}}
\caption{Total computational cost of the complete optimization run for each variant of the method for the experiment described in \cref{subsec:Expe1}, complementing the per-gradient timings of \cref{tab:diff_gradient}.}
\label{tab:cost_experiment1}
\end{table}

The summary of the experiment is collected in \cref{tab:experiment1} where we present the parameters used for each variant of the method, as well as the iterations required to converge together with the mesh quality measure and the $L^2$-norm of the shape gradient at the last iterate.

\pgfplotstableread[col sep=comma]{Data/it_table_2.txt}{\dataexperimentone}

\newcommand{\getbestvalue}[4]{
  \pgfplotstablesort[sort key=#2]{\sortedtabletmp}{#1}%
  \edef\getbestdirtmp{#4}%
  \def\getbestminwordtmp{min}%
  \ifx\getbestdirtmp\getbestminwordtmp
    \pgfplotstablegetelem{0}{#2}\of{\sortedtabletmp}%
  \else
    \pgfplotstablegetrowsof{\sortedtabletmp}%
    \pgfmathtruncatemacro{\lastrowtmp}{\pgfplotsretval-1}%
    \pgfplotstablegetelem{\lastrowtmp}{#2}\of{\sortedtabletmp}%
  \fi
  \global\expandafter\let\csname #3\endcsname\pgfplotsretval
}
\getbestvalue{\dataexperimentone}{k}{BestkExpOne}{min}
\getbestvalue{\dataexperimentone}{J}{BestJExpOne}{min}
\getbestvalue{\dataexperimentone}{norm_V_L2}{BestnormExpOne}{min}
\getbestvalue{\dataexperimentone}{mesh_quality}{BestmeshExpOne}{max}

\makeatletter
\newcommand{\boldmathifbest}[1]{%
  \pgfkeysgetvalue{/pgfplots/table/@preprocessed cell content}\rawtmp
  \edef\besttmp{#1}%
  \ifx\rawtmp\besttmp
    \pgfkeysgetvalue{/pgfplots/table/@cell content}\fmttmp
    \edef\newcontenttmp{{\noexpand\boldmath\fmttmp}}%
    \pgfkeyslet{/pgfplots/table/@cell content}\newcontenttmp
  \fi
}
\makeatother

\pgfplotstableset{
  columns/t/.style={column type={c},
    column name={$t$}, fixed, precision=4, dec sep align,
  },
  columns/k/.style={column type={c},
    column name={$\bar{k}$},
    postproc cell content/.append code={\boldmathifbest{\BestkExpOne}},
  },
  columns/J/.style={column type={c},
    column name={$J(u^{\bar{k}})$}, sci, sci precision=3,
    /pgf/number format/sci generic={mantissa sep=\times, exponent={10^{##1}}},
    postproc cell content/.append code={\boldmathifbest{\BestJExpOne}},
  },
  columns/norm_V_L2/.style={column type={c},
    column name={$\|\bm{V}^{\bar{k}}\|_{L^2}$}, sci, sci precision=3,
    /pgf/number format/sci generic={mantissa sep=\times, exponent={10^{##1}}},
    postproc cell content/.append code={\boldmathifbest{\BestnormExpOne}},
  },
  columns/mesh_quality/.style={column type={c},
    column name={$\varphi(\holdall_h^{\bar{k}})$}, sci, sci precision=3,
    /pgf/number format/sci generic={mantissa sep=\times, exponent={10^{##1}}},
    postproc cell content/.append code={\boldmathifbest{\BestmeshExpOne}},
  },
}
\begin{table}[ht]
\centering
\makebox[\linewidth][c]{%
\pgfplotstabletypeset[
  every head row/.style={before row=\toprule, after row=\midrule},
  every row no 0/.style={after row=\midrule},
  every last row/.style={after row=\bottomrule}]{\dataexperimentone}}
\caption{Summary of the experiment described in \cref{subsec:Expe1}.}
\label{tab:experiment1}
\end{table}

\begin{figure}[h]
\begin{tikzpicture}[trim left=(SP.west)]
\node[inner sep=0pt] (SP) at (0,0)
    {\includegraphics[width=.25\textwidth]{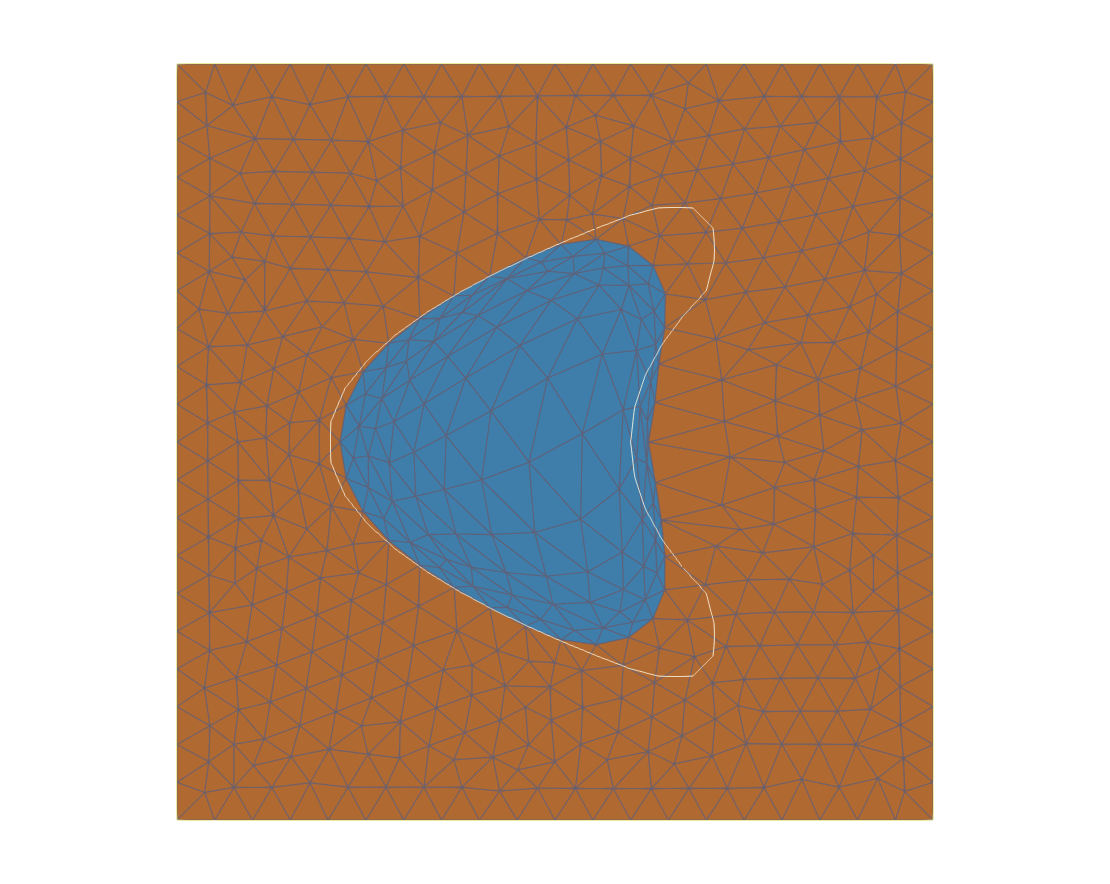}};
\node at (0,-1.7) {$u^{500}$};
\node at (0,-2.1) {SP metric};
\node[inner sep=0pt] (H1) at (3.2,0)
    {\includegraphics[width=.25\textwidth]{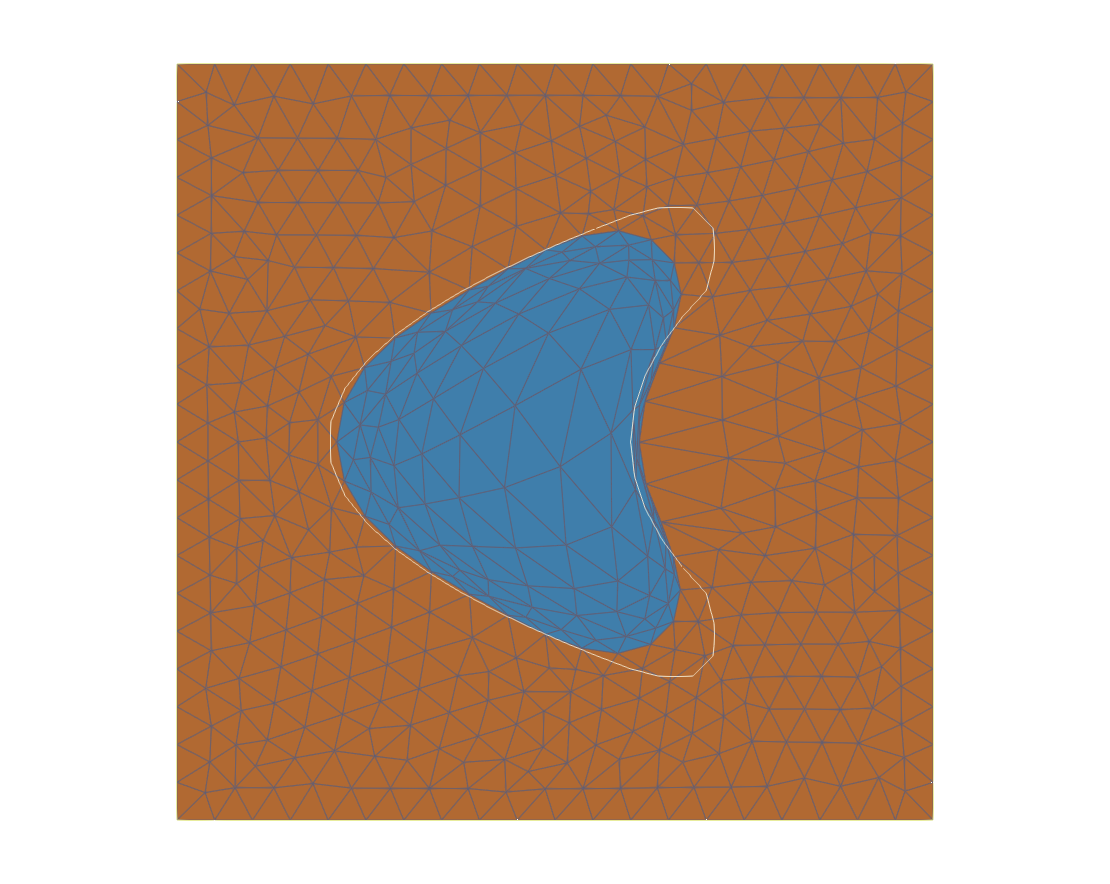}};
\node at (3.2,-1.7) {$u^{500}$};
\node at (3.2,-2.1) {$H^1$ metric};
\node[inner sep=0pt] (H2) at (6.4,0)
    {\includegraphics[width=.25\textwidth]{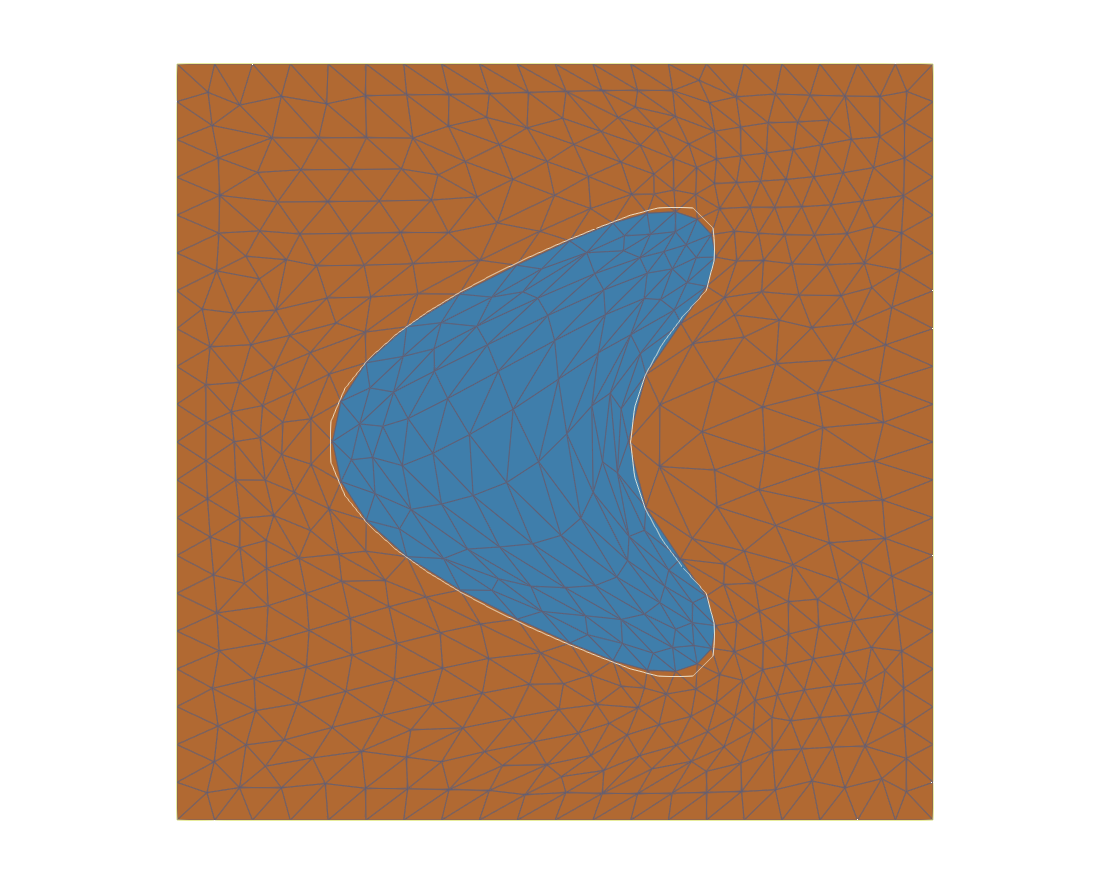}};
\node at (6.4,-1.7) {$u^{184}$};
\node at (6.4,-2.1) {$H^2$ metric};
\node[inner sep=0pt] (H3) at (9.6,0)
    {\includegraphics[width=.25\textwidth]{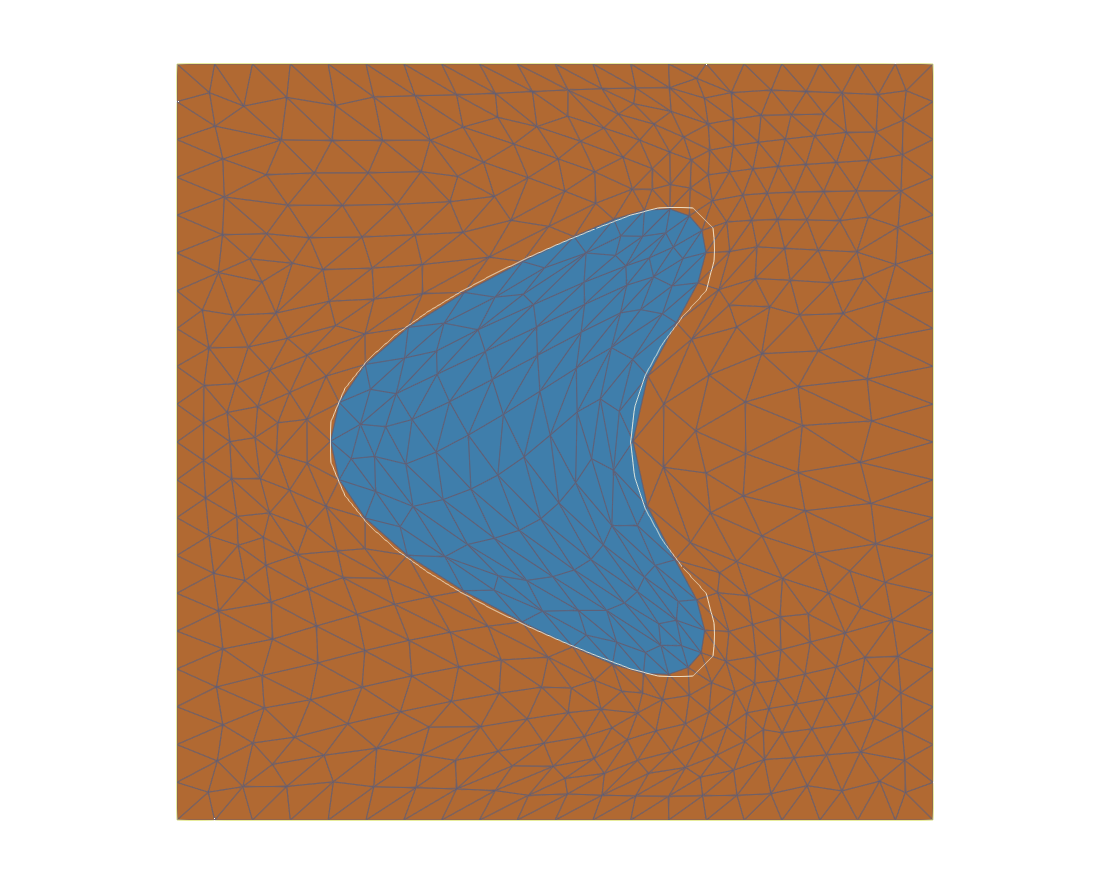}};
\node at (9.6,-1.7) {$u^{279}$};
\node at (9.6,-2.1) {$H^3$ metric};
\node[inner sep=0pt] (H4) at (12.8,0)
    {\includegraphics[width=.25\textwidth]{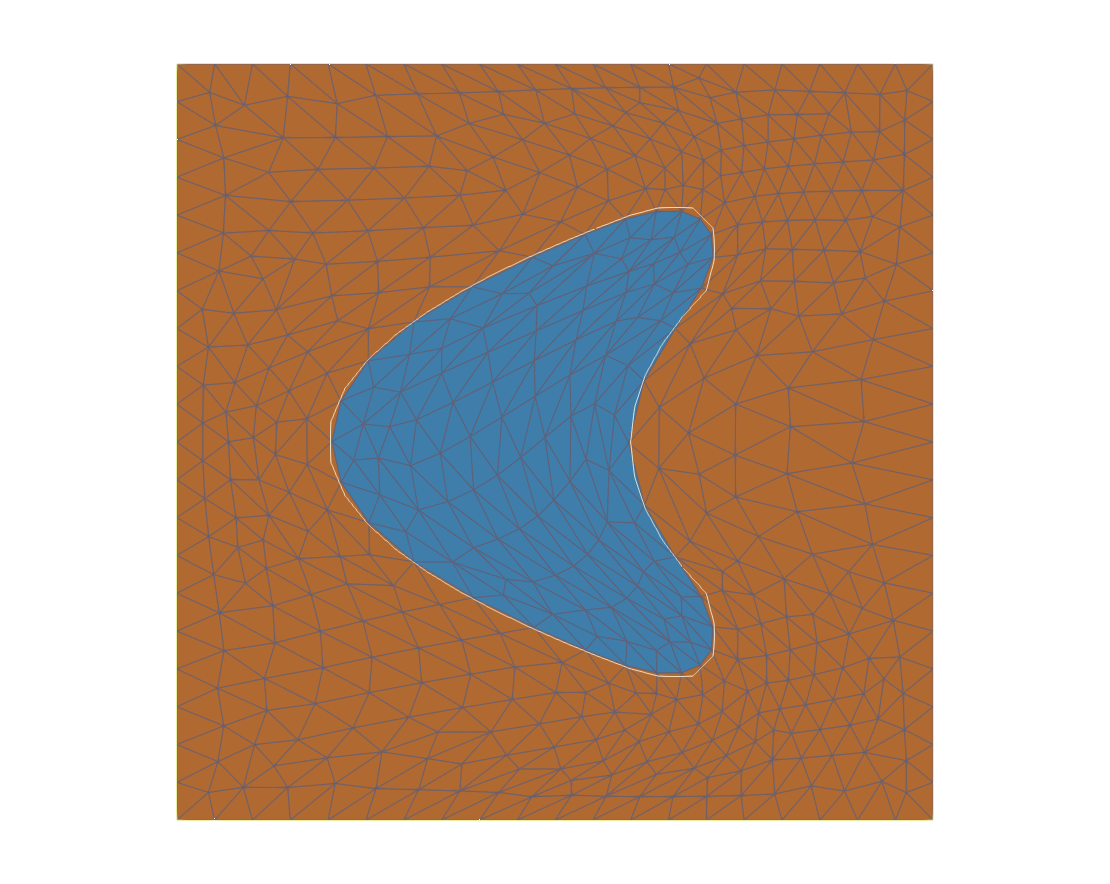}};
\node at (12.8,-1.7) {$u^{500}$};
\node at (12.8,-2.1) {$H^4$ metric};
\end{tikzpicture}
\caption{Final iterates of each variant of the method for the experiment described in \cref{subsec:Expe1}}
\label{fig:exp1_solutions}
\end{figure}

\cref{fig:exp1_solutions} shows the different solutions obtained by each one of the variants of the method. 
In white we show the ``optimal solution'', i.e., the one used to generate the desired stated and the one we are using as a ground truth. 
As already mentioned, all the variants of the method generate solutions very close to the ground truth. 
We show the mesh to showcase how the outer metrics generate solutions with overall better quality of the meshes. 

In addition to showing the final shapes generated by each variant of the method, we present in \cref{fig:exp1_norm_qual} the changes on the objective function, $L^2$-norm of the gradients and mesh quality measure along the iterations generated by each variant of the method.

\begin{figure}[h]
\pgfplotsset{/pgfplots/group/.cd, horizontal sep=2cm}
\begin{center}
\begin{tikzpicture}
\pgfplotsset{compat=1.14, width=\textwidth}
\begin{groupplot}[group style={
        {group size=3 by 1}, xlabels at=edge bottom}, height=4.3cm, width=5.1cm,
        legend style={transpose legend, legend columns=5, draw=none }, 
        xlabel=iterations]
\nextgroupplot[title= $J(u^k)$, legend to name=grouplegend]
\addplot[color=black, dashed] table [x=iter,y=objective, col sep=comma]{Data/it_history_SP.txt};\label{pgfplots:SP}
\addlegendentry{SP metric}
\addplot[color=blue, dotted] table [x=iter,y=objective, col sep=comma]{Data/it_history_H1.txt};\label{pgfplots:H1}
\addlegendentry{$H^1$ metric}
\addplot[color=magenta, mark=*, mark size = 0.25pt] table [x=iter,y=objective, col sep=comma]{Data/it_history_H2.txt};\label{pgfplots:H2}
\addlegendentry{$H^2$ metric}
\addplot[color=violet, mark=x, mark size = 0.25pt] table [x=iter,y=objective, col sep=comma]{Data/it_history_H3.txt};\label{pgfplots:H3}
\addlegendentry{$H^3$ metric}
\addplot[color=orange, mark =star, mark size = 0.25pt] table [x=iter,y=objective, col sep=comma]{Data/it_history_H4.txt};\label{pgfplots:H4}
\addlegendentry{$H^4$ metric}

\nextgroupplot[title = $\|\bm{V}(u^k)\|_{L^2}$, ymode=log]
\addplot[color=black, dashed] table [x=iter,y=norm_gradient, col sep=comma]{Data/it_history_SP.txt};
\addplot[color=blue, dotted] table [x=iter,y=norm_gradient, col sep=comma]{Data/it_history_H1.txt};
\addplot[color=magenta, mark=*, mark size = 0.25pt] table [x=iter,y=norm_gradient, col sep=comma]{Data/it_history_H2.txt};
\addplot[color=violet, mark=x, mark size = 0.25pt] table [x=iter,y=norm_gradient, col sep=comma]{Data/it_history_H3.txt};
\addplot[color=orange, mark =star, mark size = 0.25pt] table [x=iter,y=norm_gradient, col sep=comma]{Data/it_history_H4.txt};

\nextgroupplot[title = $\varphi(\holdall_h^k)$, ymode=log]
\addplot[color=black, dashed] table [x=iter,y=mesh_quality_min, col sep=comma]{Data/it_history_SP.txt};
\addplot[color=blue, dotted] table [x=iter,y=mesh_quality_min, col sep=comma]{Data/it_history_H1.txt};
\addplot[color=magenta, mark=*, mark size = 0.25pt] table [x=iter,y=mesh_quality_min, col sep=comma]{Data/it_history_H2.txt};
\addplot[color=violet, mark=x, mark size = 0.25pt] table [x=iter,y=mesh_quality_min, col sep=comma]{Data/it_history_H3.txt};
\addplot[color=orange, mark =star, mark size = 0.25pt] table [x=iter,y=mesh_quality_min, col sep=comma]{Data/it_history_H4.txt};
\end{groupplot}

\node[yshift=-35pt] at ($(group c2r1.south)!0.5!(group c2r1.south)$) {\ref{pgfplots:SP} SP metric  \ref{pgfplots:H1} $H^1$ metric  \ref{pgfplots:H2} $H^2$ metric  \ref{pgfplots:H3} $H^3$ metric  \ref{pgfplots:H4} $H^4$ metric};
\end{tikzpicture}

\caption{Behavior of different quantities along the iterations for each variant of the method for the experiment described in \cref{subsec:Expe1}.}
\label{fig:exp1_norm_qual}
\end{center}
\end{figure}

As already mentioned in~\cite{Herzog2024discretize} many previously proposed Riemannian metrics in the context of optimize-then-discretize approach fail on generating small enough shape gradients, and for this reason, most of the gradient-based methods rely on early stopping. This behavior can also be observed in the middle panel of \cref{fig:exp1_norm_qual} where the $L^2$ norm of the shape gradient  remains bounded away from zero for the SP and $H^1$ metrics. On the other hand, the norm of the shape gradients for the $H^s$ metrics, with $s\geq 2$,  tends to zero. Particularly interesting is the behavior of the method with the $H^2$ metric since its $L^2$ norm converges faster to zero than the other variants of the method.
What is more surprising is the quality of the meshes generated by the $H^s$ metrics, with $s\geq 2$, as can be seen in the right panel of \cref{fig:exp1_norm_qual}. The values of the mesh quality measure stay remarkably stable all along the optimization process even though the meshes are subject to the same kind of large deformations as in the other two cases. Finally, it is also worth highlighting that the value of the objective function obtained by the $H^s$ metrics, with $s\geq 2$, is overall smaller than with the SP metric and $H^1$ metric. 

In view of the results, we conclude that the $H^2$ metric provides the best solution in terms of its objective function value with a relatively good mesh quality in the least amount of iterations. 
It is also worth recalling that the time of execution per iteration for the computation of this shape gradient is comparable to the one required for the SP metric, thus the method is overall computationally competitive.  

Finally, we would like to report that minimum parameter tuning was required to obtain the solutions for the $H^2,H^3$ norms. However, the variant of the method using the $H^4$ norm required a very thorough tuning of the parameter $A$ and the stepsize. For this reason, we believe the best solution is obtained with $s=2$, since it is a good trade off among all the indicators of good performance of the method. 

\subsection{Experiment 2: two-dimensional optimal bridge}
\label{subsec:compliance}
We assume there are no body forces applied to the structure, i.e., $\bm{f} = (0,0)^\top$, and the boundary loads on the inhomogeneous Neumann portion of the boundary are given by $\bm{g}=(0,-0.25)^\top$.
The Lam\'e parameters are given in terms of the Young's modulus $E=1$ and the Poisson ratio $\nu =0.3$ such that  $\mu = \nicefrac{E}{2(1+\nu)}$ and  $\lambda = \nicefrac{E\nu}{(1+\nu)(1-2\nu)}$.
The parameter penalizing the volume of the structure is set to $\ell = 9.9\times 10^{-2}$. 
These parameters have been chosen such that we can reproduce the results obtained in \cite{Herzog2024discretize}. As already mentioned, we assume we have an informed topology of the structure depicted in \cref{fig:exp2_initial_top}.

\begin{figure}[ht]
    \centering
    \resizebox{!}{0.3\textheight}{%
    \begin{tikzpicture}
    \draw[thick] (0,0) -- (0,1) -- (2.5,4) -- (5,5) -- (7.5,4) -- (10,1) -- (10,0) -- (9,0) --  (5.5,0) -- (4.5,0) -- (1,0) --cycle;
    \draw[line width = 4pt, color = NavyBlue] (0,0)--(1,0);
    \draw[line width = 4pt, color = NavyBlue] (9,0)--(10,0);
    \draw[->, thick] (4.5,0) -- (4.5,-0.5);
    \draw[->, thick] (5,0) -- (5,-0.5);
    \draw[->, thick] (5.5,0) -- (5.5,-0.5);
    \draw[line width = 4pt, color = OliveGreen] (4.5,0)--(5.5,0);
    \node at (5,-0.7)  {$\bm{g} = (0,-0.25)^\top$};
    \draw[thick, orange] (2.5,1) circle (0.5);
    \draw[thick , orange] (3.5,3) circle (0.5);
    \draw[thick , orange] (6.5,3) circle (0.5);
    \draw[thick , orange] (7.5,1) circle (0.5);
    \node[left] at (0,0)  {$(0,0)$};
    \node[left] at (0,1)  {$(0,1)$};
    \node[left] at (2.5,4)  {$(2.5,4)$};
    \node[above] at (5,5)  {$(5,5)$};
    \node[right] at (7.5,4)  {$(7.5,4)$};
    \node[right] at (10,1)  {$(10,1)$};
    \node[right] at (10,0)  {$(10,0)$};
    \node[below] at (0.5,0)  {$\Gamma^D$};
    \node[below] at (9.5,0)  {$\Gamma^D$};
    \node[above] at (5,0)  {$\Gamma^N$};
    \end{tikzpicture}}
    \caption{Informed initial shape for the experiment described in \cref{subsec:compliance}. The portion of the boundary in green corresponds to the inhomogeneous Neumann boundary condition. The blue portion of the boundary corresponds to the homogeneous Dirichlet condition. The moving boundary describing the actual shape of the bridge is depicted in orange.}
    \label{fig:exp2_initial_top}
\end{figure}
The initial mesh contains $2083$ nodes and $4172$ elements. 
It is well-known for this kind of problems that keeping the same initial mesh imposes restrictions on the attainable structures. For this reason, we allow remeshing within the optimization method when the mesh quality is under $10^{-1}$. 
Owing to the existence of multiple solutions and the expected micro structures, we aim to compare the methods based on the attained values of the objective function for the given topology. We do not allow changes of topology since this contradicts the theory developed in the previous sections. Moreover, due to these observations, we cannot use the norm of the gradient as stopping criterion. Instead, we use
\[
\max_{m = 2,\ldots, 10} J(u^{k-m}) - J(u^k) < 10^{-5}
\]
as suggested in~\cite{Herzog2024discretize}, i.e., the algorithm will stop when the value of the objective has not considerably changed in the last $10$ iterations. 
As per \cref{alg:gradAlg_Diff}, the maximum number of iterations for this experiment is set to $1000$, and we recall that failure to converge within the maximum number of iterations does not imply the method would not converge later.
Moreover, we have chosen the parameters $\mu_{\textup{min}} = 5$ and $\mu_{\textup{max}} = 15$ for the SP metric variant of the method to guarantee there will be no changes on the topology. However, it is worth highlighting that in this case we use the opposite rule as in \cref{subsec:Expe1}, i.e, we impose $\mu_{\textup{min}}$ on all of the fixed boundary pieces ($\Gamma_{\textup{out}}$, $\Gamma_{D,\textup{left/right}}$, $\Gamma_N$, $\Gamma_{\textup{out,left/right}}$) and $\mu_{\textup{max}}$ on the interface $u$.

We solved the problem comparing four different variants of the Riemannian steepest descent method, i.e., we use the SP metric and the $H^2$-, $H^3$- and $H^4$ metrics. The values of the parameter $A$ for each $s$ metric are shown in \cref{tab:exp_2_data}. For all the experiments we choose $t=1$ as the initial stepsize, however, it has been updated depending on the amount of times we needed to remesh (making it smaller each time we remesh). This table also contains the final iteration count, values of the objective function, $L^2$-norm of the gradient, and mesh quality measure~\eqref{eq:mesh_quality} at the last iteration.

Remeshing is carried out with Gmsh~\cite{geuzaine2009gmsh} whenever the mesh-quality measure~\eqref{eq:mesh_quality} drops below the threshold $10^{-1}$: the deformed mesh is exported, a fresh triangulation of the same domain is regenerated from the original geometry description, and the boundary and subdomain markers are reassigned on the new mesh before the optimization resumes. 
The stepsize update is precise, the initial value $t^0=1$ is left unchanged until the first remeshing event and, after the $j$-th such event, the stepsize is set to $t^0/j$. 
The $SP$, $H^3$ and $H^4$ variants required exactly three remeshing events and ended with stepsize $t=0.333$; while the $H^2$ variant required only two remeshing events with final stepsize $t=0.5$.
Because remeshing resets the triangulation, the mesh-quality value reported at the final iterate of \cref{tab:exp_2_data} reflects the deformation accumulated since the last remeshing event rather than the smoothness of the deformation fields over the whole run; the comparison of final mesh quality across variants that was informative in \cref{subsec:interface}, where no remeshing occurred, therefore carries less weight here.
Unlike Experiment~1, where a value of $A$ for each $H^s$ variant was selected by a coarse grid search over several candidates, this compliance problem required no comparable search: the values of $A$ reported in \cref{tab:exp_2_data} were reached after a handful of manual trials guided by the heuristic already described above (smaller $s$ requires larger $A$), and none of the four variants required further re-tuning once a stable, topology-preserving descent was observed.

\pgfplotstableread[col sep=comma]{Data/com_table_3.txt}{\dataexptwo}

\newcommand{\getbestvaluetwo}[4]{
  \pgfplotstablesort[sort key=#2]{\sortedtabletmptwo}{#1}%
  \edef\getbestdirtmptwo{#4}%
  \def\getbestminwordtmptwo{min}%
  \ifx\getbestdirtmptwo\getbestminwordtmptwo
    \pgfplotstablegetelem{0}{#2}\of{\sortedtabletmptwo}%
  \else
    \pgfplotstablegetrowsof{\sortedtabletmptwo}%
    \pgfmathtruncatemacro{\lastrowtmptwo}{\pgfplotsretval-1}%
    \pgfplotstablegetelem{\lastrowtmptwo}{#2}\of{\sortedtabletmptwo}%
  \fi
  \global\expandafter\let\csname #3\endcsname\pgfplotsretval
}
\getbestvaluetwo{\dataexptwo}{k}{BestkExpTwo}{min}
\getbestvaluetwo{\dataexptwo}{J}{BestJExpTwo}{min}
\getbestvaluetwo{\dataexptwo}{norm_V_L2}{BestnormExpTwo}{min}
\getbestvaluetwo{\dataexptwo}{mesh_quality}{BestmeshExpTwo}{max}

\makeatletter
\newcommand{\boldmathifbesttwo}[1]{%
  \pgfkeysgetvalue{/pgfplots/table/@preprocessed cell content}\rawtmptwo
  \edef\besttmptwo{#1}%
  \ifx\rawtmptwo\besttmptwo
    \pgfkeysgetvalue{/pgfplots/table/@cell content}\fmttmptwo
    \edef\newcontenttmptwo{{\noexpand\boldmath\fmttmptwo}}%
    \pgfkeyslet{/pgfplots/table/@cell content}\newcontenttmptwo
  \fi
}
\makeatother

\pgfplotstableset{
  columns/k/.style={column type={c},
    column name={$\bar{k}$},
    postproc cell content/.append code={\boldmathifbesttwo{\BestkExpTwo}},
  },
  columns/J/.style={column type={c},
    column name={$J(u^{\bar{k}})$}, fixed, precision=3, dec sep align,
    postproc cell content/.append code={\boldmathifbesttwo{\BestJExpTwo}},
  },
  columns/norm_V_L2/.style={column type={c},
    column name={$\|\bm{V}^{\bar{k}}\|_{L^2}$}, sci, sci precision=3,
    /pgf/number format/sci generic={mantissa sep=\times, exponent={10^{##1}}},
    postproc cell content/.append code={\boldmathifbesttwo{\BestnormExpTwo}},
  },
  columns/mesh_quality/.style={column type={c},
    column name={$\varphi(\holdall_h^{\bar{k}})$}, sci, sci precision=3,
    /pgf/number format/sci generic={mantissa sep=\times, exponent={10^{##1}}},
    postproc cell content/.append code={\boldmathifbesttwo{\BestmeshExpTwo}},
  },
}
\begin{table}
\centering
\makebox[\linewidth][c]{%
\pgfplotstabletypeset[
  every head row/.style={before row=\toprule, after row=\midrule},
  every row no 0/.style={after row=\midrule},
  every last row/.style={after row=\bottomrule}]{\dataexptwo}}
\caption{Summary of the experiment described in \cref{subsec:compliance}.}
\label{tab:exp_2_data}
\end{table}%
\begin{figure}[h] 
\centering
\begin{tikzpicture}
\node[inner sep=0pt] (SP) at (-5.2,0)
    {\includegraphics[width=.245\textwidth]{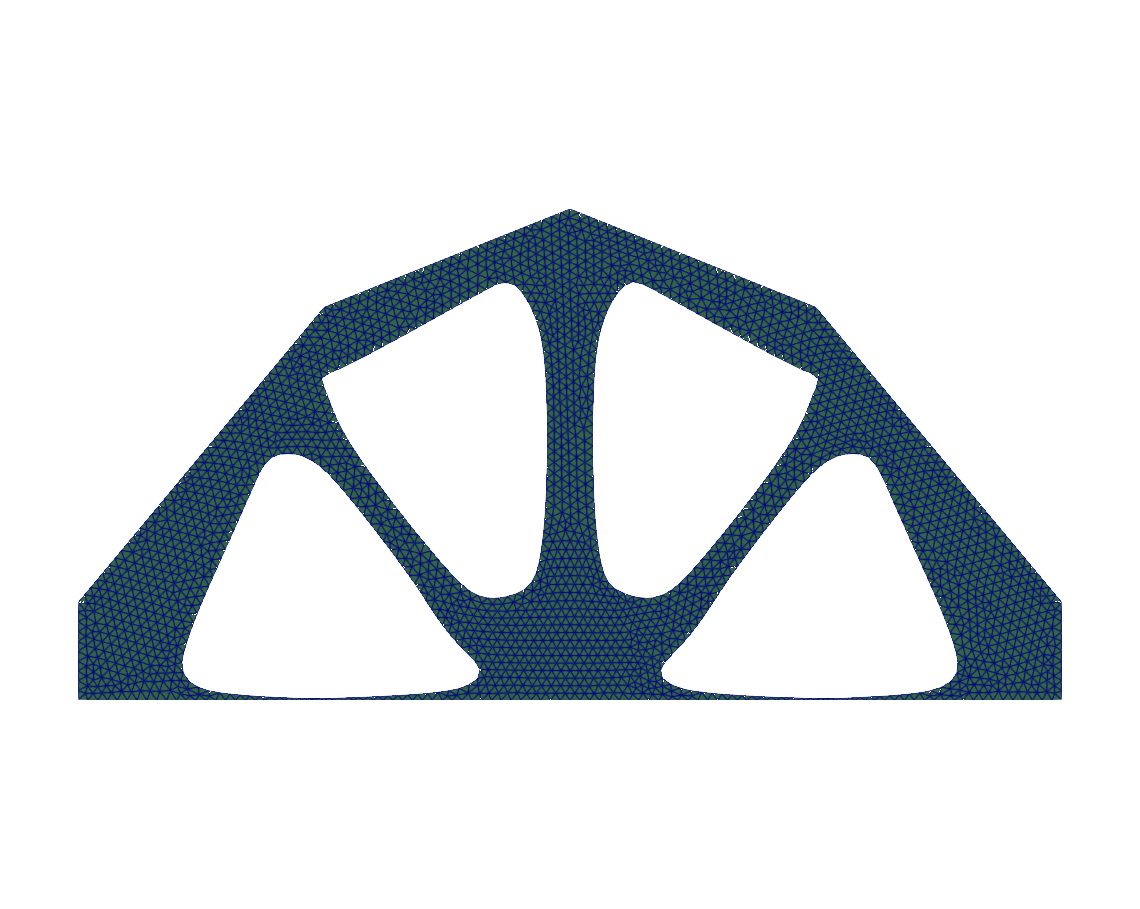}};
\node at (-5.2,-1.2) {$u^{1000}$};
\node at (-5.2,-1.6) {SP metric};
\node[inner sep=0pt] (H1) at (-1.1,0)
    {\includegraphics[width=.245\textwidth]{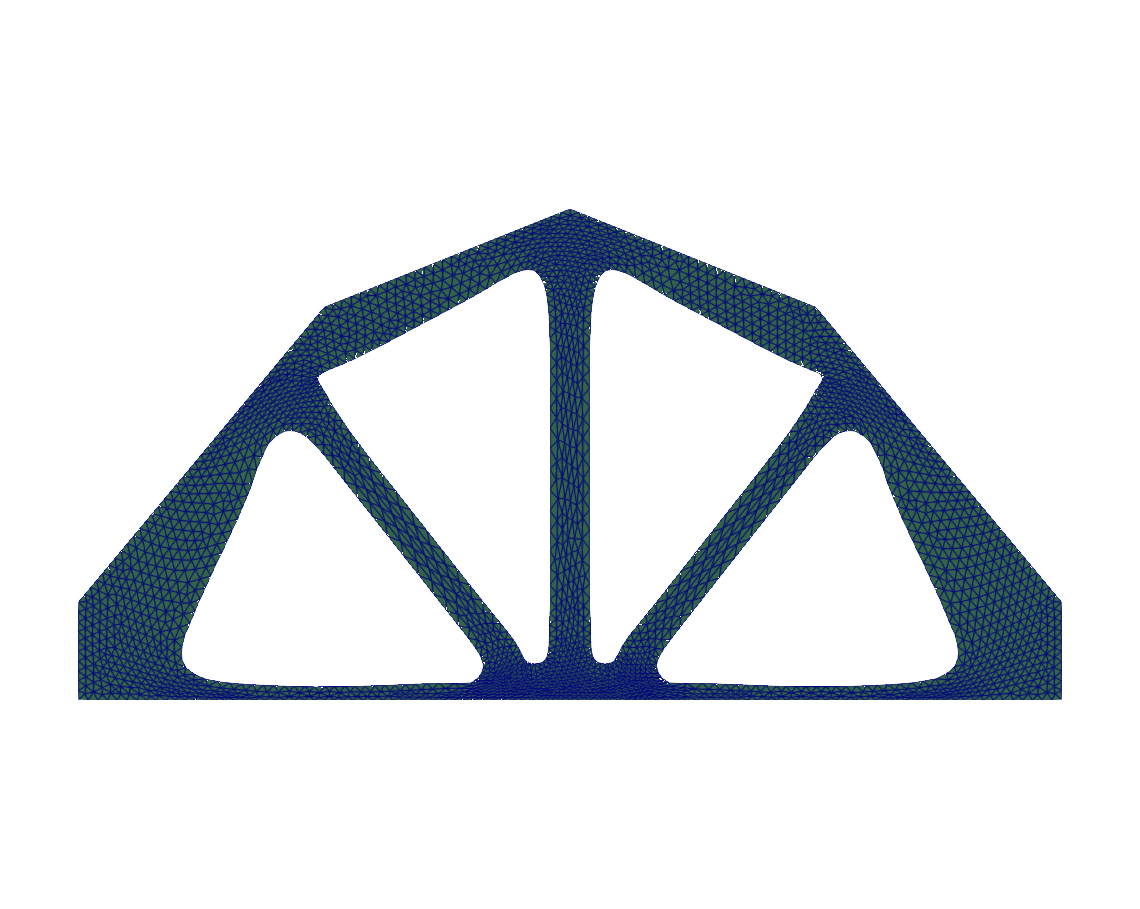}};
\node at (-1.1,-1.2) {$u^{1000}$};
\node at (-1.1,-1.6) {$H^2$ metric};
\node[inner sep=0pt] (H2) at (3,0)
    {\includegraphics[width=.245\textwidth]{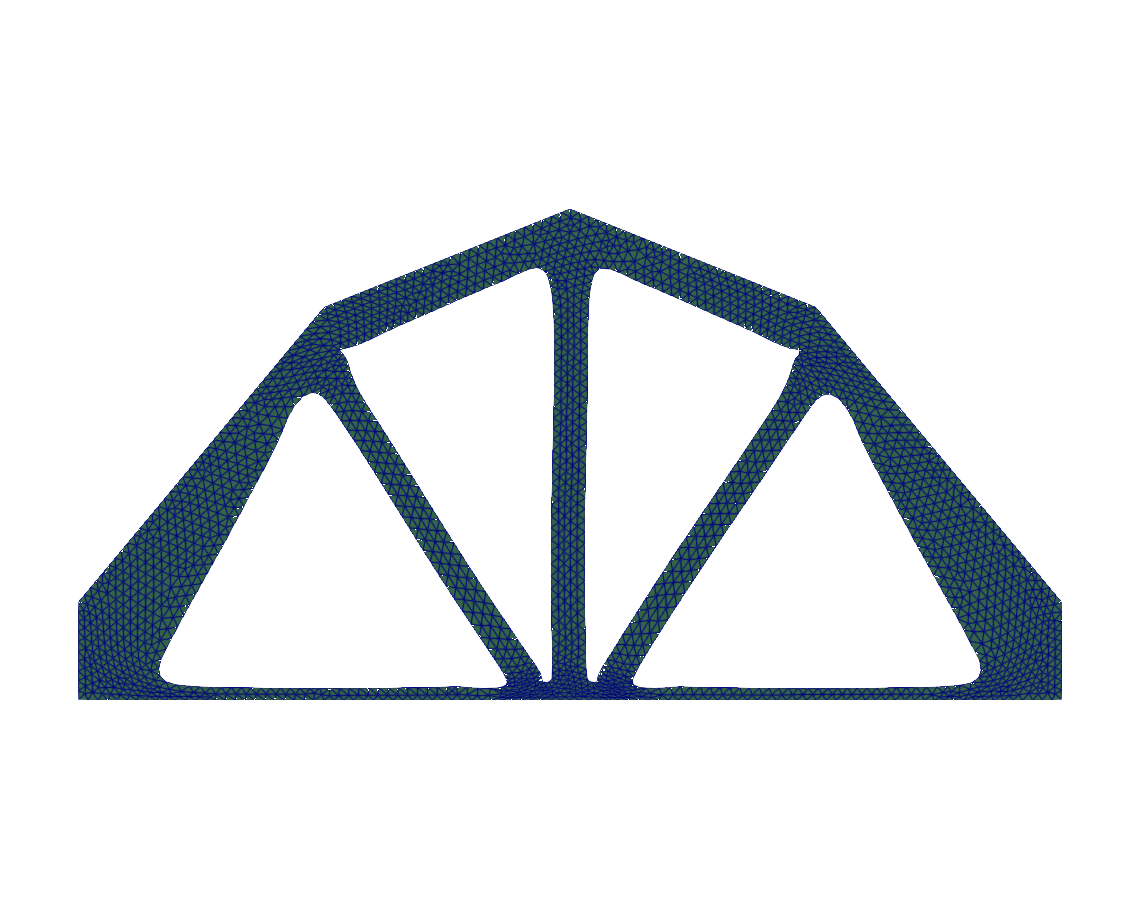}};
\node at (3,-1.2) {$u^{468}$};
\node at (3,-1.6) {$H^3$ metric};
\node[inner sep=0pt] (H3) at (7.1,0)
    {\includegraphics[width=.245\textwidth]{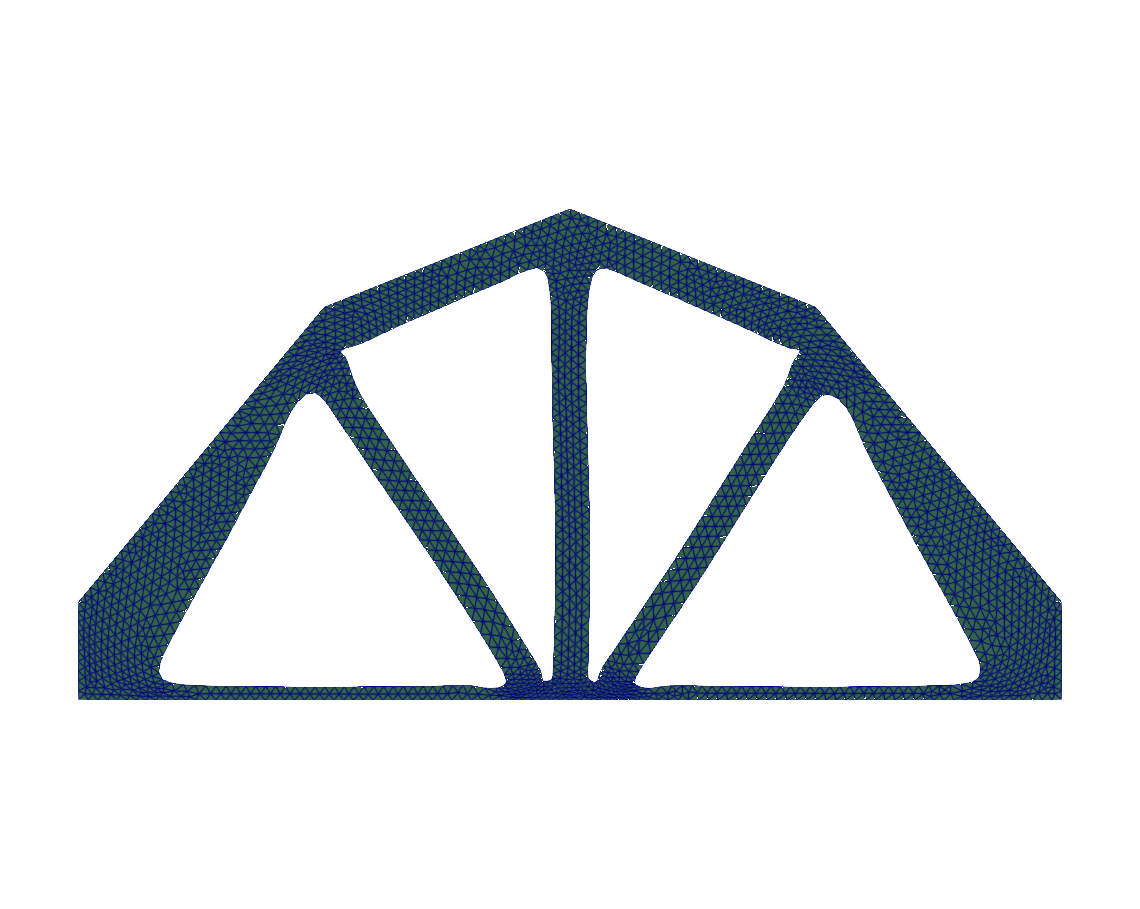}};
\node at (7.1,-1.2) {$u^{316}$};
\node at (7.1,-1.6) {$H^4$ metric};
\end{tikzpicture}
\caption{Solutions obtained with the different variants of the Riemannian descent method for the experiment described in \cref{subsec:compliance}}
\label{fig:exp_2_meshes}
\end{figure}

\Cref{tab:cost_experiment2} complements \cref{tab:exp_2_data} with the total computational cost of each complete run: the total wall-clock time, the time per iteration, and the total number of linear systems solved. As in Experiment~1, each iteration of the $H^s$ variants requires two linear solves (the state constraint of \eqref{prob:compliance} and the shape-gradient system~\eqref{eq:grad_system_reindexed}, the latter assembled as a single, monolithically coupled mixed-space system regardless of $s$), while the SP variant additionally solves a Poisson problem for the Lam\'e parameter $\mu$, for three solves per iteration in total. Unlike Experiment~1, this compliance-minimization problem does not involve a synthetic target state to be recovered, so the question of an inverse crime raised by the referee does not arise here: the objective in \eqref{prob:compliance} is evaluated directly from the computed state $\bm{y}$ and is never compared against an independently generated reference solution.

\pgfplotstableset{
  columns/iterations/.style={column type={c},
    column name={iterations $\bar k$},
  },
  columns/total_wall_clock_s/.style={column type={r},
    column name={total time}, fixed, precision=3,
    assign cell content/.code={%
      \pgfkeyssetvalue{/pgfplots/table/@cell content}%
      {\SI{##1}{\second}}%
    }},
  columns/s_per_iter/.style={column type={r},
    column name={time / iteration}, fixed, precision=4,
    assign cell content/.code={%
      \pgfkeyssetvalue{/pgfplots/table/@cell content}%
      {\SI{##1}{\second}}%
    }},
  columns/total_pde_solves/.style={column type={c},
    column name={total PDE solves}, /pgf/number format/1000 sep={},
  },
}
\begin{table}[ht]
\centering
\makebox[\linewidth][c]{%
\begin{filecontents*}{Data/com_total_cost_table.txt}
metric,iterations,total_wall_clock_s,s_per_iter,total_pde_solves,remeshing_count,final_step_size
SP,1000,125.354,0.1254,3000,3,3.333333e-01
H2,1000,174.846,0.1748,2000,2,5.000000e-01
H3,468,105.645,0.2257,936,3,3.333333e-01
H4,316,111.690,0.3534,632,3,3.333333e-01
\end{filecontents*}
\pgfplotstabletypeset[
  col sep=comma,
  columns={metric,iterations,total_wall_clock_s,s_per_iter,total_pde_solves},
  every head row/.style={before row=\toprule, after row=\midrule},
  every row no 0/.style={after row=\midrule},
  every last row/.style={after row=\bottomrule}]{Data/com_total_cost_table.txt}
}
\caption{Total computational cost of the complete optimization run for each variant of the method for the experiment described in \cref{subsec:compliance}, complementing \cref{tab:exp_2_data}. The $SP$, $H^3$ and $H^4$ variants required exactly three remeshing events and ended with stepsize $t=0.333$; while the $H^2$ variant required only two remeshing events with final stepsize $t=0.5$. }
\label{tab:cost_experiment2}
\end{table}

\begin{figure}[h]
\centering
\begin{tikzpicture}
\node[inner sep=0pt] (SP) at (-5.2,0)
    {\includegraphics[width=.245\textwidth]{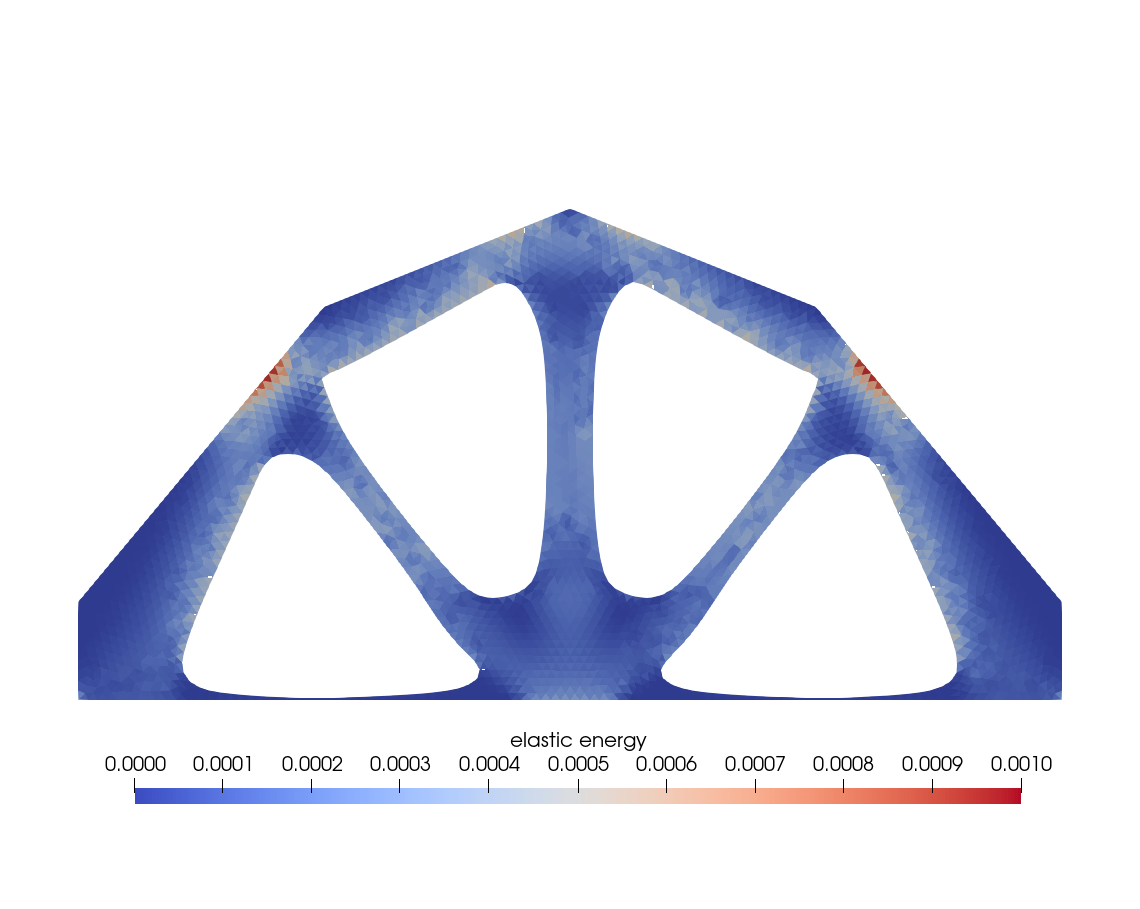}};
\node at (-5.2,-1.6) {$\mathcal{E}(u^{1000})$};
\node at (-5.2,-2) {SP metric};
\node[inner sep=0pt] (H1) at (-1.1,0)
    {\includegraphics[width=.245\textwidth]{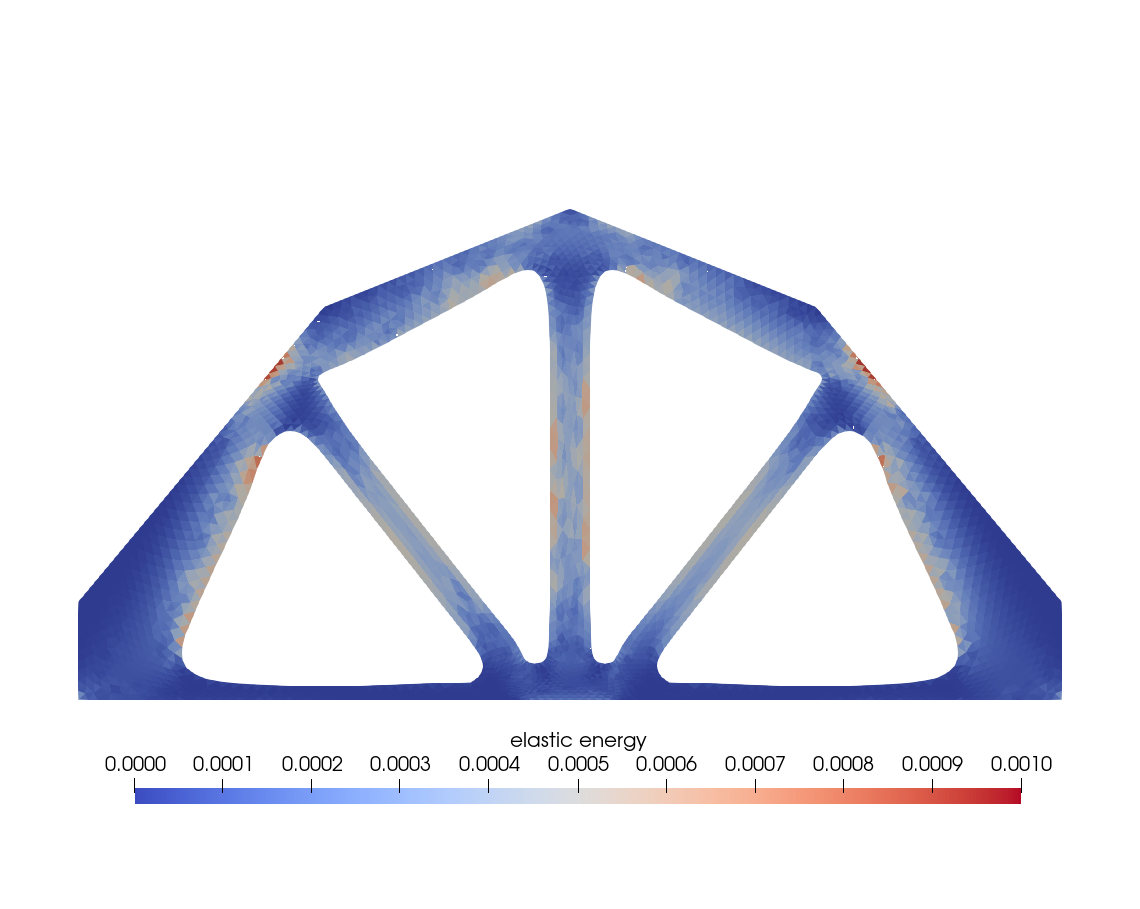}};
\node at (-1.1,-1.6) {$\mathcal{E}(u^{1000})$};
\node at (-1.1,-2) {$H^2$ metric};
\node[inner sep=0pt] (H2) at (3,0)
    {\includegraphics[width=.245\textwidth]{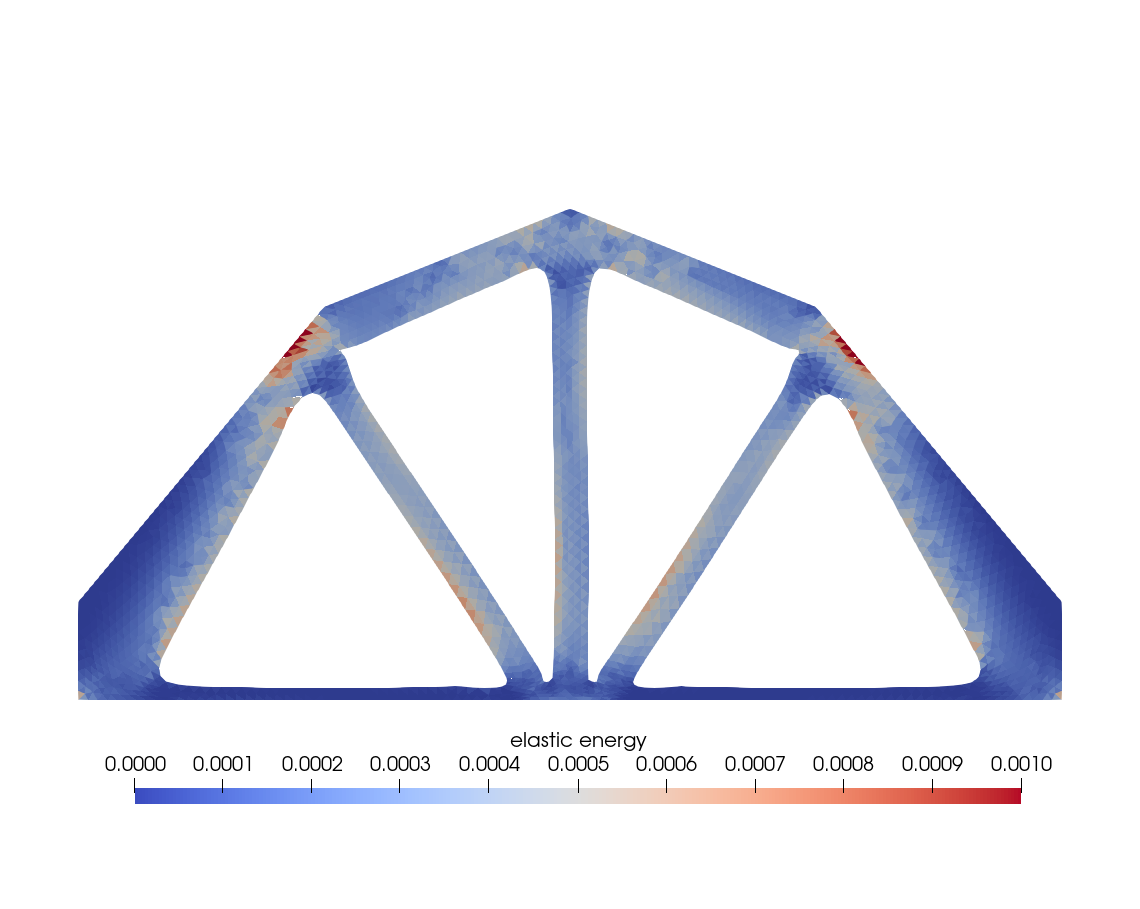}};
\node at (3,-1.6) {$\mathcal{E}(u^{468})$};
\node at (3,-2) {$H^3$ metric};
\node[inner sep=0pt] (H3) at (7.1,0)
    {\includegraphics[width=.245\textwidth]{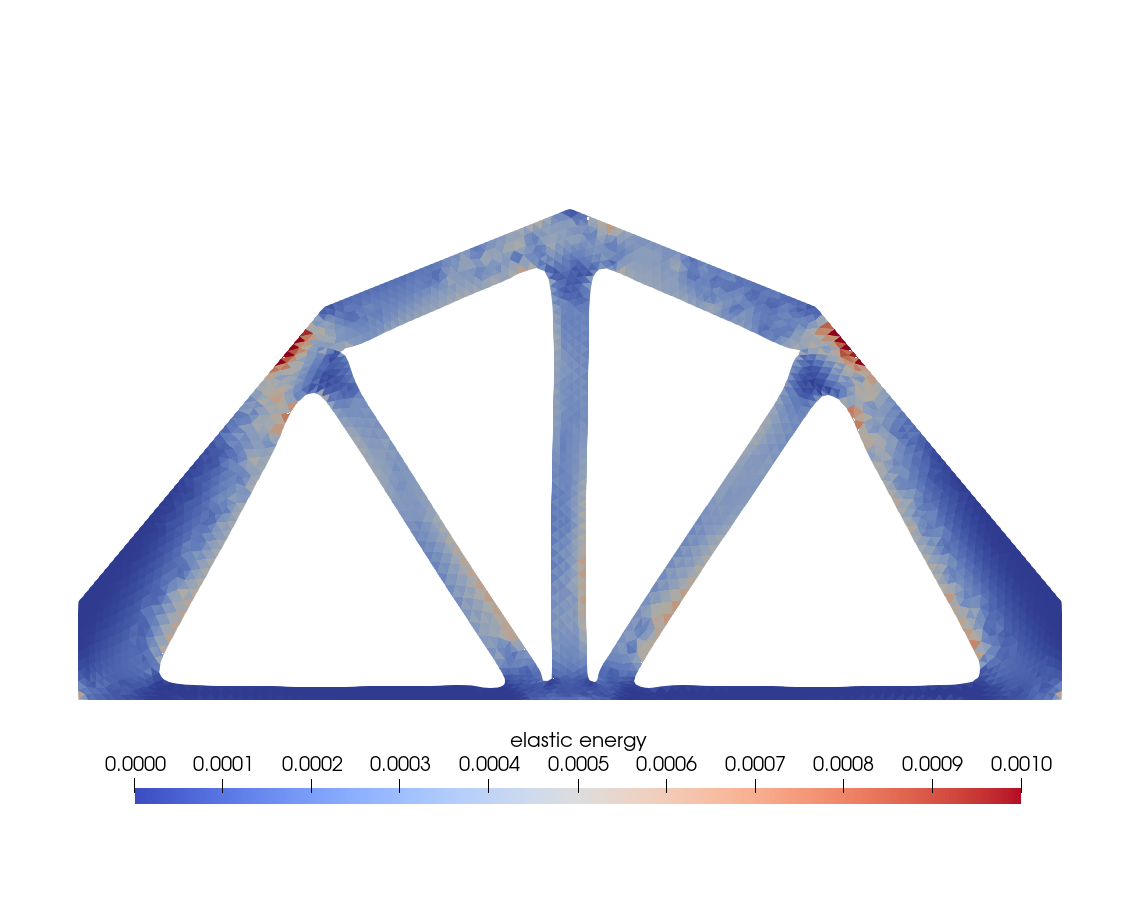}};
\node at (7.1,-1.6) {$\mathcal{E}(u^{316})$};
\node at (7.1,-2) {$H^4$ metric};
\end{tikzpicture}
\caption{Cell-wise elastic energy for each variant of the method for the experiment described in \cref{subsec:compliance}.}
\label{fig:exp_2_elastic_energy}
\end{figure}

Moreover, we compute the cell-wise contribution of the elastic energy, and it is depicted in~\cref{fig:exp_2_elastic_energy}. 

We finish this section by showing in \cref{fig:exp2_obj_qual} the evolution of the objective function, and quality of the mesh at each iteration of the algorithm. We have not plotted the values of the norm of the gradient since they all remain bounded away from zero, which indicates the ill-posedness of this problem since its optimal solution must present micro-structures, and tends to change the topology. 
However, using the objective function is valid since we aim to look for structures with similar values of compliance, and this is backed up by the cell-wise elastic energy shown in \cref{fig:exp_2_elastic_energy}, where we can see that all the structures have similar values. The scale used in this plot are all the same and range from $0$ to $10^{-3}$.  
We conclude that the best solution is again obtained with the $H^3$ metric since it attains the smallest value of the objective function while completing in less wall-clock time despite their higher per-iteration cost.
However, it must also be highlighted the good performance of the algorithm with $H^4$ metric, where the additional complexity of the computation of the shape gradient is compensated by the requirement of less iterations to satisfy the stopping criterion.  

\begin{figure}[h!]
\pgfplotsset{/pgfplots/group/.cd, horizontal sep=3cm}
\begin{center}
\begin{tikzpicture}
\pgfplotsset{compat=1.18, width=\textwidth}
\begin{groupplot}[group style={
        {group size=2 by 1}, xlabels at=edge bottom}, height=4cm, width=5cm,
        legend style={transpose legend, legend columns=4, draw=none }, 
        xlabel=iterations]
\nextgroupplot[title= $J(u^k)$, legend to name=grouplegend_comp]
\addplot[color=black, dashed] table [x=iter,y=objective, col sep=comma]{Data/com_history_SP.txt};\label{pgfplots:SP_comp}
\addlegendentry{SP metric}
\addplot[color=magenta, mark=*, mark size = 0.25pt] table [x=iter,y=objective, col sep=comma]{Data/com_history_H2.txt};\label{pgfplots:H2_comp}
\addlegendentry{$H^2$ metric}
\addplot[color=violet, mark=x, mark size = 0.25pt] table [x=iter,y=objective, col sep=comma]{Data/com_history_H3.txt};\label{pgfplots:H3_comp}
\addlegendentry{$H^3$ metric}
\addplot[color=orange, mark =star, mark size = 0.25pt] table [x=iter,y=objective, col sep=comma]{Data/com_history_H4.txt};\label{pgfplots:H4_comp}
\addlegendentry{$H^4$ metric}

\nextgroupplot[title = $\varphi(\holdall_h^k)$, ymode=log]
\addplot[color=black, dashed] table [x=iter,y=mesh_quality_min, col sep=comma]{Data/com_history_SP.txt};
\addplot[color=magenta, mark=*, mark size = 0.25pt] table [x=iter,y=mesh_quality_min, col sep=comma]{Data/com_history_H2.txt};
\addplot[color=violet, mark=x, mark size = 0.25pt] table [x=iter,y=mesh_quality_min, col sep=comma]{Data/com_history_H3.txt};
\addplot[color=orange, mark =star, mark size = 0.25pt] table [x=iter,y=mesh_quality_min, col sep=comma]{Data/com_history_H4.txt};
\end{groupplot}
\node[yshift=-35pt] at ($(group c1r1.south)!0.5!(group c2r1.south)$) {\ref{pgfplots:SP_comp} SP metric   \ref{pgfplots:H2_comp} $H^2$ metric  \ref{pgfplots:H3_comp} $H^3$ metric  \ref{pgfplots:H4_comp} $H^4$ metric};
\end{tikzpicture}
\caption{Objective function, and mesh quality behavior along the iterations for the different variants of the method for the experiment described in \cref{subsec:compliance}.}
\label{fig:exp2_obj_qual}
\end{center}
\end{figure}

\section{Conclusion}
\label{sec:conclusion}

In this paper, the properties of the Riemannian manifold $(\Diffc, H^s)$ are exploited while using a Riemannian version of a steepest descent method to applications in shape optimization from medical and engineering fields. Thanks to the geometry of $\Diffc$, we were able to prove an equivalence between the push-forward and the shape derivative.
Additionally, we leverage the well-established $H^s$ metrics on this shape space to devise different variants of the Riemannian steepest descent method. Each variant will compute shape gradients using a $H^s$ metric with $s = 1,2,3,4$.
Thus, shape gradients will be computed by solving a variational problem associated to a $2s$-th order differential operator. 
We, thus, propose a reformulation of this problem as a system of well-posed second-order problems, and characterize the field it produces in \cref{prop:split-equivalence}: its regularity, the Navier-type boundary conditions it inherits, and its relation to the $H^s$ inner product.
To study the performance of the different variants of the proposed method, we consider two problems:
the minimization of a tracking-type functional for interface identification inspired in electrical impedance tomography and a compliance minimization problem for a two-dimensional bridge. For both problems, we outlined their shape derivatives and referenced the appropriate results on their shape differentiability.

Our numerical investigations aim to compare the performance of the different variants of the Riemannian steepest descent method for the problems under consideration. In addition to comparing the results for the $H^s$ metric with different values of $s$, we also compare the results, with the well-known Steklov--Poincar\'e metric.

Overall, the variants of the method associated with the $H^s$ metrics (where $s=2,3,4$) demonstrated superior performance in terms of iteration count and objective function values, and, in the interface identification problem, in final mesh quality as well.

Two questions are left open by the present analysis. 
A rigorous account of why higher-order metrics preserve mesh quality is not attempted here; the regularity $\bm{V}\in H^{2s-1}(\holdall,\R^2)$ furnished by \cref{prop:split-equivalence} makes the deformation fields $C^1$ for $s\ge2$ and thereby suggests a mechanism, but turning this into a quantitative statement remains an interesting problem. 
Relatedly, that same regularity is what makes the update map admissible for $s\ge2$; for $s=1$ and for the SP metric, whose gradients lie only in $H^1(\holdall,\R^2)$, admissibility is enforced at the mesh level rather than following from the gradient, and a precise local result covering those cases is likewise open.

A further question is raised by the framework itself. 
The metric is defined on $\Xc(\R^2)$, whose elements have compact support, whereas the shape gradient produced by~\eqref{eq:grad_system_reindexed} vanishes on the fixed boundary $\Gamma$ but has a non-vanishing normal derivative there, by \cref{prop:split-equivalence}. Its extension by zero to $\R^2$ is therefore Lipschitz but not continuously differentiable across $\Gamma$, so that the computed deformation field does not lie in $\Xc(\R^2)$, and the identification of the update with an element of $\Diffc$ is not immediate. 
The discrepancy is confined to a neighbourhood of $\Gamma$ and does not affect the computations reported here, in which the deformation is applied at the moving boundary. 
Reconciling the two, i.e., whether by imposing further conditions on $\Gamma$, by working on a larger hold-all, or by weakening the ambient group is an interesting question that we leave open.

We recall that changes of topology are outside the scope of the method, the optimization proceeding by diffeomorphic deformation of a fixed initial configuration. 
Approaches that remove this restriction, by interleaving topological derivative steps with boundary evolution, are developed in~\cite{allaire2005structural}; combining such a mechanism with the Riemannian framework used here is an interesting question that we leave open.

In our future work, we intend to perform our computations using the exponential map instead of a retraction.
Note that the geodesic equation on $(\Diffc, H^s)$ was derived in \cite{Michor2007overview}, and that $(\Diffc,H^s)$ is geodesically and metrically complete for integer $s>2$, so that updates along geodesics would remain admissible for arbitrary step lengths.
Finally, from a computational point of view, an interesting direction for future research is to directly solve the polyharmonic equation using higher-order finite elements or mixed elements reformulations.

\section*{Acknowledgments}
This work has been partly supported by the International Centre for Mathematical Sciences (Edinburgh, UK) within the Research in Groups program. Additionally, the authors are indebted to Tim Suchan (Helmut-Schmidt-University/University of the Federal Armed Forces Hamburg, Germany) for contributing the remesh capability and supporting its incorporation into the existing program structure.


\printbibliography
\appendix

\section{Proof of \protect{\cref{prop:split-equivalence}}}
\label{sec:proof_prop_2}

Throughout,
\[
a(\bm X,\bm W)\coloneqq\int_\holdall \bm X\cdot\bm W\,\dd x
+A\int_\holdall\nabla\bm X:\nabla\bm W\,\dd x
\]
denotes the bilinear form appearing at each stage of~\eqref{eq:grad_system_reindexed}, and $\|\cdot\|_{*}$ the norm of $\big(H^1_{\Gamma,0}(\holdall,\R^2)\big)^{*}$.
We first record three auxiliary facts.

The proof of \cref{prop:split-equivalence} requires knowing that every element of the space $H^{1/2}(\Sigma,\R^2)$ is the trace of some field in $H^1_{\Gamma,0}(\holdall,\R^2)$, and not merely that such traces exist.
\Cref{lem:trace-factorisation} is not the classical trace theorem, which concerns the unconstrained space $H^1(\holdall,\R^2)$ and the whole of $\partial\holdall$; here the source is restricted by a condition on $\Gamma$ and the target is only the complementary portion.
Indeed, for a portion $\Gamma_1\subsetneq\partial\holdall$ meeting its complement, the traces of fields vanishing on $\partial\holdall\setminus\Gamma_1$ fill out the Lions--Magenes space $H^{1/2}_{00}(\Gamma_1,\R^2)$, which is strictly smaller than $H^{1/2}(\Gamma_1,\R^2)$ (cf.~\cite[Ch.~1, \S11]{lions2012non}). 
No such distinction arises here, because $\Gamma$ and $\Sigma$ are unions of connected components with $\ol\Gamma\cap\ol\Sigma=\emptyset$; there is no interface, extension by zero is bounded, and the two spaces coincide.

\begin{lemma}
\label{lem:trace-factorisation}
Let $\holdall$ be a bounded Lipschitz domain and let $\partial\holdall=\Gamma\sqcup\Sigma$ be a decomposition into unions of connected components.
Then $\gamma_\Sigma\colon H^1_{\Gamma,0}(\holdall,\R^2)\to H^{1/2}(\Sigma,\R^2)$ is bounded and surjective, with kernel $H^1_0(\holdall,\R^2)$.
Consequently it induces an isomorphism
$H^1_{\Gamma,0}(\holdall,\R^2)/H^1_0(\holdall,\R^2)\cong H^{1/2}(\Sigma,\R^2)$.
\end{lemma}

\begin{proof}
We start by proving the boundedness of the operator $\gamma_\Sigma$. 
Recall, the map $\gamma_\Sigma$ is the composition of $\gamma_0$ with restriction to
$\Sigma$.
The first factor satisfies $\|\gamma_0\bm W\|_{H^{1/2}(\partial\holdall,\R^2)}\le C_{\gamma}\|\bm W\|_{H^1(\holdall,\R^2)}$ for every $\bm W\in H^1(\holdall,\R^2)$, with $C_\gamma=C_\gamma(\holdall)$, by the trace theorem on a bounded Lipschitz domain.
For the second, $\ol\Gamma\cap\ol\Sigma=\emptyset$ gives the direct sum decomposition $H^{1/2}(\partial\holdall,\R^2)=H^{1/2}(\Gamma,\R^2)\oplus H^{1/2}(\Sigma,\R^2)$, in which extension by zero from either portion is bounded and restriction to $\Sigma$ is the associated projection, of norm $C_\Sigma$.
Combining the two estimates gives
\[
\|\gamma_\Sigma\bm W\|_{H^{1/2}(\Sigma,\R^2)}
\le C_{\Sigma}C_{\gamma}\,\|\bm W\|_{H^1(\holdall,\R^2)}
\qquad\text{for all }\bm W\in H^1_{\Gamma,0}(\holdall,\R^2),
\]
the restriction of a bounded operator to a closed subspace being bounded with the same constant.

Now, we focus on characterizing the kernel of $\gamma_\Sigma$. 
Let $\bm W $  be an element of the kernel of $\gamma_\Sigma$, i.e., $\bm W\in H^1_{\Gamma,0}(\holdall,\R^2)$ and $\gamma_\Sigma\bm W=0$. It follows automatically that $\gamma_0\bm W=0$ on $\Gamma$ by the definition of the space $H^1_{\Gamma,0}(\holdall, \R^2)$ and on $\Sigma$ by hypothesis, hence on all of $\partial\holdall$, so $\bm W\in H^1_0(\holdall,\R^2)$; the converse inclusion is clear.
Finally, we study its surjectivity, let $\bm g\in H^{1/2}(\Sigma,\R^2)$ and let $\tilde{\bm g}$ denote its extension by zero to $\partial\holdall$.
Because $\Gamma$ and $\Sigma$ are unions of connected components, they are relatively open and closed in $\partial\holdall$ and $\ol\Gamma\cap\ol\Sigma=\emptyset$; no matching condition is therefore imposed across a common interface, and $\tilde{\bm g}\in H^{1/2}(\partial\holdall,\R^2)$.
Applying a bounded right inverse $E$ of $\gamma_0$ gives $E\tilde{\bm g}\in H^1(\holdall,\R^2)$ with $\gamma_0 E\tilde{\bm g}=\tilde{\bm g}$, so that $E\tilde{\bm g}\in H^1_{\Gamma,0}(\holdall,\R^2)$ and $\gamma_\Sigma E\tilde{\bm g}=\bm g$.
Since $\gamma_\Sigma$ vanishes on $H^1_0(\holdall,\R^2)$, it factors through the quotient, inducing $\bar\gamma_\Sigma\colon H^1_{\Gamma,0}(\holdall,\R^2)/H^1_0(\holdall,\R^2)\to H^{1/2}(\Sigma,\R^2)$ with $\bar\gamma_\Sigma[\bm W]=\gamma_\Sigma\bm W$.
This map is continuous for the quotient norm, injective by the computation of the kernel, and surjective by the preceding step, hence an isomorphism of Banach spaces by the open mapping theorem.
\end{proof}

As already mentioned, in the next lemma we present the well-possednes of \cref{def:gen-conormal}. 
\begin{lemma}
\label{lem:gen-conormal}
Let $F$ and $\bm Z$ be as stated in \cref{def:gen-conormal}; then, there exists exactly one $\partial_n^{F}\bm Z\in H^{-1/2}(\Sigma,\R^2)$ which satisfies~\eqref{eq:gen-conormal}.
If moreover $\bm Z\in H^2(\holdall,\R^2)$ and $F\in L^2(\holdall,\R^2)$, then $\partial_n^{F}\bm Z=\partial_n\bm Z$.
\end{lemma}

\begin{proof}
For $\bm W\in H^1_{\Gamma,0}(\holdall,\R^2)$ set
\[
\Lambda[\bm W]\coloneqq
\int_\holdall \bm Z\cdot\bm W\,\dd x
+A\int_\holdall\nabla\bm Z:\nabla\bm W\,\dd x
-F[\bm W],
\]
so that~\eqref{eq:gen-conormal} reads $A\langle\partial_n^{F}\bm Z,\gamma_\Sigma\bm W\rangle=\Lambda[\bm W]$.
It is clear that $\Lambda$ is a bounded linear functional on $H^1_{\Gamma,0}(\holdall,\R^2)$, since $\bm Z\in H^1(\holdall,\R^2)$ and
$F\in\big(H^1_{\Gamma,0}(\holdall,\R^2)\big)^{*}$.
For $\bm\phi\in C_0^\infty(\holdall,\R^2)$ one has $\Lambda[\bm\phi]=0$, this being the statement $(\id-A\Delta)\bm Z=F$ in $\mathcal D'(\holdall,\R^2)$; since $C_0^\infty(\holdall,\R^2)$ is dense in $H^1_0(\holdall,\R^2)$ and $\Lambda$ is continuous, $\Lambda$ vanishes on $H^1_0(\holdall,\R^2)$.
It therefore factors through the quotient, inducing a bounded functional
$\bar\Lambda$ on $H^1_{\Gamma,0}(\holdall,\R^2)/H^1_0(\holdall,\R^2)$ with
$\bar\Lambda[\bm W]=\Lambda[\bm W]$.
By \cref{lem:trace-factorisation} the map $\bar\gamma_\Sigma$ is an isomorphism of
that quotient onto $H^{1/2}(\Sigma,\R^2)$, so
$\left(\nicefrac{1}{A}\right)\,\left(\bar\Lambda\circ\bar\gamma_\Sigma^{-1}\right)$ is a bounded functional on
$H^{1/2}(\Sigma,\R^2)$, that is, an element of $H^{-1/2}(\Sigma,\R^2)$.
Taking $\partial_n^{F}\bm Z$ to be that element gives~\eqref{eq:gen-conormal}. 

Uniqueness follows from the surjectivity of $\gamma_\Sigma$. 
Indeed, if two elements of $H^{-1/2}(\Sigma,\R^2)$ both satisfy~\eqref{eq:gen-conormal}, their difference $\bm\mu$ satisfies $\langle\bm\mu,\gamma_\Sigma\bm W\rangle=0$ for every $\bm W\in H^1_{\Gamma,0}(\holdall,\R^2)$, hence $\langle\bm\mu,\bm g\rangle=0$ for every $\bm g\in H^{1/2}(\Sigma,\R^2)$, and so
$\bm\mu=0$.

If in addition $\bm Z\in H^2(\holdall,\R^2)$ and $F\in L^2(\holdall,\R^2)$, then Green's formula gives
$\Lambda[\bm W]=A\int_{\partial\holdall}\partial_n\bm Z\cdot\gamma_0\bm W\,\ds$,
and the contribution of $\Gamma$ vanishes because $\gamma_0\bm W=0$ there; the two
notions therefore coincide.
\end{proof}

The regularity of the fields produced by~\eqref{eq:grad_system_reindexed} rests on a single analytical ingredient, applied $s-1$ times in the proof below: a stage with datum in $H^m$ has its solution in $H^{m+2}$, up to the boundary.
This is the classical shift theorem, but not in the form usually quoted, since textbook statements concern either the Dirichlet or the Neumann problem on the whole of $\partial\holdall$, whereas each stage carries an essential condition on $\Gamma$ and a natural one on $\Sigma$. 
It is worth highlight that what makes the mixed problem no harder than either pure one is the decomposition of $\partial\holdall$ itself. 
In our case, since $\Gamma$ and $\Sigma$ are unions of connected components, the type of the boundary condition is constant on each component and no point carries a change of type.
The stage is then a regular elliptic boundary value problem, and the shift theorem applies directly.
It is also the only point at which the regularity of $\partial\holdall$ enters the proof of \cref{prop:split-equivalence}.

\begin{lemma}
\label{lem:shift}
Let $m\ge0$ be an integer and let every connected component of $\partial\holdall$ be of class $C^{m+1,1}$.
Let $\bm f\in H^{m}(\holdall,\R^2)$ and let $\bm Z\in H^1_{\Gamma,0}(\holdall,\R^2)$ satisfy
$a(\bm Z,\bm W)=\int_\holdall\bm f\cdot\bm W\,\dd x$ for all $\bm W\in H^1_{\Gamma,0}(\holdall,\R^2)$.
Then $\bm Z\in H^{m+2}(\holdall,\R^2)$ and
\[
\|\bm Z\|_{H^{m+2}(\holdall,\R^2)}\le C(m,A,\holdall,\Gamma,\Sigma)\,\|\bm f\|_{H^{m}(\holdall,\R^2)} .
\]
\end{lemma}

\begin{proof}
The operator $\id-A\Delta$ is uniformly elliptic with constant coefficients and diagonal in the Cartesian components, so it suffices to treat each scalar component.
By the preceding discussion the stage is a regular elliptic boundary value problem, that is, the Dirichlet operator on $\Gamma$ together with the Neumann operator on $\Sigma$ is a normal system covering $\id-A\Delta$, and the shift theorem applies with the stated boundary regularity; cf.~\cite[Th.~2.5.1.1]{grisvard2011elliptic}.
No compatibility condition is imposed on $\bm f$. 
Unlike the pure Neumann problem for $-\Delta$, the form $a$ is coercive because of its zeroth-order term, so the problem is uniquely solvable for every datum.
\end{proof}

\begin{proof}[Proof of \cref{prop:split-equivalence}]

We begin with unique solvability, which is a routine application of the Lax--Milgram theorem to each stage $j$.
The form $a$ is bounded on $H^1(\holdall,\R^2)$, with constant $\max(1,A)$, and satisfies $a(\bm X,\bm X)\ge\min(1,A)\|\bm X\|^2_{H^1(\holdall,\R^2)}$; the zeroth-order term supplies the coercivity (no Poincar\'e inequality is needed) and the estimate holds on the closed subspace $H^1_{\Gamma,0}(\holdall,\R^2)$ as well.
Since~\eqref{eq:grad_system_reindexed} is triangular, we solve it in decreasing order of $j$.
For $j=s-1$, the right-hand side is bounded on $H^1_{\Gamma,0}(\holdall,\R^2)$ by \cref{assmp:dj-bounded}, and at each subsequent stage it is the field $\bm X_{j+1}\in L^2(\holdall,\R^2)$ just obtained; in both cases Lax--Milgram gives a unique solution, with
\begin{align*}
\|\bm X_{s-1}\|_{H^1}&\le\min(1,A)^{-1}\|DJ(u)\|_{*},\\
\|\bm X_j\|_{H^1}&\le\min(1,A)^{-1}\|\bm X_{j+1}\|_{L^2}
 \quad \text{ for } j\le s-2.
\end{align*}
Iterating from $j=s-1$ to $j=0$ determines the tuple, uniquely at every stage.
Note that no regularity of $\partial\holdall$ has been used here, nor anything about how $\partial\holdall$ splits into $\Gamma$ and $\Sigma$.

\medskip

We now show that the fields just constructed are considerably more regular than the formulation suggests, and that each satisfies its equation in the strong sense.
This is what gives meaning to the expression $(\id-A\Delta)^j\bm V$ appearing in the statement, and what allows the integrations by parts used later to recover the boundary conditions.
Precisely, we claim that for $k=0,\dots,s-1$, 
\begin{equation}\label{eq:strong-form}
\bm X_{s-1-k}\in H^{2k+1}(\holdall,\R^2)
\qquad\text{and}\qquad
(\id-A\Delta)\bm X_{s-1-k}=f_k ,
\end{equation}
where $f_0\coloneqq DJ(u)$ and $f_k\coloneqq\bm X_{s-k}$ for $k\ge1$.
For $k=0$ the identity holds in $\big(C_0^\infty(\holdall,\R^2)\big)^{*}$; from $k\ge1$ onwards the right-hand side is the square-integrable field $\bm X_{s-k}$, both sides are functions, and it holds almost everywhere on $\holdall$. 
We require the strong forms, rather than the weak formulations alone, since they can be composed with one another, which is what yields the identity $\bm X_j=(\id-A\Delta)^j\bm V$ below; the weak formulations cannot be composed, as each one will be tested against a different function.

We argue by induction on $k\coloneqq s-1-j$, which counts the stages of~\eqref{eq:grad_system_reindexed} in the order they are solved, from $k=0$ at $j=s-1$ to $k=s-1$ at $j=0$.

For $k=0$, testing the stage $j=s-1$ against $\bm\phi\in C_0^\infty(\holdall,\R^2)$ gives
\[
\int_\holdall \bm X_{s-1}\cdot\bm\phi\,\dd x
+A\int_\holdall\nabla\bm X_{s-1}:\nabla\bm\phi\,\dd x
=DJ(u)[\bm\phi],
\]
which is $(\id-A\Delta)\bm X_{s-1}=DJ(u)$ in $\big(C_0^\infty(\holdall,\R^2)\big)^{*}$, and $\bm X_{s-1}\in H^1(\holdall,\R^2)$ as already shown.
It is worth highlighting this is the regularity the datum $DJ(u)$ permits, and the iterated elliptic regularity argument improves on it from here.

Suppose now that~\eqref{eq:strong-form} holds for some $k\le s-2$; we prove it for $k+1$.
The field $\bm X_{s-k-2}$ solves the stage $j=s-k-2$ in $H^1_{\Gamma,0}(\holdall,\R^2)$, with datum $f_{k+1}=\bm X_{s-k-1}\in H^{2k+1}(\holdall,\R^2)$ by the induction hypothesis.
\cref{lem:shift} therefore applies with $m=2k+1$, its boundary hypothesis $C^{2k+2,1}$ following from $C^{2s-2,1}$ since $2k+2\le 2s-2$ for $k\le s-2$.
It gives $\bm X_{s-k-2}\in H^{2k+3}(\holdall,\R^2)$; both sides of the equation are then square integrable and the identity holds almost everywhere, which is~\eqref{eq:strong-form} at $k+1$.

Written in terms of $j$, the claim reads $\bm X_j\in H^{2(s-1-j)+1}(\holdall,\R^2)$ for $j=0,\dots,s-1$, and in particular $\bm V\in H^{2s-1}(\holdall,\R^2)\hookrightarrow H^{s}(\holdall,\R^2)$.
Chaining the Lax--Milgram bound with the $s-1$ applications of \cref{lem:shift} gives~\eqref{eq:prop-apriori}, with $C=\min(1,A)^{-1}\prod_{k=1}^{s-1}C(2k-1,A,\holdall,\Gamma,\Sigma)$.

The regularity and the strong forms just established are exactly what is needed to complete part~(i).  
The first makes each $(\id-A\Delta)\bm X_{j-1}$ a well-defined Sobolev field, and the second identifies it as $\bm X_j$.
Composing these identities gives, for $j=1,\dots,s-1$,
\[
\bm X_j=(\id-A\Delta)\bm X_{j-1}
=(\id-A\Delta)^2\bm X_{j-2}=\cdots
=(\id-A\Delta)^j\bm X_0=(\id-A\Delta)^j\bm V ,
\]
an identity in $H^{2(s-1-j)+1}(\holdall,\R^2)$.
Every substitution is justified as follows: for $j\le s-1$ one has $\bm X_{j-1}\in H^{2(s-j)+1}(\holdall,\R^2)\subset H^{3}(\holdall,\R^2)$, so $(\id-A\Delta)$ carries it into $H^{1}(\holdall,\R^2)$, the Sobolev exponent is reduced by exactly two at each application, which is the exponent stated.

Now we focus on the boundary conditions.
Since $2(s-1-j)+1\ge1$ for every $j\le s-1$, both sides of the identity belong to $H^1(\holdall,\R^2)$ and $\gamma_0$ may be applied.
Each $\bm X_j$ lies in $H^1_{\Gamma,0}(\holdall,\R^2)$ by construction, so $\gamma_\Gamma\bm X_j=0$, which by the identity $\bm X_j=(\id-A\Delta)^j\bm V$ eliminates the auxiliary fields and gives~\eqref{eq:bc-dirichlet}.
On $\Sigma$ the traces are in general non-zero, the conditions imposed are natural rather than essential; they are recovered below from the variational identity itself, in the sense of \cref{def:gen-conormal} where the regularity requires it.

For~\eqref{eq:bc-neumann}, let $j\le s-2$, so that $\bm X_j\in H^{2(s-1-j)+1}(\holdall,\R^2)\subset H^{3}(\holdall,\R^2)$ and $\partial_n\bm X_j\in H^{1/2}(\Sigma,\R^2)$ is classically defined.
Applying Green's formula to the stage $j$ of~\eqref{eq:grad_system_reindexed} gives, for every $\bm W\in H^1_{\Gamma,0}(\holdall,\R^2)$,
\[
\begin{split}
\int_\holdall\big((\id-A\Delta)\bm X_j\big)\cdot\bm W\,\dd x
+A\int_{\Gamma}\partial_n\bm X_j\cdot\gamma_\Gamma\bm W\,\ds
+A\int_{\Sigma}\partial_n\bm X_j\cdot\gamma_\Sigma\bm W\,\ds
\\
=\int_\holdall\bm X_{j+1}\cdot\bm W\,\dd x ,
\end{split}
\]
where the boundary integral over $\partial\holdall$ has been split, $\Gamma$ and $\Sigma$ being unions of connected components.
The integral over $\Gamma$ vanishes because $\gamma_\Gamma\bm W=0$, and $(\id-A\Delta)\bm X_j=\bm X_{j+1}$ almost everywhere by~\eqref{eq:strong-form}, so the integral over $\Sigma$ vanishes for every such $\bm W$.
By \cref{lem:trace-factorisation} every element of $H^{1/2}(\Sigma,\R^2)$ is of the form $\gamma_\Sigma\bm W$, so $\partial_n\bm X_j$ integrates to zero against every element of that space, in particular against itself.
Hence $\partial_n\bm X_j=0$ on $\Sigma$, and since $\bm X_j=(\id-A\Delta)^j\bm V$, this is~\eqref{eq:bc-neumann}.

For~\eqref{eq:bc-neumann-top}, the field $\bm X_{s-1}$ is only of class $H^1$ and its conormal derivative is not classically defined.
By~\eqref{eq:strong-form} with $k=0$ we may apply \cref{def:gen-conormal} with $F=DJ(u)$ and $\bm Z=\bm X_{s-1}$; the right-hand side of~\eqref{eq:gen-conormal} is then, for every $\bm W\in H^1_{\Gamma,0}(\holdall,\R^2)$,
$a(\bm X_{s-1},\bm W)-DJ(u)[\bm W]=0$ by the stage $j=s-1$ itself.
Hence $\partial_n^{DJ(u)}\bm X_{s-1}=0$ in $H^{-1/2}(\Sigma,\R^2)$,  which by the identity $\bm X_{s-1}=(\id-A\Delta)^{s-1}\bm V$ eliminates the auxiliary field and gives~\eqref{eq:bc-neumann-top}.
This proves~(i).

We turn to~(ii).
The forward implication is immediate from~(i). 
Indeed, testing the stage $j=s-1$ against $\bm\phi\in C_0^\infty(\holdall,\R^2)$ and using $\bm X_{s-1}=(\id-A\Delta)^{s-1}\bm V$ gives $(\id-A\Delta)^s\bm V=DJ(u)$ in $\big(C_0^\infty(\holdall,\R^2)\big)^{*}$, whose boundary conditions are~\eqref{eq:bc-dirichlet}--\eqref{eq:bc-neumann-top}.

Conversely, let $\bm V\in H^{2s-1}(\holdall,\R^2)$ satisfy that equation together with~\eqref{eq:bc-dirichlet}--\eqref{eq:bc-neumann-top}, and set $\bm X_j\coloneqq(\id-A\Delta)^j\bm V$ for $j=0,\dots,s-1$.
Then it holds that $\bm X_j\in H^{2(s-1-j)+1}(\holdall,\R^2)\subset H^1(\holdall,\R^2)$ and $\gamma_\Gamma\bm X_j=0$ by~\eqref{eq:bc-dirichlet}, so $\bm X_j\in H^1_{\Gamma,0}(\holdall,\R^2)$.
We check that this tuple solves~\eqref{eq:grad_system_reindexed}.
For $j\le s-2$ one has $\bm X_j\in H^{3}(\holdall,\R^2)$ and $(\id-A\Delta)\bm X_j=\bm X_{j+1}$ pointwise. Green's formula gives, for every $\bm W\in H^1_{\Gamma,0}(\holdall,\R^2)$,
\[
a(\bm X_j,\bm W)=\int_\holdall\bm X_{j+1}\cdot\bm W\,\dd x
+A\int_{\Gamma}\partial_n\bm X_j\cdot\gamma_\Gamma\bm W\,\ds
+A\int_{\Sigma}\partial_n\bm X_j\cdot\gamma_\Sigma\bm W\,\ds .
\]
The integral over $\Gamma$ vanishes because $\gamma_\Gamma\bm W=0$, and the integral over $\Sigma$ by~\eqref{eq:bc-neumann}, so $a(\bm X_j,\bm W)=\int_\holdall\bm X_{j+1}\cdot\bm W\,\dd x$, which is the stage $j$ of~\eqref{eq:grad_system_reindexed}.

At $j=s-1$ the field $\bm X_{s-1}$ is of class $H^1$ only, and~\eqref{eq:bc-neumann-top} takes over: the equation puts $\bm X_{s-1}$ in the setting of \cref{def:gen-conormal} with $F=DJ(u)$, and the vanishing of the generalised conormal derivative says that the right-hand side of~\eqref{eq:gen-conormal} is zero, that is, $a(\bm X_{s-1},\bm W)=DJ(u)[\bm W]$ for every $\bm W\in H^1_{\Gamma,0}(\holdall,\R^2)$.
The tuple therefore solves~\eqref{eq:grad_system_reindexed}, and by uniqueness it is the one constructed above; in particular $\bm V=\bm X_0$, which is also the uniqueness asserted in~(ii).

It remains to prove~(iii).
Let $\bm W\in C_0^\infty(\holdall,\R^2)$.
Every integral below is supported in $\operatorname{supp}\bm W\Subset\holdall$, so no extension of $\bm V$ outside $\holdall$ is needed.
Moving $i$ derivatives in $x$ and $j$ derivatives in $y$ off $\bm V$ and onto $\bm W$, which is the definition of the distributional derivative and produces no boundary terms, turns the left-hand side of the identity in~(iii), that is the form of~\eqref{eq:metric_whole_space} evaluated at $(\bm V,\bm W)$, into
\[
\sum_{\substack{i,j\ge0\\ i+j\le s}}
\frac{(-A)^{i+j}s!}{i!\,j!\,(s-i-j)!}
\big\langle\partial_x^{2i}\partial_y^{2j}\bm V,\bm W\big\rangle ,
\]
the sign $(-1)^{i+j}$ from the $i+j$ integrations by parts cancelling against the one in $(-A)^{i+j}$, the pairings on the right being those of $\big(C_0^\infty(\holdall,\R^2)\big)^{*}$ with $C_0^\infty(\holdall,\R^2)$.
By the multinomial expansion of $(\id-A\Delta)^s$ given by $\sum_{i+j\le s}\frac{(-A)^{i+j}s!}{i!\,j!\,(s-i-j)!}\, \partial_x^{2i}\partial_y^{2j}$, the sum above is $\big\langle(\id-A\Delta)^s\bm V,\bm W\big\rangle$, a duality pairing rather than an integral, since $\bm V\in H^{2s-1}(\holdall,\R^2)$ gives $(\id-A\Delta)^s\bm V\in H^{-1}(\holdall,\R^2)$.
By~(ii) it equals $DJ(u)[\bm W]$, which proves~(iii). 
\end{proof}

\end{document}